\documentclass{amsart}
\usepackage[utf8]{inputenc}

\usepackage[T1]{fontenc}
\usepackage{lmodern}
\usepackage[english]{babel}
\usepackage{amscd,amsfonts,amsmath,amssymb,amsthm,mathrsfs}
\usepackage{bbm}

\usepackage{tikz}
\usetikzlibrary{arrows.meta}

\usepackage{geometry}
\usepackage{tikz-cd}
\usepackage[matrix,arrow,curve]{xy}
\usepackage{hyperref}

\hypersetup{
  pdftitle={Du Bois cohomology of noncommutative spaces and the Hodge numbers invariance conjecture for singular varieties},
  pdfauthor={Andrei Konovalov}
}

\usepackage{xcolor}

\theoremstyle{plain}
\newtheorem{prop}{Proposition}[section]
\newtheorem{theorem}[prop]{Theorem}
\newtheorem{corollary}[prop]{Corollary}
\newtheorem{lemma}[prop]{Lemma}
\newtheorem{obs}[prop]{Observation}
\theoremstyle{definition}
\newtheorem{defi}[prop]{Definition}
\newtheorem{ques}[prop]{Question}
\newtheorem{remark}[prop]{Remark}
\newtheorem{example}[prop]{Example}

\newcommand{\lra}{\longrightarrow}
\newcommand{\ra}{\rightarrow}
\newcommand{\cA}{{\mathcal A}}
\newcommand{\cB}{{\mathcal B}}
\newcommand{\cC}{{\mathcal C}}
\newcommand{\cE}{{\mathcal E}}
\newcommand{\bbA}{\mathbb{A}}
\newcommand{\bbC}{\mathbb{C}}
\newcommand{\bbF}{\mathbb{F}}
\newcommand{\bbN}{\mathbb{N}}
\newcommand{\bbQ}{\mathbb{Q}}
\newcommand{\bbS}{\mathbb{S}}
\newcommand{\bbZ}{\mathbb{Z}}
\newcommand{\Catperf}{\mathrm{Cat}^{\mathrm{perf}}}
\newcommand{\Perf}{\mathrm{Perf}}
\newcommand{\Motloc}{\mathrm{Mot}^{\mathrm{loc}}}
\newcommand{\MotlocA}{\mathrm{Mot}^{\mathrm{loc},\mathbb A^1}}
\newcommand{\Motloccyc}{\mathrm{Mot}^{\mathrm{loc},\mathrm{cyc}}}
\newcommand{\Motlocwg}{\mathrm{Mot}^{\mathrm{loc},\mathrm{wg}}}
\newcommand{\Uloc}{\mathcal U_{\mathrm{loc}}}
\newcommand{\TC}{\mathrm{TC}}
\newcommand{\HH}{\mathrm{HH}}
\newcommand{\HC}{\mathrm{HC}}
\newcommand{\HP}{\mathrm{HP}}

\newcommand{\fib}{\mathrm{fib}}
\newcommand{\cof}{\mathrm{cof}}

\newcommand{\coker}{\mathrm{coker}}

\newcommand{\bbP}{\mathbb{P}}

\newcommand{\K}{\mathrm{K}}
\newcommand{\KH}{\mathrm{KH}}
\newcommand{\DK}{\mathrm{DK}}
\newcommand{\KU}{\mathrm{KU}}
\newcommand{\THH}{\mathrm{THH}}
\newcommand{\TR}{\mathrm{TR}}
\newcommand{\DHH}{\mathrm{DHH}}
\newcommand{\DHC}{\mathrm{DHC}}

\newcommand{\dk}{\mathrm{dk}}
\newcommand{\kh}{\mathrm{kh}}
\newcommand{\tc}{\mathrm{tc}}

\newcommand{\kinv}{\mathrm{k}^{\mathrm{inv}}}

\newcommand{\Coh}{\mathrm{Coh}}
\newcommand{\QCoh}{\mathrm{QCoh}}
\newcommand{\Sg}{\mathrm{Sg}}
\newcommand{\Mod}{\mathrm{Mod}}
\newcommand{\Sp}{\mathrm{Sp}}
\newcommand{\CycSp}{\mathrm{CycSp}}
\newcommand{\Sch}{\mathrm{Sch}}
\newcommand{\Fun}{\mathrm{Fun}}
\newcommand{\Ind}{\mathrm{Ind}}
\newcommand{\Ho}{\mathrm{Ho}}

\DeclareMathOperator{\Hom}{Hom}
\DeclareMathOperator{\End}{End}
\DeclareMathOperator{\Res}{Res}
\DeclareMathOperator{\Spec}{Spec}
\DeclareMathOperator{\Bl}{Bl}
\DeclareMathOperator{\Gr}{Gr}
\DeclareMathOperator{\gr}{gr}
\DeclareMathOperator{\Fil}{Fil}
\DeclareMathOperator{\id}{id}
\DeclareMathOperator{\im}{im}
\DeclareMathOperator{\ch}{ch}
\DeclareMathOperator{\td}{td}
\DeclareMathOperator{\sig}{sig}
\DeclareMathOperator{\Cl}{Cl}
\DeclareMathOperator{\Pic}{Pic}

\newcommand{\DB}{\mathrm{DB}}
\newcommand{\NS}{\mathrm{NS}}
\newcommand{\gen}{\mathrm{gen}}
\newcommand{\cl}{\mathrm{cl}}
\newcommand{\der}{\mathrm{der}}
\newcommand{\an}{\mathrm{an}}
\newcommand{\cdh}{\mathrm{cdh}}
\newcommand{\pcdh}{\mathrm{pcdh}}
\newcommand{\dR}{\mathrm{dR}}

\DeclareMathOperator{\Sing}{Sing}

\newcommand{\sft}{\mathrm{sft}}

\newcommand{\E}{\mathrm{E}}

\author{Andrei Konovalov}
\title[Du Bois invariants and mixed Hodge numbers]{Du Bois cohomology of noncommutative spaces and the Hodge numbers invariance conjecture for singular varieties}
\email{kaa.akonovalov@gmail.com}

\begin{document}

\begin{abstract}
We extend the definition of cdh-sheafified Hochschild homology from varieties to stable categories linear over a field of characteristic zero, thereby proving its Morita-invariance. The dimensions of its homology groups give meaningful numerical invariants of a category; we combine them with invariants coming from topological K-theory and Poincar\'e duality to prove equality of mixed Hodge numbers for singular complex projective Fourier--Mukai partners in several classes in dimensions two, three, and four, including normal surfaces and threefolds with rational singularities and equal irregularity.

Beyond these applications, we construct Du Bois versions of topological cyclic homology of stable categories together with appropriate multiplicative structures, prove a Geller--Weibel-type result for the variants of negative cyclic homology we introduce, and give categorical characterizations of cohomological purity of proper varieties and the $\bbQ$-factoriality defect of complex projective threefolds with rational singularities.
\end{abstract}
\maketitle

\tableofcontents

\section{Introduction}

The perfect derived category of an algebraic variety is known to retain a great deal of information about the variety.
In particular, many invariants of a variety, such as its algebraic K-theory and Hochschild homology, can be defined in terms of its perfect derived category alone. Still, some information is lost:
there are many examples of nonisomorphic varieties with equivalent perfect derived categories. It is therefore natural to ask which geometric information is retained by the perfect derived category.

The Hodge numbers invariance conjecture, motivated by the homological mirror symmetry program suggested in the 1994 talk of M.~Kontsevich \cite{Kon95}, predicts that smooth projective complex varieties with equivalent derived categories have the same Hodge numbers. The conjecture is known in dimensions up to four; see \cite{PS11} for threefolds, \cite{Ant26} for fourfolds, and \cite{Abu26} for some further progress. For singular proper varieties, Betti cohomology groups carry mixed Hodge structures, numerically described by mixed Hodge numbers
\[
h_n^{a,b}(X):=
\dim_{\bbC}\Gr_F^a\Gr_{a+b}^W H^n(X^{\an},\bbC).
\]
We ask to what extent these numbers are determined by the perfect derived category $\Perf(X)$.

In the smooth proper case, Hochschild homology supplies an important part of the answer: the Hochschild--Kostant--Rosenberg decomposition expresses its dimensions as column sums of Hodge numbers. For singular varieties, ordinary Hochschild homology need not have finite-dimensional total homology, but for general proper varieties one can use cdh-sheafified Hochschild homology. Our main construction extends the latter from varieties to stable enhanced categories, thereby proving its Morita-invariance. We combine the resulting numerical invariants with invariants from topological $\K$-theory and the Poincar\'e duality map to recover mixed Hodge numbers for several classes of singular varieties.

The principal geometric results are the following
(Theorems~\ref{thmnormalsurfacereconstruction},
\ref{thmthreefoldreconstruction}, and
\ref{thmfourfoldreconstruction}).

\begin{theorem}\label{thmintroHodgeinvariance}
Let $X$ and $Y$ be connected projective complex varieties which are Fourier--Mukai partners. Then
\[
h_n^{a,b}(X)=h_n^{a,b}(Y)
\qquad\text{for every }n,a,b
\]
in each of the following cases:
\begin{enumerate}
\item
$X$ and $Y$ are normal surfaces.

\item
$X$ and $Y$ are threefolds with rational singularities and
\[
h^1(X,\mathcal O_X)=h^1(Y,\mathcal O_Y).
\]

\item
$X$ and $Y$ are fourfolds which are rational homology manifolds, and
\[
h^{0,1}_1(X)=h^{0,1}_1(Y),
\qquad
h^{0,2}_2(X)=h^{0,2}_2(Y).
\]
\end{enumerate}
\end{theorem}

In this paper, by Fourier--Mukai partners, we mean  varieties having equivalent enhanced categories of perfect complexes; some equivalent formulations are recalled in Subsection~\ref{derivedcategories}.
Assuming that the corresponding numbers vanish on both partners in the threefold and fourfold statements already gives many interesting classes of examples. In Section~8, we give explicit formulas describing the mixed Hodge numbers in each case in terms of categorical invariants and the numbers appearing in the hypotheses, and discuss some examples.

\medskip

We first describe the invariant replacing Hochschild homology in the singular setting. Write $\Catperf_k$ for the category of small idempotent-complete $k$-linear stable categories and $\Mod_k^{BS^1}$ for $k$-modules with a circle action. For a categorical invariant $F$ and a scheme $X$, we write $F(X):=F(\Perf(X))$.

\begin{theorem}\label{thmHHDBdefining}
Let $k$ be a field of characteristic zero. There is a morphism of truncating localizing invariants
\[
R_{HH}:\HH^\mathrm{E}_{\DB}(-/k)\ra\HH^\cl_{\DB}(-/k):
\Catperf_k\longrightarrow\Mod_k^{BS^1}
\]
and natural comparisons
$\HH^{?}(-/k)\to\HH^?_{\DB}(-/k)$, $?\in\{\E, \cl\}$, with the following properties.
Put
\[
\HC^{-,?}_{\DB}(-/k):=\HH^?_{\DB}(-/k)^{hS^1}.
\]

\begin{enumerate}
\item
On every separated finite-type $k$-scheme or, more generally, a tame (classical) ANS stack $X$ of finite type over $k$, $R_{HH}$ is an equivalence and there are canonical equivalences
\[
\HH^?_{\DB}(X/k)\simeq\HH_{\cdh}(X/k),
\qquad
\HC^{-,?}_{\DB}(X/k)\simeq\HC^-_{\cdh}(X/k).
\]
The first is an equivalence of mixed complexes, and the second is an equivalence of $\HC^-(k/k)$-modules; $\HC^-(k/k)\simeq k[[u]]$.
In particular, the comparison with ordinary Hochschild homology is an equivalence on smooth $k$-schemes.

\item
For every extension of fields $k\subset k'$ and every separated finite-type $k$-scheme $X$, there are natural equivalences
\[
\HH^?_{\DB}(X/k)\otimes_k k'
\simeq
\HH^?_{\DB}(X_{k'}/k')
\]
compatible with the circle actions.

\item
For $k=\bbQ$, there are natural commutative squares
\[
\begin{tikzcd}[column sep=large]
\K(T)\arrow[r,"\ch"]\arrow[d]
  &\HC^{-,?}(T/\bbQ)\arrow[d]\\
\KH(T)\arrow[r]
  &\HC^{-,?}_{\DB}(T/\bbQ).
\end{tikzcd}
\]
The fibers of the vertical maps are $\bbQ[[u]]$-modules, and the natural map between them is the unit of derived $u$-completion. In particular, the square is cartesian for $T=\Perf(X)$, where $X$ is a separated finite-type $\bbQ$-scheme.

\item The comparison map $\HH^{E}(-/k)\to\HH^E_{\DB}(-/k)$ is a map of lax symmetric monoidal functors.

\item $\HH_\DB^\cl$ is finitary, i.e. it preserves filtered colimits.
\end{enumerate}
\end{theorem}

Here $\K$ is nonconnective algebraic $\K$-theory, $\KH$
is homotopy $\K$-theory, $\HH$ is Hochschild homology, $\HC^{-,\cl}$ is classical negative cyclic homology, and $\HC^{-,\E}$ is its Efimov counterpart. A localizing invariant is \emph{truncating} if it sends
$R\to\pi_0R$ to an equivalence for every connective
$k$-algebra $R$ \cite{LT19}.

The theorem combines the construction of Subsection~4.1
with Theorems~\ref{DBstacksdescent} and~\ref{thmcdhTCDB} and the geometric base-change properties of the Du Bois complex. For $X\in\Sch^\sft_k$, the cdh comparison gives an identification of underlying $k$-modules
\[
\HH_{\DB}(X/k)\simeq
\bigoplus_{p\geq0}
\Gamma(X,\underline\Omega^p_{X/k})[p],
\]
with circle action corresponding to the de Rham differential.
For a proper $X\in\Sch^\sft_k$, this has finite-dimensional total homology.
Consequently the \emph{Du Bois--Hochschild numbers}
\[
hh_i(X):=\dim_k\HH_{\DB,i}(X/k)
\]
satisfy
\[
hh_i(X)
=
\sum_{p-q=i}\underline h^{p,q}(X)
=
\sum_{\substack{n,a,b\\2a-n=i}}h_n^{a,b}(X),
\]
where $\underline h^{p,q}$ are the Du Bois numbers.
These sums are therefore invariant under Fourier--Mukai equivalence. Moreover, the Hodge--de Rham spectral sequence from $\HH_{\DB}(X/k)$ to $\HP_{\DB}(X/k)$ degenerates at
$E_2$ (Corollary~\ref{DBHodgeformula}). We also show (Theorem~\ref{thmpurity}) that the numbers $hh_i$ detect cohomological purity: cohomology in each degree is pure of the expected weight if and only if
\[
hh_i(X)=hh_{-i}(X)
\qquad\text{for every }i.
\]
In particular, having pure cohomology is a Morita-invariant property of proper varieties.

These statements are not restricted to varieties over $\bbC$.
Subsection~\ref{secDBrecollections} defines the Du Bois and mixed Hodge numbers of proper varieties over an arbitrary field of characteristic zero, using the Hodge and weight filtrations on derived de Rham cohomology, similarly to \cite{Bha12}. The statements above stay true in this generality.

\medskip

To separate individual mixed Hodge numbers within the sums above, we need additional linear equations. To obtain them, we use the comparison map
\[
\Theta_X:
\K^{\mathrm{top}}(\Perf(X))
\longrightarrow
\K^{\mathrm{top}}(\Coh(X)).
\]
Here $\K^{\mathrm{top}}$ is Blanc's topological $\K$-theory
\cite{Bla16}, and $\Coh(X)$ denotes the enhanced bounded derived category of complexes with coherent cohomology. For proper schemes, the inclusion
$\Perf(X)\to\Coh(X)$ is determined by $\Perf(X)$.
Thus $\Theta_X$ supplies an invariant of Fourier--Mukai partners of its own interest.

For a projective complex variety, the comparisons with topological $\K$-cohomology and Borel--Moore
$\K$-homology identify $\Theta_X$, under rational Atiyah-Hirzebruch isomorphisms, with
\[
-\cap\td_*(X):
H^{\mathrm{ev/odd}}(X,\bbQ)
\longrightarrow
H_{\mathrm{ev/odd}}(X,\bbQ),
\]
where $\td_*(X)$ is the Baum--Fulton--MacPherson Todd class of $\mathcal O_X$
(Theorem~\ref{proptopcomparison}).
The comparison identifies the action coming from $\Perf(X)^\otimes$-module structure on $\Coh(X)$ with the topological cap-product action. Its compatibility with pushforward to a point expresses the categorical index pairing in terms of the same Todd class.

For an irreducible $X$, the top component of $\td_*(X)$ is the fundamental class, so the diagonal blocks of this map are the cap-product maps
\[
-\cap[X]:H^j(X,\bbQ)\longrightarrow H_{2\dim X-j}(X,\bbQ).
\]
This allows us to use the categorical construction above to measure the failure of Poincar\'e duality.

The resulting \emph{topological defect numbers}
\[
d_i(X):=
\dim_{\bbQ}\ker(\pi_i\Theta_X\otimes\bbQ)
\]
depend only on the parity of $i$ and are invariant under Fourier--Mukai equivalence. For connected normal projective surfaces, we prove
\[
d_0(X)=\dim_{\bbQ}W_1H^2(X,\bbQ),
\qquad
d_1(X)=b_3(X)-b_1(X).
\]
Together with the $hh_i$, these determine all mixed Hodge numbers. For connected projective threefolds with rational singularities, we instead obtain
\[
d_0(X)=\sigma(X),
\qquad
\sigma(X):=
\dim_{\bbQ}(\Cl(X)/\Pic(X))_{\bbQ}.
\]
Thus the $\bbQ$-factoriality defect is a derived invariant within this class. Combining it with the results of Park--Popa \cite{PP25a,PP25b} recovers all mixed Hodge numbers once the irregularity is fixed.

For projective rational homology manifolds, Poincar\'e duality holds and the defect numbers vanish. The index pairing nevertheless provides additional information.
For a connected projective complex rational homology manifold of complex dimension $2d$, its symmetrization
$S_X$ satisfies
\[
\sig(S_X)=(-1)^d\sig_{\mathrm{top}}(X).
\]
This extends the signature argument of \cite{Ant26}
to the singular setting. Intersection cohomology supplies the Hodge--Riemann relations needed to express the signature in terms of Hodge numbers. In dimension four, this gives the additional equation allowing us to describe all the mixed Hodge numbers in terms of the numbers $\underline{h}^{0,1}=h^{0,1}_1$ and $\underline{h}^{0,2}=h^{0,2}_2$ and categorical invariants.

\medskip

We now describe the construction underlying the Du Bois invariants. The first invariant for which we define a Du Bois counterpart is topological cyclic homology. This invariant comes in three flavors: genuine, Nikolaus-Scholze's, and Efimov's. The last version was introduced in \cite{Efi25b} where rigidity of the category $\Motloc_k$ of localizing motives was proved and used, together with the recollement relating the $\bbA^1$-invariant and cyclotomic subcategories of $\Motloc_k$, to compare it with the Nikolaus-Scholze version.

We use the multiplicative structure of the recollement to prove the following result, which is Corollary~\ref{multiplicativeinvariants}
together with Proposition~\ref{TCDBtruncating}.

\begin{theorem}\label{thmintromultiplicative}
Let $k$ be a
commutative ring spectrum.
There is a truncating lax symmetric monoidal localizing invariant $\TC^{\mathrm E}_{\DB}$ fitting into a cartesian square of lax symmetric monoidal localizing invariants
\[
\begin{tikzcd}
\K\arrow[r]\arrow[d]
  &\TC^{\mathrm E}\arrow[d]\\
\KH\arrow[r]
  &\TC^{\mathrm E}_{\DB}.
\end{tikzcd}
\]
In particular, evaluation on a symmetric monoidal
$k$-linear category gives a cartesian square of commutative ring spectra.

Set
\[
\DK:=\fib(\K\to\KH),
\qquad
\K^{\mathrm{inv},\mathrm E}:=
\fib(\K\to\TC^{\mathrm E}).
\]
For every $M\in\Motloc_k$, the triangles
\[
\DK(M)\longrightarrow
\TC^{\mathrm E}(M)\longrightarrow
\TC^{\mathrm E}_{\DB}(M)
\]
and
\[
\K^{\mathrm{inv},\mathrm E}(M)\longrightarrow
\KH(M)\longrightarrow
\TC^{\mathrm E}_{\DB}(M)
\]
are exact triangles of $\TC^E(k)$-modules and
$\KH(k)$-modules, respectively. When $k$ is connective, $\TC^E(k)\simeq \TC(k)$.
\end{theorem}

Here $\TC^{\mathrm E}$ is Efimov's topological cyclic homology\footnote{It is sometimes also referred to as \emph{motivic} topological cyclic homology. However, by \cite{RSW25} (see Remark~\ref{omegaoneextension}), all the (topological) cyclic homology theories we consider are motivic in the sense that they factor through $\Motloc_k$, so reserving this flattering term for only one of them feels unfair.}.
This statement was also independently obtained by T.~Nikolaus and M.~Ramzi.

The square of lax symmetric monoidal invariants above demonstrates importance of working with Efimov's topological cyclic homology. Using Efimov's version, we also construct genuine and Nikolaus--Scholze versions of $\TC_{\DB}$. All three are truncating and agree on weakly geometric motives, in particular on motives of schemes. On schemes over prime fields, the square above recovers the cdh-descent square appearing in \cite{EM26}. 

For a $\bbQ$-dg-category $T\in\Catperf_\bbQ$, its Hochschild homology over $\bbQ$ can be recovered from its Nikolaus-Scholze topological cyclic homology as the cofiber of multiplication by $u\in\bbQ[[u]]\simeq\TC(\bbQ)$. We combine this observation with the base change property of Hochschild homology and define for a field $k$ of characteristic zero and $T\in\Catperf_k$
\[
\HH^?_{\DB}(T/k):=
\TC^?_{\DB}(T)
\otimes_{\HC^-(k/\bbQ)}\HH(k/k),
\qquad ?\in\{\mathrm E,\cl\},
\]
with its canonical circle action; here $\cl:=\NS$. We check that the resulting invariant recovers cdh-sheafified homology on schemes which essentially allows us to prove Theorem~\ref{thmHHDBdefining} and first Hodge numbers reconstruction statements.

Homotopy fixed points, homotopy orbits, and Tate invariants give negative cyclic, cyclic, and periodic cyclic homology, respectively. On schemes, the Efimov and classical versions agree, so we omit the flavour superscripts there.

These theories satisfy the following extension of the Geller--Weibel theorem \cite{GW89}.

\begin{theorem}\label{thmintroperiodic}
Let $k$ be a field of characteristic zero.
For both the classical and Efimov versions, with or without the Du Bois modification, periodic cyclic homology is $\bbA^1$-invariant:
\[
\HP^?_\dagger(T/k)
\xrightarrow{\sim}
\HP^?_\dagger(T\otimes_k k[t]/k),
\qquad
?\in\{\mathrm E,\cl\},
\quad
\dagger\in\{\varnothing,\DB\}.
\]
The $\bbA^1$-homotopifications of
$\HH^?_\dagger$ and $\HC^?_\dagger$ vanish.
Consequently,
\[
\HC^{-,?}_\dagger(-/k)
\longrightarrow
\HP^?_\dagger(-/k)
\]
is the $\bbA^1$-homotopification map.
\end{theorem}

The nontrivial point of Theorem~\ref{char0GW} is $\bbA^1$-invariance of $\HP^E_\dagger$; the proof again uses Efimov's recollement and his corepresentability theorem for $\TR$ of connective algebras.
On schemes, the classical and Du Bois periodic theories agree:
\[
\HP(X/k)\xrightarrow{\sim}\HP_{\DB}(X/k).
\]
For proper schemes over
$\bbQ$ we show that there is a split triangle
\[
\DK(X)\longrightarrow\Sigma\HC(X/\bbQ)
\longrightarrow\Sigma\HC_{\DB}(X/\bbQ)
\]
of $\bbQ$-modules
(Corollary~\ref{DKsplittingvarieties}), which is useful for computations of algebraic K-theory of singular varieties.

We also define a Du Bois version of Blanc's Deligne cohomology:
\[
H_{\mathcal D,\DB}(T):=
\K^{\mathrm{top}}(T)
\mathop{\times}_{\HP^{\cl}_{\DB}(T/\bbC)}
\HC^{-,\cl}_{\DB}(T/\bbC).
\]
This is a truncating localizing invariant receiving a Chern character from $\KH$. On proper complex schemes, its rationalization recovers periodized classical Deligne--Beilinson cohomology and the comparison from Blanc's Deligne cohomology to its Du Bois counterpart is induced, weight by weight, by the filtered map from Hodge-completed derived de Rham cohomology to the Du Bois complex (Corollary~\ref{DeligneDBrational}).

\medskip

\subsubsection*{Structure of the paper}

Section~2 recalls the categorical, motivic, and Hodge-theoretic background, including some consequences of Efimov's rigidity theorem \cite{Efi25b} used in this paper. In Section~3 we define the Du Bois versions of topological cyclic homology, together with their multiplicative structures, for categories over an arbitrary base. In Section~4 we specialize to characteristic zero, define Du Bois Hochschild and cyclic homology, and study their properties. Section~5 computes the Du Bois invariants on schemes and stacks and shows, in particular, that the column sums of Du Bois numbers of a proper variety are determined by its perfect derived category.
Section~6 uses the Du Bois--Hochschild numbers to prove the purity criterion and first reconstruction results for surfaces. Section~7 uses topological K-theory and Poincar\'e duality to define the topological defect numbers and to give a categorical description of the topological signature.
Section~8 combines these results to reconstruct mixed Hodge numbers for large classes of surfaces, threefolds, and fourfolds, and concludes by showing that individual Hodge numbers need not be determined by the underlying localizing motive alone.

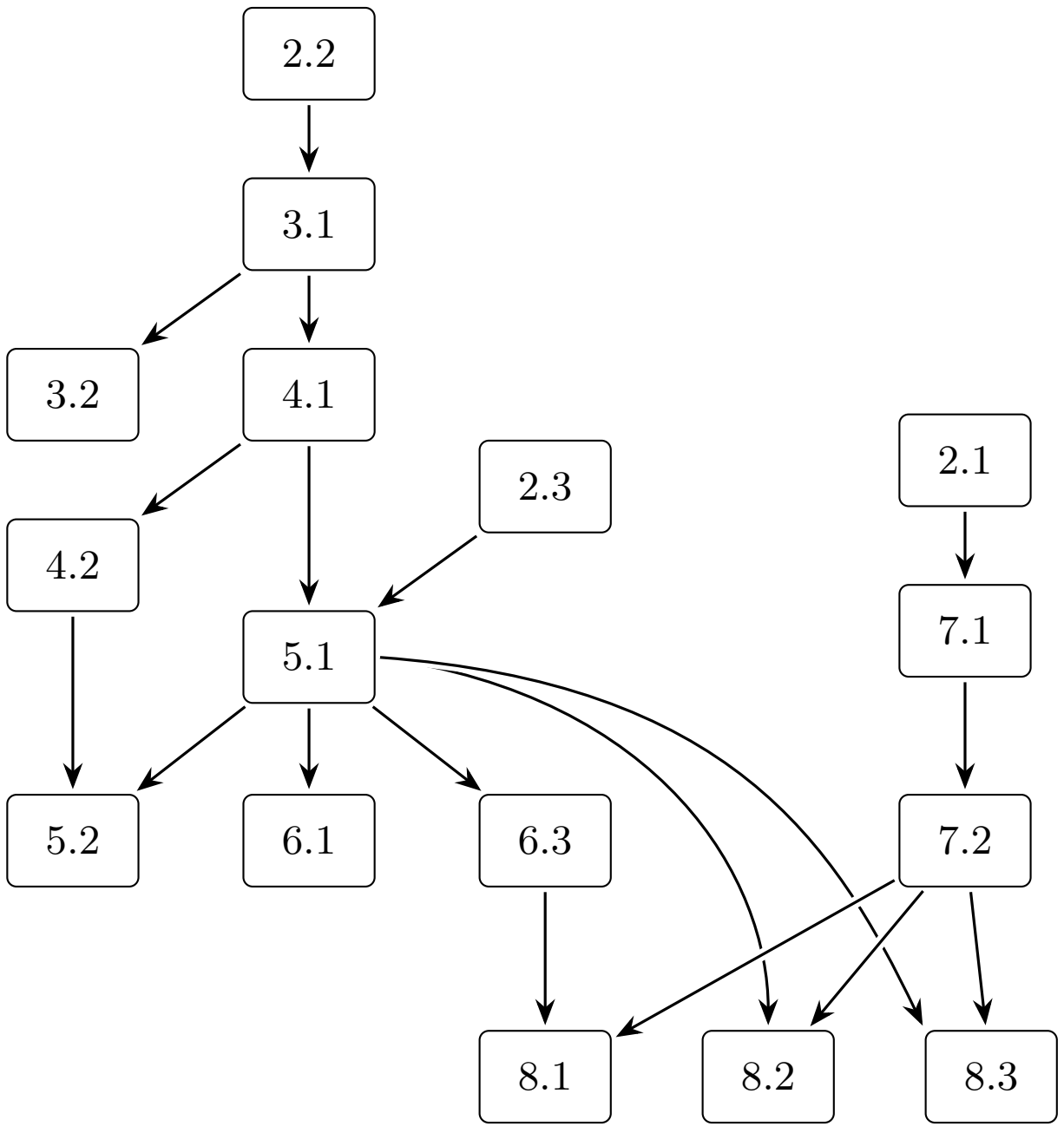
\begin{figure}[htbp]
\centering
\begin{tikzpicture}[
  x=1cm,
  y=1cm,
  subsection/.style={
    draw,
    rounded corners=2pt,
    fill=white,
    minimum width=10mm,
    minimum height=7mm,
    inner sep=2pt,
    font=\small
  },
  dependency/.style={
    -{Stealth[length=2mm]},
    semithick,
    shorten >=1pt,
    shorten <=1pt,
    preaction={draw=white,line width=2.5pt}
  }
]

\node[subsection] (s22) at (0,0)     {2.2};
\node[subsection] (s31) at (0,-1.3)  {3.1};
\node[subsection] (s32) at (-1.8,-2.6) {3.2};
\node[subsection] (s41) at (0,-2.6)  {4.1};

\node[subsection] (s42) at (-1.8,-3.9) {4.2};
\node[subsection] (s23) at (1.8,-3.3)  {2.3};
\node[subsection] (s51) at (0,-4.6)    {5.1};

\node[subsection] (s21) at (5,-3.1) {2.1};
\node[subsection] (s71) at (5,-4.4) {7.1};
\node[subsection] (s72) at (5,-6)   {7.2};

\node[subsection] (s52) at (-1.8,-6) {5.2};
\node[subsection] (s61) at (0,-6)    {6.1};
\node[subsection] (s63) at (1.8,-6)  {6.3};

\node[subsection] (s81) at (1.8,-7.8) {8.1};
\node[subsection] (s82) at (3.5,-7.8) {8.2};
\node[subsection] (s83) at (5.2,-7.8) {8.3};

\draw[dependency] (s22) -- (s31);
\draw[dependency] (s31) -- (s32);
\draw[dependency] (s31) -- (s41);
\draw[dependency] (s41) -- (s42);
\draw[dependency] (s41) -- (s51);
\draw[dependency] (s23) -- (s51);

\draw[dependency] (s42) -- (s52);
\draw[dependency] (s51) -- (s52);
\draw[dependency] (s51) -- (s61);
\draw[dependency] (s51) -- (s63);

\draw[dependency]
  (s51.east)
  .. controls (2.7,-4.9) and (3.5,-6.3) ..
  (s82.north);
\draw[dependency]
  (s51.east)
  .. controls (3.2,-4.8) and (4.0,-6.0) ..
  (s83.north west);
\draw[dependency] (s63) -- (s81);

\draw[dependency] (s21) -- (s71);
\draw[dependency] (s71) -- (s72);
\draw[dependency] (s72) -- (s81);
\draw[dependency] (s72) -- (s82);
\draw[dependency] (s72) -- (s83);

\end{tikzpicture}
\caption{Main dependencies between subsections, omitting Subsections~6.2 and~8.4.}
\label{figsubsectiondependencies}
\end{figure}

\medskip
\subsubsection*{Conventions}

By a category, we mean an $(\infty,1)$-category.

We use $\Hom$ for Hom-spectra. $\underline{\Hom}$ is used for internal hom objects.

For a dg-algebra or a classical algebra $A$, we think of $A$-dg-modules as $HA$-modules in $\Sp$. In particular, by a "complex" we mean $H\bbZ$-module.

Therefore all tensor products, Homs, pushforwards, pullbacks, completions are "derived".

We freely switch between homological and cohomological notation. By default, all degrees are  homological; for a complex, $\pi_i=H_i=H^{-i}$ and an element there is said to have the degree $|-|=i$. 

For filtrations a superscript comes with decreasing filtrations and a subscript with increasing. So the weight filtration comes with a subscript, even though its graded pieces are related to cohomological grading.

We will usually work over a base $E_\infty$-ring spectrum $k$. We write $\Catperf_k$ for the  category of small idempotent-complete $k$-linear stable  categories, to which we will usually refer to as $k$-categories. When $k$ is an Eilenberg-Maclane ring spectrum, the corresponding theory is equivalent to Morita theory of enhanced $k$-linear triangulated categories, for example, to the theory of $k$-dg-categories up to Morita equivalence.

We write $\Sch^\sft_k$ for the 1-category of separated schemes of finite type over a classical ring $k$. By a $k$-variety, we mean a reduced separated scheme of finite type over a field $k$. We do not assume irreducibility.

We write
\[
b_i(X):=\dim_{\bbQ}H^i(X^\an,\bbQ)
\]
for the Betti numbers of a complex variety $X$.

\subsubsection*{Acknowledgements}

I am especially grateful to Alexander Efimov and Maxim Ramzi for stimulating conversations that helped me develop my intuition about localizing motives and appreciate the power of the rigidity theorem. I am also grateful to Christopher Brav, Dmitry Kaledin, Marc Levine, and Vladimir Sosnilo for helpful discussions.

\medskip

\section{Recollections}

We recall the categorical and geometric constructions used throughout the paper and fix our notation.

In the first subsection, we discuss derived categories and their equivalences. We compare several natural notions of Fourier--Mukai partners, justifying our choice. The material here is classical.

In the second subsection, we recall the theory of localizing motives, introduced by Blumberg--Gepner--Tabuada, together with some modern developments, most importantly Efimov's rigidity theorem. We recall Efimov's $\bbA^1$-recollement and use it to define several important localizing invariants. The material of this subsection is essentially contained in \cite{Efi25b}; a straightforward but important novelty is the introduction of the subcategory of weakly geometric motives.
As this paper demonstrates, this category is a natural domain on which many localizing invariants behave especially well.

Finally, in the third subsection, we recall cdh descent and the Du Bois complex, with its two filtrations. We make a point of working over an arbitrary field of characteristic zero, replacing singular cohomology by derived de Rham cohomology as in \cite{Bha12}. We show that the associated Du Bois and mixed Hodge numbers can be defined in this generality, that their definitions neither require nor depend on an embedding of the field into $\bbC$, and that they have the expected properties, including the Hodge symmetry of mixed Hodge numbers. These facts will be used in later sections.

\subsection{Derived categories and Fourier-Mukai partners}\label{derivedcategories}

In this subsection, $k$ is an excellent Noetherian ring of finite Krull dimension.

For a qcqs scheme $X/k$, write $\QCoh(X)$ for the presentable $k$-linear stable category of complexes with quasi-coherent cohomology, and
$\Perf(X)$ for its compact objects. Equivalently, perfect complexes locally are bounded complexes of finite-rank vector bundles. Pullbacks preserve perfect complexes, so $\Perf$ defines a contravariant functor in $X$. If $F$ is an invariant of small stable categories, we write
$F(X)=F(\Perf(X))$.

For a Noetherian $X/k$, let $\Coh(X)$ denote the category of bounded complexes with coherent cohomology. There is a fully faithful inclusion $\Perf(X)\subset\Coh(X)$, which is an equivalence when $X$ is regular.

The following statement is somewhat folklore and exact attribution depends on assumptions; for proper $k$-schemes we can refer to the results of Neeman and Canonaco--Neeman--Stellari; see also \cite[Lemma~4.7]{Bal11}
and \cite[Remark 3.14]{KL15}
for the projective over a field case.

\begin{prop}\label{propderivedequivalence}
Let $X$ and $Y$ be proper $k$-schemes. The following are equivalent:
\begin{enumerate}
    \item[(i)] there exists an exact $k$-linear equivalence
    $\Perf(X)\simeq\Perf(Y)$;
    \item[(ii)] there exists an exact $k$-linear equivalence
    $\Coh(X)\simeq\Coh(Y)$;
    \item[(iii)] there exists a $k$-linear triangulated equivalence of 1-categories
    $\Ho(\Perf(X))\simeq
    \Ho(\Perf(Y))$;
    \item[(iv)] there exists a $k$-linear triangulated equivalence of 1-categories
    $\Ho(\Coh(X))\simeq
    \Ho(\Coh(Y))$.
\end{enumerate}
\end{prop}

\begin{proof}

The perfect objects are precisely those
$P\in\Coh(X)$ for which
\[
\bigoplus_{j\in\mathbb Z}
\Hom_{\Ho(\Coh(X))}(P,E[j]) 
\qquad\text{is a finitely generated $k$-module for every }E\in\Coh(X);
\]
(for example, by
\cite[Example~1.7(iii)]{Nee26}), hence (iv)$\Rightarrow$(iii).

(iii)$\Rightarrow$(i) 
follows from uniqueness of enhancements
(e.g. \cite[Theorem~B]{CNS22}).

(i)$\Rightarrow \QCoh(X)\simeq \Ind\Perf(X)\simeq \Ind\Perf(Y)\simeq \QCoh(Y)\Rightarrow$(ii), since $\Coh(X)$ can be intrinsically characterized in terms of $\QCoh(X)$
(\cite[Remark~5.8]{Nee26}, see also \cite[Remark 3.14]{KL15}).

(ii)$\Rightarrow$(iv) is immediate, finishing the proof.

\end{proof}

\begin{defi}\label{defiFM}
    We say that $X/k$ and $Y/k$ are Fourier-Mukai partners if there exists an exact $k$-linear equivalence $\Perf(X)\simeq \Perf(Y)$.
\end{defi}

In such a situation, we always assume such an equivalence chosen.
By \cite[Theorem~8.9 and Corollary~8.12]{Toe07}, it is determined by a Fourier--Mukai kernel
$K\in\QCoh(X\times_kY)$ (the product here is derived; if $X/k$ is flat, e.g. if $k$ is a field, it is isomorphic to the classical product).
The possible choices of an equivalence differ by an autoequivalence. An autoequivalence of a category induces automorphisms on its invariants, hence equalities on numerical invariants.

\medskip

Homological projective duality provides many examples of singular Fourier--Mukai partners\footnote{We thank A.~Kuznetsov for pointing out this source of examples.}. We recall one classical construction.

\begin{example}[Kuznetsov]\label{examplePfaffian}
Let $V$ be a seven-dimensional $\bbC$-vector space and let $L$ be a seven-dimensional subspace in $\bigwedge^2V^*$. Consider the intersections in $\Sch^\sft_\bbC$ 
\begin{equation}\label{Kuznetsovcap}
X_L=\Gr(2,V)\times^{\cl}_{\bbP(\bigwedge^2V)}\bbP(L^\perp),
\qquad
Y_L=\operatorname{Pf}_{\leq4}(V^*)\times^{\cl}_{\bbP(\bigwedge^2V^*)}\bbP(L),
\end{equation}
where $\operatorname{Pf}_{\leq4}(V^*)\subset
\bbP(\bigwedge^2V^*)$ is the Pfaffian variety of skew forms of rank at most four.
Whenever
\begin{equation}\label{Kuzdim}
\dim X_L=\dim Y_L=3,
\qquad
\bbP(L)\cap\Gr(2,V^*)=\varnothing,
\end{equation}
\cite[Theorem~11.4]{Kuz06}, combined with Proposition~\ref{propderivedequivalence}, gives an equivalence
\[
\Perf(X_L)\simeq\Perf(Y_L).
\]
The assumptions impose the expected dimensions but do not require transverse intersections; suitable choices of $L$ give positive-dimensional singular loci. The dimension assumptions ensure that the classical intersections in \eqref{Kuznetsovcap} agree with their derived counterparts. Both schemes are local complete intersections with trivial dualizing sheaves.
The construction and equivalence also hold over any algebraically closed field of characteristic zero.
\end{example}

Further singular examples obtained from homological projective duality include expected-dimensional intersections of two Grassmannians \cite[Theorem~1.1]{BCP20} and of two spinor varieties \cite[Proposition~4.2 and the subsequent remark]{Man19};
see also \cite[Theorems~6.1 and~6.3]{KP21} for their formulation using derived intersections.
Other sources include flops of threefolds with terminal Gorenstein singularities \cite[Theorem~1.1]{Che02} and fine compactified Jacobians of a fixed reduced connected projective curve with planar singularities \cite[Theorem~A]{MRV19}.

\subsection{Localizing motives and Efimov's $\bbA^1$-recollement}

Throughout, $k$ is an $E_\infty$-ring unless a stronger hypothesis is stated.
We write $\Catperf_k$ for the  category of small idempotent-complete
$k$-linear stable  categories.

\begin{defi}
A diagram $\cA\to\cB\to\cC$ in $\Catperf_k$ is an \emph{exact} or
\emph{Karoubi} triple if the square
\[
\begin{tikzcd}
\cA \arrow[r] \arrow[d] & \cB \arrow[d] \\
0 \arrow[r] & \cC
\end{tikzcd}
\]
is bicartesian.
\end{defi}

Equivalently, $\cA\to\cB$ is fully faithful, the composite
$\cA\to\cC$ is zero, and the induced functor from the Verdier quotient
$\cB/\cA$ to $\cC$ is an equivalence after idempotent completion.

\begin{example}\label{Nisnevich}\label{exlocalization}
Let $X$ be a qcqs $k$-scheme, $j\colon U\subset X$ a quasi-compact open subscheme, and $Z=X\setminus U$ the closed complement on underlying topological spaces.  Then
\[
\Perf_Z(X)\longrightarrow\Perf(X)\xrightarrow{j^*}\Perf(U)
\]
is a Karoubi triple. This localization sequence implies that localizing invariants (Definition~\ref{locdefi}), restricted to qcqs schemes, satisfy Nisnevich (in particular, Zariski) descent.
\end{example}

\begin{example}\label{Sgcategory}
Let $X$ be a Noetherian classical $k$-scheme.
It is often convenient to consider its (idempotent-completed) singularity category
\[
\Sg(X)
 :=\bigl(\Coh(X)/\Perf(X)\bigr)^{\kappa}.
\]
Thus $\Perf(X)\to \Coh(X)\to \Sg(X)$ is an exact triple in $\Catperf_k$. The proof of Proposition~\ref{propderivedequivalence} shows that this exact triple is completely determined by $\Perf(X)$, hence Fourier-Mukai partners produce equivalent triples. This observation is used in Section~\ref{sectopologicalcomparison}.
\end{example}

\begin{defi}\label{locdefi}
Let $\cE$ be a stable  category.  A functor
$F\colon\Catperf_k\to\cE$ is
\emph{localizing} if $F(0)=0$ and $F$ sends exact triples to exact triangles.  We say it is  \emph{finitary} if it preserves filtered colimits.
\end{defi}

In this paper, we use several localizing invariants without recalling their definitions and general theory. For Hochschild and cyclic homology we refer to \cite{Kel99}; for topological cyclic homology and its genuine and Nikolaus--Scholze versions, to \cite{NS18}; for topological restriction homology, to \cite{KN17}; and for topological $\K$-theory and its relative version, to
\cite{Bla16,Mou19}. We recall below the motivic descriptions of algebraic and homotopy $\K$-theory and introduce the new invariants used here. For the first treatment of algebraic $\K$-theory of stable categories, we refer to \cite{Sch11}, and for the classical theory to \cite{Wei13}.

\medskip

By Blumberg--Gepner--Tabuada \cite{BGT13} there is a universal finitary localizing invariant
\[
\Uloc:=\Uloc^k\colon\Catperf_k\longrightarrow\Motloc_k.
\]
The target is a presentably
symmetric monoidal stable category, hence closed. We will omit $k$ from the $\Uloc$ notation.

\begin{remark}\label{omegaoneextension}
By Ramzi--Sosnilo--Winges \cite{RSW25}, every $\omega_1$-finitary localizing invariant factors essentially uniquely through the category of localizing motives, as an exact functor preserving $\omega_1$-filtered colimits.
The invariants used in the paper satisfy this condition. We therefore occasionally use the same notation for these invariants on categories and for their extensions to motives.
\end{remark}

We write $\Uloc(X):=\Uloc(\Perf(X))$ for the motive of a scheme.

\begin{example}\label{lciblowups}
Let $Z\hookrightarrow X$ be a closed embedding of smooth schemes of codimension $c\geq1$, and let $E$ be the exceptional divisor. The blowup and projective bundle formulas give
\[
\Uloc(\Bl_ZX)\simeq
\Uloc(X)\oplus\Uloc(Z)^{\oplus(c-1)},\qquad
\Uloc(E)\simeq\Uloc(Z)^{\oplus c};
\]
see \cite{Orl93}. These decompositions are compatible with the maps in the blow up square.
\end{example}

Extension and restriction of scalars induce an adjunction
\begin{equation}\label{basemotloc}
k\otimes -:\Motloc_\bbS\rightleftarrows\Motloc_k:\Res_k,
\qquad k\otimes -\dashv \Res_k.
\end{equation}

Here the left adjoint is symmetric monoidal, hence rigidity implies that the right adjoint preserves colimits and the projection formula holds:
\[
\Res_k\bigl(M\otimes_k(k\otimes N)\bigr)
\simeq \Res_kM\otimes_\bbS N.
\]
See, for example,
\cite[Propositions 1.9(i) and 1.12(iv)]{Efi25a}.

We use the notation
\[
\Gamma(M):=\Hom_{\Motloc_k}(\mathbbm 1,M).
\]
The corepresentability theorem of \cite{BGT13} gives
\[
\Gamma\Uloc(T)\simeq \K(T),
\]
which we will use as the definition of algebraic K-theory.
A.~Efimov proved that $\Motloc_k$ is rigid in the sense of Gaitsgory-Rozenblyum \cite{Efi25b}.  It implies, in particular, that every $\Sp$-valued finitary localizing invariant can be realized as $\Motloc(\mathbbm{1}, \Uloc(-)\otimes M)$ for some incarnating motive $M\in\Motloc_k$. For general localizing invariants, situation is more subtle, but we will see that in many cases of interest one can still get useful descriptions in terms of $\Motloc$  .

We introduce an important subcategory of $\Motloc_k$.

\begin{defi}
The \emph{weakly geometric} subcategory or the category of \emph{weakly geometric motives}
$\Motlocwg_k\subset\Motloc_k$ is the smallest thick stable subcategory containing the motives of connective $E_1$-$k$-algebras.
\end{defi}

Note that Example~\ref{exlocalization} implies that the motives of qcqs schemes are weakly geometric.
 On $\Motlocwg_k$, the bestiary of cyclic-homology-like theories simplifies dramatically.

\begin{remark}\label{LevySosnilo}
    $\Motlocwg_k$ can also be described as the thick subcategory generated by the motives of $k$-linear $c$-categories of \cite{LS25}.
\end{remark}

\medskip

For the standard algebraic simplex write

\[
\bbS[\Delta^n]
 = \bbS[t_1,\ldots,t_n] = \Sigma^\infty(\bbN^n)_+; \ k[\Delta^n]:=\bbS[\Delta^n]\otimes k .
\]
We refer to §8 in \cite{Efi25b} for the formulas defining all the maps $k[\Delta^m] \ra k[\Delta^n]$ corresponding to order-preserving maps of posets. Put
\[
\Delta_k^n:=\Uloc\bigl(\Perf(k[\Delta^n])\bigr),
\qquad
\widetilde\Delta_k^n:=\fib(\mathbbm 1\longrightarrow\Delta_k^n),
\qquad
A:=\bigl|\Delta_k^\bullet\bigr|,
\qquad
J:=\fib(\mathbbm 1\longrightarrow A)\simeq\bigl|\widetilde\Delta_k^\bullet\bigr|.
\]

As explained in \cite[Sections 1.3 and 8]{Efi25b}, $A$ is an idempotent $E_\infty$-algebra and $J$ a coidempotent $E_\infty$-coalgebra.

\medskip

Recall that the $\bbA^1$-homotopification of a localizing invariant
(or, more generally, an arbitrary functor with values in a category admitting geometric realizations) can be computed by
\begin{equation}\label{A1ization}
F^{\bbA^1}(-)
:=|[n]\longmapsto
F\bigl(-\otimes_k\Perf(k[\Delta^n])\bigr)|.
\end{equation}

If $F$ is a finitary $\Sp$-valued localizing invariant, choose
$M_F\in\Motloc_k$ such that
$F(-)\simeq\Gamma(\Uloc(-)\otimes_kM_F)$. Then
\begin{equation}\label{finitaryA1ization}
F^{\bbA^1}(-)\simeq
|\Gamma(\Uloc(-)\otimes_kM_F\otimes_k\Delta_k^\bullet)|
\simeq\Gamma(\Uloc(-)\otimes_kA\otimes_kM_F).
\end{equation}
For non-finitary invariants see Remark~\ref{nonfinitarity}.

The idempotent algebra $A$ gives a recollement
(\cite[Section 1.3 and Remark 9.6]{Efi25b})
\begin{equation}\label{recollement}
\xymatrix{
\MotlocA_k \ar[rr]|{i_*} &&
\Motloc_k
 \ar@<-2ex>[ll]|{i^*}
 \ar@<2ex>[ll]|{i^!}
 \ar[rr]|{j^*} &&
\Motloccyc_k
 \ar@<-2ex>[ll]|{j_!}
 \ar@<2ex>[ll]|{j_*}.}
\end{equation}
The two recollement triangles can be written
\begin{equation}\label{recol1}
\dk:=(-)\otimes J =j_!j^*
\longrightarrow
\id_{\Motloc_k}
\longrightarrow
i_*i^* = (-)\otimes A =:\kh
\end{equation}
\begin{equation}\label{recol2}
\kinv:=\underline{\Hom}_{\Motloc_k}(A, -) = i_*i^!
\longrightarrow
\id_{\Motloc_k}
\longrightarrow
j_*j^* = \underline{\Hom}_{\Motloc_k}(J, -)=:\tc.
\end{equation}

We will use lower-case notation for endofunctors of motives. On weakly geometric motives, it corresponds to upper-case notation for the spectrum-valued invariants.

\medskip

The cyclotomic structure on $\THH$ gives symmetric monoidal realizations
\[
\Motloc_\bbS\xrightarrow{\,\THH^{\gen}\,}
\CycSp^{\gen}
\longrightarrow\CycSp^{\NS}.
\]
The image of $A$ is zero: in the realization of the polynomial simplex the two endpoint maps identify $0$ and $1$ in $\pi_0$. Thus the image of
$J\to\mathbbm1$ is an equivalence of coalgebras. Mapping out of $J$ and then out of its image gives natural lax symmetric monoidal comparisons
\begin{equation}\label{cyclotomicrealizations}
\TC^{\mathrm{E}}=\Hom_{\Motloc_\bbS}(J,-)
\longrightarrow \TC^{\gen}\longrightarrow \TC^{\NS},
\end{equation}
compatible with the cyclotomic trace. The composite is an equivalence on connective $E_1$-rings by \cite[Theorem 9.3]{Efi25b}. In fact, the proof can be adapted to see that the first map is an equivalence in the same assumptions, but we do not need this: the second arrow is an equivalence on connective $E_1$-rings by the bounded-below comparison in \cite{NS18}.
Consequently both arrows are equivalences on weakly geometric motives.

\begin{defi}\label{defrecolinvariants}
For $T\in\Catperf_k$ set
\[
\begin{aligned}
\DK(T)&:=\Gamma\,\dk(\Uloc(T)),
& \K^{\mathrm{inv},\mathrm{E}}(T)&:=\Gamma\,\kinv(\Uloc(T)),
& \TC^{\mathrm{E}}(T)&:=\Gamma\,\tc(\Uloc(T)).
\end{aligned}
\]
\end{defi}

\begin{obs}[Efimov's comparison]\label{obscomparisonclassical}
Consider the recollement \eqref{recollement} on $\Motloc_k$. Applying $\Gamma$ to the triangles \eqref{recol1}, \eqref{recol2} returns the fiber sequences
\[
\DK\longrightarrow \K\longrightarrow \KH
\]
and
\[
\K^{\mathrm{inv},\mathrm{E}} \lra \K \lra \TC^{\mathrm{E}}
\]
On the weakly geometric subcategory there are canonical equivalences
\[
\K^{\mathrm{inv},\mathrm{E}}(T)\simeq \K^{\mathrm{inv}}(T),
\qquad
\TC^{\mathrm{E}}(T)\simeq \TC(T),
\]
compatible with the cyclotomic trace.
\end{obs}

\begin{proof}
Write subscripts for the base. We have $k\otimes A_\bbS\simeq A_k$ and
$k\otimes J_\bbS\simeq J_k$. Hence adjunction gives
\[
\Hom_{\Motloc_k}(J_k,M)\simeq\Hom_{\Motloc_\bbS}(J_\bbS,\Res_kM),\qquad
\Hom_{\Motloc_k}(A_k,M)\simeq\Hom_{\Motloc_\bbS}(A_\bbS,\Res_kM).
\]
Similarly $\Gamma_k\simeq\Gamma_\bbS \Res_k$, and the projection formula identifies $\Gamma_k(M\otimes_kJ_k)$ and
$\Gamma_k(M\otimes_kA_k)$ with their absolute versions.
These equivalences respect the exact triangle $J\to\mathbbm1\to A$, hence both triangles and the map $\DK\to \TC^{\mathrm{E}}$ are independent of $k$.
The first triangle now follows from corepresentability and
\eqref{finitaryA1ization}; the second comparison follows from
\eqref{cyclotomicrealizations} and exactness.
\end{proof}

\medskip

\subsection{Cdh-topology and the Du Bois complex}\label{secDBrecollections}

Let $k$ be a field of characteristic zero, and let $\Sch^{\sft}_k$
denote separated schemes of finite type over $k$. An abstract blow-up square is a cartesian square
\begin{equation}\label{abs}
\begin{tikzcd}
E\arrow[r]\arrow[d]&\widetilde X\arrow[d,"p"]\\
Z\arrow[r,"i"]&X
\end{tikzcd}
\end{equation}
in which $i$ is a closed immersion, $p$ is proper, and
$\widetilde X\setminus E\to X\setminus Z$ is an isomorphism. The main example is the blow-up $\widetilde X=\Bl_Z X$, with
$E=\widetilde X\times_X Z$ the exceptional locus.
The cdh topology is generated by Nisnevich covers (which include Zariski covers) and the covering families $\{\widetilde X\to X,Z\to X\}$ of these squares. For a presheaf $G$ with values in a stable category $\cE$, $G_{\cdh}$
denotes its cdh sheafification aka cdh hypercohomology. We refer to \cite{Nis89,Voe10,CHSW08} for further details. 

Resolution of singularities allows one to compute cdh sheafification using smooth hyperresolutions. We will say that a presheaf on
$\Sch^{\mathrm{sft}}_k$ \emph{satisfies the smooth blowups condition} if it maps blowup squares with smooth ambient scheme and smooth center to cocartesian squares. If $G$ is a Nisnevich sheaf satisfying this condition, then $G_{\cdh}\simeq G$ on smooth $k$-schemes
\cite[Lemmas~4.3--4.6]{Voe10}. In particular, this applies to presheaves obtained from localizing invariants
(Example~\ref{Nisnevich}, Example~\ref{lciblowups}). For such a presheaf, a bounded cubical hyperresolution of $X$ computes $G_{\cdh}(X)$ as the limit of a finite diagram of its values on smooth schemes.

If $X$ in \eqref{abs} is proper, then so are the other terms. Consequently, the values of $G_{\cdh}$ on proper $k$-schemes are determined by the restriction of $G$ to smooth proper $k$-schemes, and can be computed using finite diagrams of smooth proper $k$-schemes. For the rest of this subsection we focus on proper $k$-schemes, whose category is denoted
$\Sch^{\mathrm{pr}}_k$.

For constructions on $\Sch^{\mathrm{pr}}_k$, we use the
\emph{proper cdh topology}: it is generated by abstract blowup covering families within
$\Sch^{\mathrm{pr}}_k$ and the empty covering family of the empty scheme. We denote its sheafification by
$(-)_{\pcdh}$. If $G$ is a Nisnevich smooth blowup condition sheaf, then $G_{\cdh}|_{\Sch^{\mathrm{pr}}}\simeq (G|_{\Sch^{\mathrm{pr}}})_{\pcdh}$.

\medskip

A smooth proper hyperresolution $\epsilon:X_\bullet\to X$ defines the filtered Du Bois complex
\[
\underline\Omega^\bullet_{X/k}
 =\epsilon_*\Omega^\bullet_{X_\bullet/k},\qquad
\underline\Omega^p_{X/k}
 =\Gr_F^p\underline\Omega^\bullet_{X/k}[p]
 =\epsilon_*\Omega^p_{X_\bullet/k}.
\]
The construction is independent of the hyperresolution as a filtered complex; its decreasing filtration is induced by degree of differential forms \cite{DB81,CHW08,Bha12}.
Pcdh-descent gives a convenient description of its global sections $\Gamma\underline\Omega^\bullet_{-/k}$ and of the two important filtrations on it.

Endowing the derived Hodge-truncations of the global sections of the Hodge-completed derived de Rham complex with their (increasing) canonical truncations filtration gives a functor
\[
\Gamma(-, \widehat{\dR}_{-/k}):=(\Gamma(-, \widehat{\dR}_{-/k}), F_{\der}, W^{\mathrm{dec}}): (\Sch^{\mathrm{pr}}_k)^{\mathrm{op}} \lra \Fun(\bbN\times\bbN^{\mathrm{op}}, \Mod_k); \]\[ W^{\mathrm{dec}}_mF_{\der}^p\Gamma(X, \widehat{\dR}_{X/k}):=\tau^{\le m}F_{\der}^p\Gamma(X, \widehat{\dR}_{X/k}).
\]

\begin{defi}\label{DBcomplex}
    We call the pcdh-sheafification of the functor above \emph{Deligne-Du Bois complex}
    \[
    \Gamma\underline\Omega^\bullet_{-/k}:=(\Gamma\underline\Omega^\bullet_{-/k}, W^{\mathrm{dec}}, F): (\Sch^{\mathrm{pr}}_k)^{\mathrm{op}} \lra \Fun(\bbN\times\bbN^{\mathrm{op}}, \Mod_k).
    \]
    We call the increasing filtration $W^{\mathrm{dec}}$ the \emph{d\'ecalaged weight filtration} and the decreasing filtration $F$ the \emph{Hodge-Deligne} filtration.
\end{defi}

More explicitly, if $X_\alpha$ is the finite diagram of smooth proper schemes in a bounded cubical hyperresolution of $X$, then
\begin{equation}\label{eqDBbifilteredhyperresolution}
W^{\mathrm{dec}}_mF^p\Gamma\underline\Omega^\bullet_{X/k}
 \simeq
 \operatorname*{lim}_{\alpha}
 \tau^{\le m}\Gamma
       (X_\alpha,\Omega^{\geq p}_{X_\alpha/k}).
\end{equation}

Comparing the definitions and using that global sections functor commutes with finite limits, one sees that after forgetting the d\'ecalaged weight filtration one indeed gets $\Gamma\underline\Omega^\bullet_{-/k}$ considered earlier with its decreasing filtration, so the notation is consistent.

By Proposition 5.2 \cite{Bha12}, after forgetting both filtrations, the canonical map $\Gamma\widehat{\dR}_{-/k}\ra \Gamma\underline\Omega^\bullet_{-/k}$ becomes an equivalence. We see that the cdh-sheafification in Definition~\ref{DBcomplex} changes the filtrations while preserving the underlying object of $\Mod_k$.

\begin{remark}
One could instead start with global sections of the naive de Rham complex, endowed with the filtration by form degree and with the canonical truncation of each filtration step, and apply proper cdh sheafification. The comparison map is an equivalence on smooth proper schemes, so this procedure gives the same bifiltered object;
cdh-sheafification, however, does not preserve even the underlying complex for a general proper $k$-scheme.
\end{remark}

\medskip

Over $\mathbb C$, the underlying analytic Du Bois complex is quasi-isomorphic to $\mathbb C_{X^{\an}}$. For proper $X$, the spectral sequence
\[
E_1^{p,q}=\mathbb H^q(X,\underline\Omega^p_{X/\mathbb C})
 \Longrightarrow H^{p+q}(X^{\an},\mathbb C)
\]
degenerates at $E_1$, and its filtration is Deligne's Hodge filtration
\cite{Del74,DB81}.
In particular,
\begin{equation}\label{eqDBHdR}
\mathbb H^q(X,\underline\Omega^p_{X/\mathbb C})
 \simeq\Gr_F^pH^{p+q}(X,\mathbb C).
\end{equation}
We will often omit the superscript $\an$ from the notation.

The rational version of the increasing filtration is obtained from
$Y\mapsto\tau^{\leq s}\Gamma(Y^{\an},\mathbb Q)$ on smooth proper complex schemes by pcdh-descent.
Under Betti--de Rham comparison, its complexification is the increasing filtration $W^{\mathrm{dec}}$ on $\Gamma\underline\Omega^\bullet_{-/\bbC}$. Setting
\[
W_sH^m(X,\mathbb Q):=
\im\bigl(
 H^m(W^{\mathrm{dec}}_s\Gamma(X,\mathbb Q))
 \longrightarrow H^m(X,\mathbb Q)\bigr)
\]
recovers Deligne's weight filtration. Together with $F$, it endows
$H^m(X,\mathbb Q)$ with its mixed Hodge structure.

Define
\[
\underline h^{p,q}(X):=
 \dim_{\mathbb C}\mathbb H^q(X,\underline\Omega^p_{X/\mathbb C}),
\qquad
h_n^{a,b}(X):=
 \dim_{\mathbb C}\Gr_F^a
 \Gr_{a+b}^WH^n(X,\mathbb C).
\]
The induced weight filtration on $\Gr_F^pH^n$ is
\[
W_s\Gr_F^pH^n
 =\frac{(W_sH^n\cap F^pH^n)+F^{p+1}H^n}{F^{p+1}H^n}.
\]
Taking $\Gr_F$ and $\Gr^W$ commutes. Thus
\eqref{eqDBHdR} and the weight bounds for proper varieties give
\begin{equation}\label{eqHodgepiecedim}
\underline h^{p,q}(X)=\sum_{b=0}^{q}h_{p+q}^{p,b}(X).
\end{equation}
Here $h_n^{a,b}=0$ unless $a,b\geq0$ and $a+b\leq n$, and complex conjugation gives Hodge symmetry $h_n^{a,b}=h_n^{b,a}$.

\medskip

Definition~\ref{DBcomplex} suggests that the dimension numbers above might be definable for any field of characteristic zero. Now we show that it is indeed possible. Put
\[
H^n_{\DB}(X/k):=
 \mathbb H^n(X,\underline\Omega^\bullet_{X/k}),
\]
and give this vector space the filtrations
\[
\begin{aligned}
W_sH^n_{\DB}(X/k)
&:=\im\bigl(H^n(W^{\mathrm{dec}}_s
     \Gamma\underline\Omega^\bullet_{-/k})
     \to H^n_{\DB}(X/k)\bigr),\\
F^pH^n_{\DB}(X/k)
&:=\im\bigl(H^n(F^p\Gamma\underline\Omega^\bullet_{-/k})
     \to H^n_{\DB}(X/k)\bigr).
\end{aligned}
\]

\begin{prop}\label{DBHodgecharzero}
Let $X/k$ be a proper scheme over a field of characteristic zero.
The bifiltered complex $\Gamma\underline\Omega^\bullet_{-/k}$ and the induced filtrations on
$H^n_{\DB}(X/k)$ are functorial in $X$ and commute with extension of the ground field. The spectral sequence
\[
E_1^{p,q}=\mathbb H^q(X,\underline\Omega^p_{X/k})
 \Longrightarrow H^{p+q}_{\DB}(X/k)
\]
degenerates at $E_1$. Moreover,
\[
H^n(W^{\mathrm{dec}}_s\Gamma\underline\Omega^\bullet_{-/k})\longrightarrow H^n_{\DB}(X/k),
\qquad
H^n(F^p\Gamma\underline\Omega^\bullet_{-/k})\longrightarrow H^n_{\DB}(X/k)
\]
are injective.

\end{prop}

\begin{defi}\label{numbersdefi}

We define
\[
\underline h^{p,q}(X):=
 \dim_k\mathbb H^q(X,\underline\Omega^p_{X/k}),
\qquad
h_n^{a,b}(X):=
 \dim_k\Gr_F^a\Gr_{a+b}^W
 H^n_{\DB}(X/k).
\]

We call them \emph{Du Bois numbers} and \emph{mixed Hodge numbers} respectively.
    
\end{defi}

\begin{prop}\label{numbersoverk}

These are finite nonnegative integers, unchanged by extension of $k$, and they satisfy
\[
\underline h^{p,q}(X)=\sum_{b=0}^{q}h_{p+q}^{p,b}(X),
\qquad h_n^{a,b}(X)=h_n^{b,a}(X),
\]
as well as
\[
h_n^{a,b}(X)=0\quad\text{unless}\quad
a,b\geq0\ \text{and}\ a+b\leq n.
\]
\end{prop}

\begin{proof}[Proof of Propositions~\ref{DBHodgecharzero}
and~\ref{numbersoverk}]

Choose a bounded cubical hyperresolution by smooth proper schemes.
For every field extension $k\subset k'$, proper coherent base change for the finitely many complexes $\Omega^{\geq p}_{X_\alpha/k}$ gives
\[
\Gamma(X_\alpha,\Omega^{\geq p}_{X_\alpha/k})\otimes_k k'
 \simeq
\Gamma((X_\alpha)_{k'},
        \Omega^{\geq p}_{(X_\alpha)_{k'}/k'}).
\]
Extension of scalars is exact, so it commutes with canonical truncations, finite limits, cohomology, and images. Base change for the bifiltered complex and its cohomological filtrations follows from \eqref{eqDBbifilteredhyperresolution}. Finiteness of the Du Bois and mixed Hodge numbers follows from \eqref{eqDBbifilteredhyperresolution} and proper coherent finiteness.
Functoriality follows from the descent construction.

The scheme and its finite hyperresolution by smooth proper varieties
lift to a subfield $k_0\subset k$ of finite transcendence degree over $\bbQ$; enlarging $k_0$ comes with the base-changed hyperresolution. Choose an embedding
$k_0\hookrightarrow\mathbb C$. Degeneration, the injectivity statements of Proposition~\ref{DBHodgecharzero}, and the numerical equalities of Proposition~\ref{numbersoverk} hold after extension of scalars to $\bbC$. Now, faithful flatness
lifts the vanishing of the spectral-sequence differentials and the injectivity statements to
$k_0$; base change then gives them over $k$.

Finally, for $n=p+q$, degeneration and the interchange of associated graded pieces give
\[
\begin{aligned}
\underline h^{p,q}(X)
 &=\dim_k\Gr_F^pH^n_{\DB}(X/k)\\
 &=\sum_s\dim_k\Gr_F^p
                  \Gr_s^WH^n_{\DB}(X/k)\\
 &=\sum_{b=0}^{q}h_{p+q}^{p,b}(X),
\end{aligned}
\]
where the last range follows from the weight bounds  established above.
\end{proof}

\begin{remark}
We think that the above discussion establishes precisely which part of the mixed Hodge theory for (individual) proper varieties is present over an arbitrary characteristic zero field. The semilinear involution coming from complex conjugation and the conjugate filtration on cohomology seem to be complex numbers phenomena. Of course, choosing an embedding $\rho:k_0\subset \bbC$ allows one to lift them to a finite transcendence degree extension $k_0\subset k_1\subset\bbC$, but this lift seems to depend on the choice of an embedding $\rho$.
\end{remark}

\section{The Du Bois and $\bbA^1$-invariant topological cyclic homology}

In this section, we introduce Du Bois topological cyclic homology, together with its multiplicative and $\TC$-module structures, using Efimov's recollement. These invariants are closely related to algebraic and homotopy K-theory and to topological cyclic homology.
We arrange the definitions so that all the maps and relevant structures are already present at the level of endofunctors of $\Motloc_k$.

Like topological cyclic homology, its Du Bois counterpart comes in three flavors: Efimov's, genuine, and Nikolaus--Scholze's.
We show that the three theories agree on weakly geometric motives and are truncating.

The definition of $\tc_{\DB}$, the endofunctor of $\Motloc_k$
corresponding to Du Bois topological cyclic homology, suggests a relation to $\bbA^1$-homotopification. In the second subsection, we study the $\bbA^1$-homotopifications of $\tc$ and $\tc_{\DB}$.
This justifies our choice of construction of $\tc_{\DB}$ and demonstrates both the non-finitarity of the new theories and the subtleties of $\bbA^1$-homotopification for non-finitary invariants.
It also provides lifts of the invariants considered in Section~4.2 to the level of
endofunctors.

\subsection{Du Bois topological cyclic homology}

The two recollement triangles \eqref{recol1} and \eqref{recol2} fit into the following octahedron diagram of exact triangles:
\begin{equation}\label{eqoctahedronendofunctors}
\begin{tikzcd}[column sep=large,row sep=large]
 & \dk \arrow[r,"\id"] \arrow[d]
 & \dk \arrow[d] \\
\kinv \arrow[r] \arrow[d,"\id"']
 & \id \arrow[r] \arrow[d]
 & \tc \arrow[d] \\
\kinv \arrow[r]
 & \kh \arrow[r]
 & \kh\,\tc .
\end{tikzcd}
\end{equation}

\begin{defi}\label{deftcdb}
We denote the lower-right corner by
\[
\tc_{\DB}:=\kh\circ \tc=i_*i^*j_*j^*
\]
and define \emph{Du Bois (Efimov's) topological cyclic homology} by
\[
\TC^{\mathrm{E}}_{\DB}(T):=\Gamma\,\tc_{\DB}(\Uloc(T)),\qquad T\in\Catperf_k.
\]
\end{defi}

The octahedron diagram above follows, for instance, from exactness of the right column. Indeed, the right column can be obtained from the exact triple~\eqref{recol1} (middle column) by applying it to $\tc$. Applying $\kh$ to \eqref{recol2} gives the bottom row.

For convenience, we include the multiplication table. Two semi-orthogonal decompositions in the recollement\eqref{recollement} imply the following composition rules (see also Remark~\ref{dualDB}):
\[
\boxed{
\begin{array}{c|cccc}
\circ & \dk & \kh & \kinv & \tc\\
\hline
\dk
& \dk
& 0
& 0
& \dk
\\[2mm]
\kh
& 0
& \kh
& \kinv
& \tc_{\DB}
\\[2mm]
\kinv
& \tc_{\DB}[-1]
& \kh
& \kinv
& 0
\\[2mm]
\tc
& \tc
& 0
& 0
& \tc
\end{array}}
\]
Here the entry in row $F$ and column $G$ denotes the composition
$F\circ G$.

After evaluating the right column and the bottom row of \eqref{eqoctahedronendofunctors} on $T$ and applying $\Gamma$, we obtain
\begin{equation}\label{eqDBtwotrianglesglobal}
\DK(T)\longrightarrow \TC^{\mathrm{E}}(T)\longrightarrow \TC^{\mathrm{E}}_{\DB}(T)
\end{equation}
and
\[
\K^{\mathrm{inv},\mathrm{E}}(T)\longrightarrow \KH(T)\longrightarrow \TC^{\mathrm{E}}_{\DB}(T).
\]

Here are some corollaries of the above picture.

\begin{prop}\label{multiplicativeoctahedron}
The lower-right square in \eqref{eqoctahedronendofunctors} is a square of commutative algebra objects in the Day-convolution symmetric monoidal category of endofunctors of $\Motloc_k$.  The right column is an exact triangle of
$\tc$-modules and the bottom row is an exact triangle of $\kh$-modules.
\end{prop}

\begin{proof}
The idempotent algebra $A$ and the complementary coidempotent object $J$ carry their canonical commutative algebra and cocommutative coalgebra structures, respectively. They induce canonical lax symmetric monoidal structures on
$\kh=A\otimes-$ and $\tc=\underline{\Hom}_{\Motloc_k}(J,-)$.
For the latter, the
structure map is
\[
\underline{\Hom}_{\Motloc_k}(J,M)\otimes
\underline{\Hom}_{\Motloc_k}(J,N)
\longrightarrow
\underline{\Hom}_{\Motloc_k}(J\otimes J,M\otimes N)
\longrightarrow
\underline{\Hom}_{\Motloc_k}(J,M\otimes N).
\]
The unit of $A$ and the counit of $J$ induce monoidal natural transformations $\id\to \kh$ and $\id\to \tc$.
Together with their composites, these give the indicated square of lax symmetric monoidal functors.

Taking the fiber of $\tc\to \kh\,\tc$ in $\tc$-modules and the fiber of $\kh\to \kh\,\tc$ in $\kh$-modules gives the asserted exact triangles: after forgetting the module structures, the fibers are
$\dk\,\tc\simeq \dk$ and
$\kh\,\kinv\simeq \kinv$, respectively.
\end{proof}

\begin{corollary}\label{multiplicativeinvariants}
The diagram
\[
\begin{tikzcd}
\K\arrow[r]\arrow[d]&\TC^{\mathrm{E}}\arrow[d]\\
\KH\arrow[r]&\TC^{\mathrm{E}}_{\DB}
\end{tikzcd}
\]
is a cartesian square of lax symmetric monoidal localizing invariants. In particular, evaluation on a symmetric monoidal $k$-linear category gives a cartesian square of commutative ring spectra. For every $M\in\Motloc_k$,
\[
\DK(M)\longrightarrow \TC^{\mathrm{E}}(M)\longrightarrow \TC^{\mathrm{E}}_{\DB}(M),\qquad
\K^{\mathrm{inv},\mathrm{E}}(M)\longrightarrow \KH(M)\longrightarrow \TC^{\mathrm{E}}_{\DB}(M)
\]
are exact triangles over $\TC^{\mathrm{E}}(k)$ and $\KH(k)$, respectively.
For connective $k$ the coefficient ring $\TC^{\mathrm{E}}(k)$ agrees with $\TC(k)$.
\end{corollary}

\begin{remark}
The statement of Corollary~\ref{multiplicativeinvariants} has  been independently obtained by Thomas~Nikolaus and Maxim~Ramzi. 
\end{remark}

\begin{remark}
The author first learned about the $\TC(k)$-module structure on $\DK(T)$ from a conversation with A.~Efimov; see also Theorem 9.5 in \cite{Efi25b} for a similar statement.
\end{remark}

\begin{remark}
Evaluating the square of Corollary~\ref{multiplicativeinvariants}
on $\Perf(X)$, for a qcqs scheme $X$ over a prime field $F$, recovers the cartesian square considered in \cite{EM26}. There it is refined to a cartesian square of filtered commutative ring spectra
\[
\begin{tikzcd}[column sep=large]
\Fil^{\bullet}_{\mathrm{mot}}\K(X)
  \arrow[r]\arrow[d]
&
\Fil^{\bullet}\TC(X)\arrow[d]
\\
\Fil^{\bullet}_{\mathrm A}\KH(X)\arrow[r]
&
(\Fil^{\bullet}\TC)_{\cdh}(X),
\end{tikzcd}
\]
where the right-hand filtration is the HKR filtration on
$\TC(X)\simeq \HC^-(X/\bbQ)$ if $F=\bbQ$, and the BMS filtration if $F=\bbF_p$; see
\cite[Definitions 4.9 and 4.28]{EM26}.

Taking fibers in filtered spectra therefore equips $\DK(X)$
and $\K^{\mathrm{inv}}(X)$ with compatible filtrations and refines the octahedron to an octahedron of filtered spectra.
The right column and bottom row are moreover triangles of filtered modules over the filtered $\TC(X)$ and $\KH(X)$, respectively.

For an explicit description of the graded pieces, we use the notation $\bbZ(j)^{\TC}$ of Tess~Bouis:
\[
\bbZ(j)^{\TC}(X)\simeq
\begin{cases}
\Gamma\bigl(X,\widehat{L\Omega}_{-/\bbQ}^{\ge j}\bigr),
  & F=\bbQ,\\
\bbZ_p(j)^{\mathrm{syn}}(X),
  & F=\bbF_p.
\end{cases}
\]
Thus
$\gr^j\TC(X)\simeq\bbZ(j)^{\TC}(X)[2j]$.
Writing $\bbZ(j)^{\mathrm{mot}}$ and $\bbZ(j)^{\mathrm A}$
for the motivic complexes of \cite{EM26}, for every $j\in\bbZ$ we obtain
\[
\begin{aligned}
\gr^j\DK(X)
&\simeq
\fib\bigl(
\bbZ(j)^{\mathrm{mot}}(X)\longrightarrow
\bbZ(j)^{\mathrm A}(X)
\bigr)[2j]\\
&\simeq
\fib\bigl(
\bbZ(j)^{\TC}(X)\longrightarrow
\bbZ(j)_{\cdh}^{\TC}(X)
\bigr)[2j],
\\[1mm]
\gr^j\K^{\mathrm{inv}}(X)
&\simeq
\fib\bigl(
\bbZ(j)^{\mathrm{mot}}(X)\longrightarrow
\bbZ(j)^{\TC}(X)
\bigr)[2j]\\
&\simeq
\fib\bigl(
\bbZ(j)^{\mathrm A}(X)\longrightarrow
\bbZ(j)_{\cdh}^{\TC}(X)
\bigr)[2j].
\end{aligned}
\]
Here the motivic and syntomic complexes vanish in negative weights, and
$\widehat{L\Omega}^{\ge j}=\widehat{L\Omega}$ for $j\le0$.

\end{remark}

\begin{remark}\label{dualDB}
There is a canonical equivalence of endofunctors
\begin{equation}\label{eqdualdescription}
\Sigma\bigl(\kinv\circ \dk\bigr)\simeq \tc_{\DB}.
\end{equation}
or, equivalently,
\[
\Sigma\,i_*i^!j_!j^*\simeq i_*i^*j_*j^*,
\]
which can be seen, for instance, by applying $\kinv=i_*i^!$ to \eqref{recol1}.
$\Omega \tc_{\DB}$ lives in a "dual" square
\begin{equation}\label{eqdualsquare}
\begin{tikzcd}[column sep=large,row sep=large]
\Omega \tc_{\DB} \arrow[r] \arrow[d]
  & \dk \arrow[d] \\
\kinv \arrow[r]
  & \id.
\end{tikzcd}
\end{equation}
\end{remark}

We finish by defining the genuine and Nikolaus--Scholze versions.
By \eqref{cyclotomicrealizations} and Corollary~\ref{multiplicativeinvariants}, the maps
$\DK\to \TC^{\mathrm{E}}\to \TC^{\gen}\to \TC^{\NS}$ are $\TC^{\mathrm{E}}(k)$-linear.

\begin{defi}
For $?\in\{\mathrm{E},\gen,\NS\}$ put
\[
\TC^?_{\DB}:=\cof(\DK\longrightarrow \TC^?).
\]
\end{defi}

These are $\Mod_{\TC^{\mathrm{E}}(k)}$-valued invariants and the natural maps between them are $\TC^{\mathrm{E}}(k)$-linear.

\begin{prop}\label{TCDBtruncating}
The three theories agree on $\Motlocwg_k$ and are truncating.
\end{prop}

\begin{proof}
Agreement follows by taking cofibers in
\eqref{cyclotomicrealizations}. For a connective
$E_1$-$k$-algebra $R$, the bottom row of
\eqref{eqoctahedronendofunctors} becomes
\[
\K^{\mathrm{inv}}(R)\longrightarrow \KH(R)
\longrightarrow \TC^{\mathrm{E}}_{\DB}(R).
\]

The first term is truncating by
\cite[Theorem 7.0.0.2]{DGM13}; the second is truncating by the proof of \cite[Proposition 3.14]{LT19}, which adapts straightforwardly to an arbitrary base $k$.
Hence $\TC^{\mathrm{E}}_{\DB}$ is truncating. Agreement on connective algebras gives the assertion for the other two versions.

\end{proof}

These invariants usually are not finitary. For $?\in\{\gen, \NS\}$, finitarity of $\DK$ together with non-finitarity of $\TC^?$ implies non-finitarity of $\TC^?_{\DB}$. For $\TC^{\mathrm{E}}$ and $\TC^{\mathrm{E}}_{\DB}$, see Remark~\ref{nonfinitarity}.

\subsection{tc and $\bbA^1$-invariance}

\begin{defi}
Consider the endofunctors of $\Motloc_k$
\begin{align}
\tc^{\mathrm{ac}}(-)
&:=|[n]\longmapsto
\underline{\Hom}_{\Motloc_k}
\bigl(J,-\otimes\widetilde\Delta_k^n\bigr)|,
\label{eqtccycdef}\\
\tc^{\bbA^1}(-)
&:=|[n]\longmapsto
\underline{\Hom}_{\Motloc_k}
\bigl(J,-\otimes\Delta_k^n\bigr)|.
\label{eqtcA1def}
\end{align}

We call the composite functors

\[
\TC^{\mathrm{E}, \bbA^1}(-):=\Gamma\,\tc^{\bbA^1}(\Uloc(-)),
\qquad
\TC^{\mathrm{E}, \mathrm{ac}}(-):=\Gamma\,\tc^{\mathrm{ac}}(\Uloc(-)).
\]

\emph{$\bbA^1$-invariant} and \emph{$\bbA^1$-acyclic Efimov's topological cyclic homology}.

\end{defi}

Note that this is compatible with \ref{A1ization}.

We have an exact triangle of endofunctors

\begin{equation}\label{actcA1}
\tc^{\mathrm{ac}}\longrightarrow \tc\longrightarrow \tc^{\bbA^1}.
\end{equation}

\begin{remark}\label{nonfinitarity}
For non-finitary invariants, \eqref{finitaryA1ization}
need not hold. In fact, we have already seen that

\[
\tc\circ \kh=j_*j^*i_*i^*=0,
\qquad
\kh\circ \tc=\tc_{\DB}.
\]
For example, over $k=\bbQ$, Efimov's comparison and
\ref{char0cyclictriangle} (Geller-Weibel) give
\[
\Gamma\,\tc^{\bbA^1}(\mathbbm1)\simeq \HP(\bbQ/\bbQ),
\qquad
\Gamma\,\tc_{\DB}(\mathbbm1)\simeq \HC^-(\bbQ/\bbQ).
\]
Here $\DK(\bbQ)=0$.
Thus $\tc^{\bbA^1}$ agrees with neither $\tc\circ \kh$
nor $\kh\circ \tc$ in general. In particular, over $\bbQ$,
$\TC^{\mathrm{E}}$ is not finitary: otherwise
\eqref{finitaryA1ization} would identify its
$\bbA^1$-homotopification with $\Gamma\,\tc\circ \kh=0$.
Consequently $J$ is not compact, since
$\Gamma\,\tc=\Hom_{\Motloc_k}(J,-)$; also $\tc$, $\TC^{\mathrm{E}}$, and
$\TC^{\mathrm{E}}_{\DB}$ are not finitary.

Note also that, in general, neither $\tc^{\bbA^1}$ nor
$\TC^{\mathrm{E},\bbA^1}$ factors through
$i^*:\Motloc_k\ra\MotlocA$.
Indeed, such a factorization would imply
$\Gamma\,\tc^{\bbA^1}(J)=0$, whereas

\begin{align}
    \tc^{\bbA^1}(J)
&=|[n]\longmapsto
\underline{\Hom}_{\Motloc_k}
\bigl(J,J\otimes\Delta_k^n\bigr)|\simeq |[n]\longmapsto
(\tc \dk)(\Delta_k^n)|\simeq |[n]\longmapsto
\tc(\Delta_k^n)|\simeq \tc^{\bbA^1}(k).
\end{align}
Since $\TC^{\mathrm{E}, \bbA^1}$ is lax symmetric monoidal, the factorization exists if and only if $\TC^{\mathrm{E}, \bbA^1}=0$, which is not always the case, as the example over $\bbQ$
above shows.

\end{remark}

Starting from the next section, we will only consider the functors defined above in the characteristic 0 case. Here we just recall the result showing how different the situation can be over a finite field: over $\bbF_p$ the factorization discussed above does happen.

\begin{prop}\label{FpA1} (Elmanto, Theorem 3.1.1 \cite{Elm20})
Let $k=\bbF_p$.  Then
\[
\Gamma\,\tc^{\bbA^1}\simeq0,
\]
Hence the canonical map
\[
\Gamma\,\tc^{\mathrm{ac}}\longrightarrow\Gamma\,\tc
\]
is an equivalence.
\end{prop}

\begin{proof}

Since $\Gamma\,\tc^{\bbA^1}$ is lax-symmetric monoidal, it is enough to check the vanishing on the unit.
By Efimov's corepresentability theorem for TC and \cite{Elm20}, we have
\[
\Gamma\,\tc^{\bbA^1}(k)\simeq \TC^{\mathrm{E},\bbA^1}(k)\simeq \TC^{\NS,\bbA^1}(k)=0.
\]
\end{proof}

We finish by introducing a ``linear order'' on the
    $\tc$-module endofunctors we consider.

Applying $\kh$ to \eqref{actcA1} gives an exact triangle

\[
(\tc_{\DB})^{\mathrm{ac}}\simeq \tc_{\DB}^{\mathrm{ac}}:=\kh\,\tc^{\mathrm{ac}}
\longrightarrow \tc_{\DB}\longrightarrow
\kh\,\tc^{\bbA^1}=:\tc_{\DB}^{\bbA^1}\simeq (\tc_{\DB})^{\bbA^1}.
\]
The outer isomorphisms follow from finitarity of $\kh$.

The morphism $\eta:\id\to \tc$, induced by $J\to\mathbbm1$, gives a commutative square
\[
\begin{tikzcd}
\id\arrow[r,"\eta"]\arrow[d]&\tc\arrow[d]\\
\kh=|{-}\otimes\Delta_k^\bullet|
\arrow[r,"|\eta|"']&\tc^{\bbA^1}=|\tc({-}\otimes\Delta_k^\bullet)|.
\end{tikzcd}
\]
The map on fibers is explicitly
\begin{equation}\label{dktotcac}
\dk(M)\simeq|M\otimes\widetilde\Delta_k^\bullet|
\xrightarrow{|\eta|}
|\tc(M\otimes\widetilde\Delta_k^\bullet)|=\tc^{\mathrm{ac}}(M).
\end{equation}
Using $\dk\,\tc\simeq \dk$ and preservation of colimits by $\dk$, we obtain
\[
\dk\,\tc^{\mathrm{ac}}\simeq \dk\,\dk\simeq \dk,
\qquad \dk\,\tc^{\bbA^1}\simeq \dk\,\kh=0.
\]
Thus $\tc^{\bbA^1}$ is $\kh$-local. The map
$\tc\to \tc^{\bbA^1}$ factors through $\tc_{\DB}$, and there is a commutative diagram with exact rows and columns:
\begin{equation}\label{acDBdiagram}
\begin{tikzcd}[column sep=large,row sep=large]
\dk\arrow[r]\arrow[d,equal]&\tc^{\mathrm{ac}}\arrow[r]\arrow[d]
  &\tc_{\DB}^{\mathrm{ac}}\arrow[d]\\
\dk\arrow[r]&\tc\arrow[r]\arrow[d]&\tc_{\DB}\arrow[d]\\
&\tc^{\bbA^1}\arrow[r,equal]&\tc^{\bbA^1}.
\end{tikzcd}
\end{equation}

\begin{corollary}\label{acyclicDBcorollary}
There is a canonical diagram
\[
\dk\longrightarrow \tc^{\mathrm{ac}}\longrightarrow \tc
\longrightarrow \tc_{\DB}\longrightarrow \tc^{\bbA^1}
\]
where composites of length 3 are canonically zero.
Moreover,
\[
\tc_{\DB}^{\bbA^1}\simeq \tc^{\bbA^1},\qquad
\fib(\tc^{\mathrm{ac}}\to \tc_{\DB}^{\mathrm{ac}})
\simeq \dk\simeq\fib(\tc\to \tc_{\DB}).
\]
\end{corollary}

\section{Du Bois localizing invariants in characteristic 0}

In this section, we build on the constructions of the previous section to define the main invariant of the paper -- Du Bois Hochschild homology. In the first subsection, we construct it together with its circle action in two flavours: Efimov's and classical.
We give a first application of this invariant.
We also define the corresponding two flavours of Du Bois negative cyclic homology and study the properties of these invariants and their relations.

In the second subsection, we introduce the
(``positive'') cyclic and periodic Du Bois theories and study their properties. In particular, we prove a version of the Geller--Weibel theorem, identifying periodic cyclic homology with the
$\bbA^1$-homotopification of negative cyclic homology for the theories considered here.
Over $\bbC$, we construct a Du Bois counterpart of Blanc's Deligne cohomology.

\medskip

From now on, $k$ is a connective commutative
$\bbQ$-algebra and $?\in\{\mathrm{E},\cl\}$;
see below.

\subsection{Du Bois Hochschild and negative cyclic homology}

We denote by
\[
\HC^{-,\cl}(-/k):
\Catperf_k\longrightarrow\Mod_{k^{hS^1}}
\]
the classical negative cyclic homology functor on $k$-categories. Similarly, we usually include the superscript $\cl$
in Hochschild, cyclic, and periodic cyclic homology to distinguish the classical theories from their Efimov versions. When the two versions agree, we omit this superscript.

To use uniform notation, we also write
\[
\TC^{\cl}:=\TC^{\NS},
\qquad
\TC^{\cl}_{\DB}:=\TC^{\NS}_{\DB}.
\]

For $\bbQ$-algebras, one has the canonical identifications

\[
\TC^{\gen}(T)\simeq \TC^{\NS}(T)
\simeq \HC^{-, \cl}(T/\bbQ),
\qquad
\K^{\mathrm{inv}}(T)\simeq \K^{\mathrm{inf}}(T),
\]
where $\K^{\mathrm{inf}}$ is the fiber of the Chern character
$\K\to \HC^{-,\cl}(-/\bbQ)$, due to the vanishing of finite-group Tate constructions on rational spectra. On weakly geometric motives, Observation~\ref{obscomparisonclassical}
further identifies
\[
\HC^{-, \mathrm{E}}(T/\bbQ)\simeq\TC^{\mathrm{E}}(T)
\simeq \HC^{-,\cl}(T/\bbQ),
\qquad
\K^{\mathrm{inv},\mathrm{E}}(T)\simeq \K^{\mathrm{inf}}(T).
\]

By Corollary~\ref{multiplicativeinvariants}, we have exact triangles of $\HC^-(k/\bbQ)$-modules
\begin{equation}\label{eqDKHCTCDBmodules}
\DK(T)\longrightarrow \TC^?(T)\longrightarrow \TC^?_{\DB}(T).
\end{equation}

Consider the following maps of lax symmetric monoidal functors:
\begin{equation}\label{TCEtoHH}
\TC^?\longrightarrow \HC^{-,?}(-/\bbQ)
\longrightarrow \HH^?(-/\bbQ)\longrightarrow \HH^?(-/k).
\end{equation}

The following observation motivates the coming definitions.

\begin{lemma}\label{equivariantscalarreduction}
The map \eqref{TCEtoHH} provides an equivalence of lax symmetric monoidal functors
\begin{equation}\label{eqequivariantscalarreduction}
\TC^{\cl}\otimes_{\HC^-(k/\bbQ)}\HH(k/k)
\simeq \HH^{\cl}(-/k):
\Catperf_k\longrightarrow\Mod_k^{BS^1}.
\end{equation}
\end{lemma}

\begin{proof}
We have
\begin{align*}
\TC^{\cl}\otimes_{\HC^-(k/\bbQ)}\HH(k/k)
&\simeq
\HC^{-, \cl}(-/\bbQ)\otimes_{\HC^-(k/\bbQ)}\HH(k/\bbQ)
                    \otimes_{\HH(k/\bbQ)}\HH(k/k)\\
&\simeq \HH^{\cl}(-/\bbQ)\otimes_{\HH(k/\bbQ)}\HH(k/k)\\
&\simeq \HH^{\cl}(-/k).
\end{align*}
The first two equivalences follow from the discussion above, and the last is base change for Hochschild homology. The intermediate tensor product carries its natural circle action, which is compatible with the base change.
\end{proof}

\begin{defi}\label{DBcyclicdefinitions}
Using \eqref{TCEtoHH}, we define \emph{classical Du Bois Hochschild homology} and \emph{Efimov--Du Bois Hochschild homology}
by
\[
\HH^?_{\DB}:=
\TC^?_{\DB}\otimes_{\HC^-(k/\bbQ)}\HH(k/k):
\Catperf_k\longrightarrow\Mod_k^{BS^1}.
\]
We define the corresponding negative cyclic homology theories using homotopy fixed points:
\[
\HC^{-,?}_{\DB}:=(\HH^?_{\DB})^{hS^1}:
\Catperf_k\longrightarrow\Mod_{\HC^-(k/k)}.
\]

When $k$ is a field and dimensions of $\HH^{\cl}_\DB(T/k)$ and $\HH^{\mathrm{E}}_\DB(T/k)$ as $k$-vector spaces agree and are finite, we will refer to them as \emph{Du Bois-Hochschild numbers} of the $k$-category $T$.
\end{defi}

The coefficient algebra is
\[
\HC^-(k/k)\simeq \HH(k/k)^{hS^1}\simeq k[[u]],
\qquad |u|=-2.
\]

Put
\[
\DHH(-/k):=\DK\otimes_{\HC^-(k/\bbQ)}\HH(k/k).
\]
Taking homotopy fixed points gives
\begin{equation}\label{DHCcompletion}
\begin{aligned}
\DHC^-(-/k)
&:=(\DHH(-/k))^{hS^1}\\
&\simeq
(
\DK\otimes_{\HC^-(k/\bbQ)}\HC^-(k/k)
)^\wedge_u.
\end{aligned}
\end{equation}
Here $(-)^\wedge_u$ denotes $u$-adic (derived) completion. In particular,
\[
\DHC^-(-/k)/u\simeq \DHH(-/k).
\]

We therefore have the following commutative diagrams of localizing invariants on $\Catperf_k$, with exact columns:
\begin{equation}\label{negativecycliccorrectiondiagram}
\begin{gathered}
\begin{tikzcd}[column sep=large,row sep=large]
\DK\arrow[r]\arrow[d]
 &\DHC^-(-/k)\arrow[r,equal]\arrow[d]
 &\DHC^-(-/k)\arrow[d]\\
\K\arrow[r]\arrow[d]
 &\HC^{-,\mathrm{E}}(-/k)\arrow[r]\arrow[d]
 &\HC^{-,\cl}(-/k)\arrow[d]\\
\KH\arrow[r]
 &\HC^{-,\mathrm{E}}_{\DB}(-/k)\arrow[r]
 &\HC^{-,\cl}_{\DB}(-/k).
\end{tikzcd}
\qquad
\begin{tikzcd}[column sep=large,row sep=large]
\DHH(-/k)\arrow[r,equal]\arrow[d]
 &\DHH(-/k)\arrow[d]\\
\HH^{\mathrm{E}}(-/k)\arrow[r]\arrow[d]
 &\HH^{\cl}(-/k)\arrow[d]\\
\HH^{\mathrm{E}}_{\DB}(-/k)\arrow[r]
 &\HH^{\cl}_{\DB}(-/k).
\end{tikzcd}
\end{gathered}
\end{equation}
Thus
\[
\begin{aligned}
\DHC^-(-/k)
&\simeq
\fib\bigl(\HC^{-,\cl}(-/k)\longrightarrow \HC^{-,\cl}_{\DB}(-/k)\bigr)\\
&\simeq
\fib\bigl(\HC^{-,\mathrm{E}}(-/k)\longrightarrow \HC^{-,\mathrm{E}}_{\DB}(-/k)\bigr).
\end{aligned}
\]
In particular, it does not depend on the version of negative cyclic homology.

\begin{remark}
When $\bbQ\subset k\subset\overline{\bbQ}$ is an algebraic field extension,
\[
\HH(k/\bbQ)\simeq \HH(k/k)\simeq k\simeq \THH(k).
\]
Thus the corresponding absolute and relative Hochschild theories agree, as do their cyclic variants.
\end{remark}

\begin{remark}
For an arbitrary $M\in\Motloc_\bbQ$, the comparison
\[
\TC^?_{\DB}(M)\longrightarrow
\HC^{-,?}_{\DB}(M/\bbQ)
\]
need not be an equivalence. It is an equivalence if and only if
\[
\DK(M)\longrightarrow\DHC^-(M/\bbQ)
\simeq\DK(M)^\wedge_u
\]
is so (for the Efimov version, we use the derived
$u$-completeness of $\TC^{\mathrm E}(M)$, which follows from the argument in the proof of Theorem~\ref{char0GW}), which is not always the case.

We can use rigidity of $\Motloc_\bbQ$
to incarnate classical cyclic homology by a motive
$M_{\HC}$, so that there is an equivalence of functors
\[
\Gamma(M_{\HC}\otimes_\bbQ -)
\simeq\HC^{\cl}(-/\bbQ).
\]
Since $\HC^{\cl,\bbA^1}=0$ by Geller-Weibel (\cite{GW89}),
$\HC^{\cl}(J_\bbQ/\bbQ)\simeq\HC^{\cl}(\bbQ/\bbQ)$.
Consequently, we  have an equivalence of $\bbQ[[u]]$-modules
\[
\DK(M_{\HC}) \simeq\HC^{\cl}(J_\bbQ/\bbQ)\simeq\HC^{\cl}(\bbQ/\bbQ)
\simeq\Omega^2\cof\bigl(\HC^{-,\cl}(\bbQ/\bbQ)\longrightarrow \HP^{\cl}(\bbQ/\bbQ)\bigr),
\]
where the last equivalence is given by the Connes exact triangle \eqref{ConnesDB}.
Its derived $u$-completion is $\Omega\HC^{-,\cl}(\bbQ/\bbQ)$, so the unit of the derived completion is not an equivalence.
\end{remark}

\begin{prop}\label{HHDBfinitarytruncating}
$\HH^{\cl}_{\DB}$ is finitary. Both versions of $\HH_{\DB}$, and hence their negative cyclic, cyclic, and periodic cyclic theories, are truncating. The versions agree on $\Motlocwg_k$.
\end{prop}

The cyclic and periodic cyclic versions $\HC^?_{\DB}$ and
$\HP^?_{\DB}$ are defined in the usual way, using homotopy orbits and Tate invariants; see Definition~\ref{HCHPDBdefi}. Alternatively, they can be recovered from $\HC^{-,?}_{\DB}$ using $\bbA^1$-homotopification, see Theorem~\ref{char0GW}.

\begin{proof}
The triangle
\[
\DHH\longrightarrow \HH^{\cl}\longrightarrow \HH^{\cl}_{\DB}
\]
proves the first assertion: $\DK$ and $\HH^{\cl}$ are finitary, and scalar extension preserves colimits. Truncation follows by base change from Proposition~\ref{TCDBtruncating}; taking homotopy fixed points, homotopy orbits, or Tate invariants preserves the resulting equivalences. Agreement follows in the same way from agreement of the $\TC_{\DB}$ theories on weakly geometric motives.
\end{proof}

\medskip

We give a simple application of $\HH_{\DB}$, there will be more in the next sections.

\begin{defi}\label{defittg}
    We call a motive $M\in \Motloc_k$ a \emph{truncation fiber} if it is isomorphic to a motive of the form $\fib(\Uloc(\Perf(R))\ra\Uloc(\Perf(\pi_0R)))$ for a connective $k$-algebra $R$.
    
    We call a motive \emph{ttg} if it belongs to the smallest thick subcategory of $\Motloc_k$ containing motives of all smooth proper $k$-schemes and all the truncation fiber motives. We say that a category $T\in\Catperf_k$ is ttg if so is $\Uloc(T)$.
\end{defi}

\begin{obs}\label{tgfiniteness}
Assume $k$ is a field. If $M$ is ttg, then $\HH_{\DB}(M/k)$ is a perfect
$k$-module.
\end{obs}

\begin{proof}
$\HH_{\DB}$ is perfect on motives of smooth proper varieties and vanishes on truncation fibers. Perfectness is preserved by finite cofibers, shifts, and retracts.
\end{proof}

\begin{prop}\label{PerfZXnotttg}
Let $k$ be a field and $X$ a smooth proper $k$-scheme. Let $U\subset X$ be an affine open subscheme of positive dimension, and let
$Z=X - U$ be the closed complement. Then $\Perf_Z(X)$ is not ttg.
\end{prop}

\begin{proof}
The localization exact triple of Example~\ref{Nisnevich} gives
\[
\Uloc(\Perf_Z(X))\longrightarrow\Uloc(X)\longrightarrow\Uloc(U).
\]
Since $X$ and $U$ are smooth, $\DK$ vanishes on both, hence also on $\Perf_Z(X)$. Thus $\HH_{\DB}$ agrees with $\HH$ on all three terms.

The $k$-vector space $\HH_0(U)=H^0(U,\mathcal{O}_U)$ is infinite-dimensional. Thus, by Observation~\ref{tgfiniteness}, perfectness of $\HH(X/k)$ implies that $\Uloc(\Perf_ZX)$ is not ttg.

\end{proof}

\medskip

\subsection{Further versions of cyclic homology in characteristic 0 and Deligne cohomology}\label{seccyclicDeligne}

We can define versions of ("positive") cyclic and periodic cyclic homology in a similar fashion.
For the Efimov theories, we use the scalar-extension realization
\[
\HH^{\mathrm{E}}(-/k):=\TC^{\mathrm{E}}\otimes_{\HC^-(k/\bbQ)}\HH(k/k),
\qquad \HC^{-,\mathrm{E}}(-/k):=\HH^{\mathrm{E}}(-/k)^{hS^1},
\]
with its canonical circle action, as in Definition~\ref{DBcyclicdefinitions}. Its cyclic and periodic versions are obtained by homotopy orbits and Tate invariants.

\begin{defi}\label{HCHPDBdefi}
We define Du Bois and Efimov--Du Bois cyclic and periodic cyclic homology using homotopy orbits and Tate invariants:
\[
\HC^?_{\DB}:=(\HH^?_{\DB})_{hS^1}:
\Catperf_k\longrightarrow\Mod_k,
\]
\[
\HP^?_{\DB}:=(\HH^?_{\DB})^{tS^1}:
\Catperf_k\longrightarrow\Mod_{\HP(k/k)}.
\]
\end{defi}

Here $\HP(k/k)\simeq k((u))$, and cyclic homology retains its natural
$u$-action given by Connes periodicity. As usual, one has the Connes exact triangle
\begin{equation}\label{ConnesDB}
\Sigma \HC^?_\dagger\longrightarrow \HC^{-,?}_\dagger
\longrightarrow \HP^?_\dagger,
\end{equation}
where $\dagger$ is either empty or $\DB$.

We prove the following version of the Geller-Weibel theorem.

\begin{theorem}\label{char0GW}\label{char0cyclictriangle}
Let $?\in\{\mathrm{E},\cl\}$ and let $\dagger$ be either empty or $\DB$.
The $\bbA^1$-homotopifications of $\HH^?_\dagger(-/k)$ and
$\HC^?_\dagger(-/k)$ vanish, and $\HP^?_\dagger(-/k)$ is
$\bbA^1$-invariant. Consequently, the canonical map
\[
\HC^{-,?}_\dagger(-/k)\longrightarrow \HP^?_\dagger(-/k)
\]
is the $\bbA^1$-homotopification map, and
\[
(\HC^{-,?}_\dagger)^{\mathrm{ac}}\simeq\Sigma \HC^?_\dagger.
\]
\end{theorem}

\begin{remark}
    As one can see from the proof, the key nontrivial part is $\bbA^1$-invariance of new versions of periodic cyclic homology. This part uses Efimov's recollement once again and his result about corepresentability of TR of connective algebras; the rest essentially boils down to the classical Geller-Weibel theorem.
\end{remark}

\begin{proof}
The invariants $\HH^?$ are lax symmetric monoidal and agree with ordinary Hochschild homology on the connective algebras
$k[\Delta^n]$. Geller--Weibel's polynomial-simplex calculation
\cite{GW89} gives
\[
\HH^{\bbA^1}(k/k)=0.
\]
The homotopifications $\HH^{?, \bbA^1}$ therefore have zero unit and vanish.
The multiplicative comparison $\HH^{\mathrm{E}}\to \HH^{\mathrm{E}}_{\DB}$ makes the latter a module over the former; its homotopification vanishes as well.
The diagram \eqref{negativecycliccorrectiondiagram}
then gives the assertion for $\HH^{\cl}_{\DB}$.
Homotopy orbits commute with realizations, proving the assertion for $\HC^?_\dagger$. It remains to prove that $\HP^?_\dagger$ is $\bbA^1$-invariant.

We show this for $k=\bbQ$;
the general case follows by compatibilities with the base change built into the definitions. Put
\[
J_{\TR,\bbQ}:=\cof\bigl(\mathbbm1\longrightarrow\Uloc(\bbQ[t])\bigr).
\]
The
reduced affine-line motive is a counital coalgebra, and
\begin{equation}\label{JTRandJTC}
J_\bbQ\otimes J_{\TR,\bbQ}\simeq J_{\TR,\bbQ}\end{equation} (see \cite[sections 8 and 9.1]{Efi25b}); in particular, $J_{\TR,\bbQ}$ lies in the essential image of $j_!:\Motloccyc_\bbQ\ra \Motloc_\bbQ$ and has a natural action of $\End(J_\bbQ)\simeq \TC^{\mathrm{E}}(\bbQ)\simeq \bbQ[[u]]$ compatible with \eqref{JTRandJTC}. $\bbA^1$-invariance of $\HP^?_\dagger$ amounts to vanishing of $\HP^?_\dagger(M\otimes J_{\TR,\bbQ})$ for every $M\in \Motloc_\bbQ$. Similarly to \eqref{negativecycliccorrectiondiagram}, we have exact triangles
\[
(\DK^\wedge_u)[u^{-1}]
\simeq \DHC^-(-/\bbQ)[u^{-1}]
\longrightarrow \HP^?(-/\bbQ)
\longrightarrow \HP^?_{\DB}(-/\bbQ).
\]

Classical periodic cyclic homology is $\bbA^1$-invariant by \cite{GW89}, so
$\HP^{\cl}(M\otimes J_{\TR,\bbQ}/\bbQ)=0$. 
Thus it suffices to prove the vanishing of
$\bigl(\DK(M\otimes J_{\TR,\bbQ})^\wedge_u\bigr)[u^{-1}]$
and $\HP^{\mathrm{E}}(M\otimes J_{\TR,\bbQ}/\bbQ)$.

Recall that
\[
\bigl(\DK(M\otimes J_{\TR,\bbQ})^\wedge_u\bigr)[u^{-1}]
\simeq
\bigl(\Gamma(M\otimes J_{\TR,\bbQ}\otimes J_\bbQ)^\wedge_u\bigr)
[u^{-1}];
\]
\[
\HP^{\mathrm{E}}(M\otimes J_{\TR,\mathbb Q}/\mathbb Q)
\simeq
\bigl(\TC^{\mathrm{E}}(M\otimes J_{\TR,\mathbb Q})^\wedge_u\bigr)[u^{-1}]\simeq  \bigl(\Hom_{\Motloccyc_\bbQ}(J_\bbQ, M\otimes J_{\TR,\bbQ}\otimes J_\bbQ)^\wedge_u\bigr)[u^{-1}].
\]

By Corollary~\ref{multiplicativeinvariants} and by the proof of Proposition~\ref{multiplicativeoctahedron}, we know that the action of $u$ on the  invariants comes from the action of $\End(J_\bbQ)\simeq \TC^{\mathrm{E}}(\bbQ)\simeq \bbQ[[u]]$ on the corresponding motives. We claim that $u$ acts by zero on $J_{\TR,\bbQ}$ which implies the desired vanishings.

By \eqref{JTRandJTC}, the counit of $J_{\TR, \bbQ}$ factors uniquely through $J_\bbQ$, hence the action map
\[\End(J_\bbQ) \lra \End(J_{\TR, \bbQ})\]

factors through 
\[
\Hom_{\Motloc_\bbQ}(J_{\TR,\bbQ},J_\bbQ)\simeq 
\Hom_{\Motloc_\bbQ}(J_{\TR,\bbQ},\mathbbm 1) \simeq \TR(\bbQ)
\]
where the last equivalence is provided by
\cite[Theorem 9.1]{Efi25b}. $\TR(\bbQ)$ is known to have $\pi_{i\ne 0}=0$, hence the multiplication by $u$ is nullhomotopic.\end{proof}

\begin{remark}
Alternatively one could use vanishing of finite-group Tate constructions on rational spectra together with \cite[Theorem 9.1]{Efi25b} to see that $\End_{\Motloc_\bbQ}(J_{\TR,\bbQ}) \simeq \cof(\TR(\bbQ) \ra \TR(\bbQ[t]))$ has $\pi_{-2}=0$ (in fact, one can see that $\pi_{0,1}$ are the only nontrivial homotopy groups), hence it can have no nonzero $u$-action. This allows one to replace the action factoring in the last paragraph by an explicit computation.
\end{remark}

The vanishing of the action of $u$ on $J_{\TR,k}$ has another consequence.

\begin{prop}
For every $M\in\Motloc_k$, one has
\[
\DK(M)[u^{-1}]\simeq0.
\]
\end{prop}

\begin{proof}
\cite[Theorem 8.1(i)]{Efi25b} holds after replacing $\bbS$ by an arbitrary commutative base ring $k$: $j_!\Motloccyc_k$ is the localizing tensor ideal generated by $J_{\TR,k}$.
The proof of Theorem~\ref{char0GW} shows that $u$ acts by zero on $J_{\TR,k}$, and hence on
$N\otimes_kJ_{\TR,k}$ for every $N\in\Motloc_k$.
Inverting $u$ is exact, preserves colimits and commutes with $\Gamma$.
It therefore vanishes on the entire ideal
$j_!\Motloccyc_k$, hence
\[
\DK(M)[u^{-1}]
\simeq
\Gamma\bigl((M\otimes_kJ_k)[u^{-1}]\bigr)\simeq\Gamma(0)
\simeq0.
\]
\end{proof}

\medskip

We note also the following splitting statement (cf. Corollary~\ref{DKsplittingvarieties}).

\begin{prop}\label{DHCsplitting}
Let $M\in\Motloc_k$. Suppose that the canonical comparison
\[
\HP^?(M/k)\longrightarrow \HP^?_{\DB}(M/k)
\]
is an equivalence and that the circle action on $\HH^?_{\DB}(M/k)$
is trivial. Then the triangle
\[
\DHC^-(M/k)\longrightarrow\Sigma \HC^?(M/k)
\longrightarrow\Sigma \HC^?_{\DB}(M/k)
\]
is exact and splits in $\Mod_k$. Equivalently, its connecting map
\[
\HC^?_{\DB}(M/k)\longrightarrow \DHC^-(M/k)
\]
vanishes.
\end{prop}

\begin{proof}
Let
$\delta:\HH^?_{\DB}(M/k)\to\Sigma \DHH(M/k)$ be the Hochschild counterpart of the last map. The assumption on $\HP^?_\dagger$ implies $\DHH(M/k)^{tS^1}=0$, so the norm map
$\Sigma \DHH(M/k)_{hS^1}\to \DHH(M/k)^{hS^1}$ is an equivalence. This verifies the exactness of the triangle.

By functoriality, the square
\[
\begin{tikzcd}[column sep=large]
\Sigma \HH^?_{\DB}(M/k)_{hS^1}
 \arrow[r,"\Sigma\delta_{hS^1}"]\arrow[d,']
 &\Sigma^2\DHH(M/k)_{hS^1}\arrow[d,"\sim"]\\
\HH^?_{\DB}(M/k)^{hS^1}\arrow[r,"\delta^{hS^1}"']
 &(\Sigma \DHH(M/k))^{hS^1}
\end{tikzcd}
\] is commutative.

Triviality of the circle action on $\HH^?_{\DB}$ implies that the left vertical map is zero. Hence so is $\delta_{hS^1}$.
\end{proof}

\ \ 

We finish this section with a Du Bois version of A.~Blanc's Deligne cohomology. Recall that Blanc defined a
2-periodic localizing invariant
\[ \K^{\mathrm{top}}: \Catperf_\bbC \ra \Mod_{\KU}\] together with a morphism (the \emph{Chern character}) $\K^{\mathrm{top}}\ra \HP^{\cl}(-/\bbC)$ (\cite[Definition 4.13 and section 4.4]{Bla16}). Blanc also defined Deligne cohomology of smooth proper $\bbC$-categories (\cite[Definition 4.33]{Bla16}) as follows:
\[
H_{\mathcal D}(T):=
\K^{\mathrm{top}}(T)\mathop{\times}_{\HP^{\cl}(T/\bbC)}
\HC^{-,\cl}(T/\bbC).
\]

This definition, in fact, makes sense for an arbitrary $T\in \Catperf_\bbC$ and can be thought of as a derived Hodge extension of Deligne cohomology of smooth proper $\bbC$-schemes (see subsection 5.2 for further discussion). We suggest its Du Bois counterpart.

Consider the composite map
\[
\K^{\mathrm{top}}(T)\xrightarrow{\ch}\HP^{\cl}(T/\bbC)
\longrightarrow \HP^{\cl}_{\DB}(T/\bbC).
\]

\begin{defi}\label{DBDelignedefinition}
For $T\in\Catperf_\bbC$, define \emph{Du Bois-Deligne cohomology}
\[
H_{\mathcal D,\DB}(T):=
\K^{\mathrm{top}}(T)\mathop{\times}_{\HP^{\cl}_{\DB}(T/\bbC)}
\HC^{-,\cl}_{\DB}(T/\bbC).
\]
\end{defi}

The following proposition collects its main properties.

\begin{prop}\label{Deligneproperties}
Du Bois-Deligne cohomology
\[
H_{\mathcal D,\DB}: \Catperf_\bbC \lra \Sp
\]
is a truncating localizing invariant. It receives maps from homotopy K-theory and Blanc's Deligne theory completing the commutative square
\[
\begin{tikzcd}
\K\arrow[r]\arrow[d]&H_\mathcal{D}\arrow[d]\\
\KH\arrow[r]&H_{\mathcal{D},\DB}.
\end{tikzcd}
\]
The right vertical map is an equivalence whenever $\DK(T)=0$, in particular on the thick subcategory generated by motives of smooth proper $\bbC$-schemes.

\end{prop}

\begin{proof}
Both $\HC^{-,\cl}_{\DB}$ and $\HP^{\cl}_{\DB}$ are truncating localizing invariants by Proposition~\ref{HHDBfinitarytruncating}.
The invariant $\K^{\mathrm{top}}$ is localizing by
\cite[Proposition 4.15]{Bla16} and truncating by
\cite[Theorem 5.5]{Kon21}. Their pullback is therefore a truncating localizing invariant.

Compatibility of Blanc's topological Chern character with the negative cyclic Chern character gives the natural map
$\K\to H_{\mathcal D}$. The comparisons of cyclic theories give $H_{\mathcal D}\to H_{\mathcal D,\DB}$.

Since $\K^{\mathrm{top}}$ is $\bbA^1$-invariant, the natural map
$\K\to \K^{\mathrm{top}}$ factors canonically through $\KH$.
On the other hand, the bottom row of
\eqref{negativecycliccorrectiondiagram} gives
\[
\KH\longrightarrow \HC^{-,\cl}_{\DB}(-/\bbC).
\]
The resulting two maps from $\KH$ to $\HP^{\cl}_{\DB}(-/\bbC)$
agree: their restrictions along $\K\to \KH$ agree by compatibility of the Chern characters, and
$\HP^{\cl}_{\DB}$ is $\bbA^1$-invariant by Theorem~\ref{char0GW}.
This provides the commutative square.

Finally, if $\DK(T)=0$, then $\DHH(T/\bbC)=0$, so the comparisons on negative and periodic cyclic homology are equivalences. The same is consequently true of the Deligne comparison, this finishes the proof.
\end{proof}

The geometric interpretation of Du Bois-Deligne cohomology is given in Section 5.2.

\section{Du Bois localizing invariants in the commutative setting}\label{seccommutativeDB}

We now evaluate the invariants constructed above on schemes and stacks. The main purpose of this section is to identify them with familiar geometric theories, justifying the name ``Du Bois''. These comparisons also connect the categorical constructions of Sections~3--4 with the Hodge theory recalled in Subsection~2.3.

In the first subsection, we prove cdh descent for Du Bois invariants on ANS stacks and identify them with cdh-sheafifications on tame ANS stacks of finite type. We then give an independent argument for schemes.
In particular, Du Bois Hochschild homology recovers cdh-sheafified Hochschild homology.
The HKR comparison expresses it in terms of the Du Bois differential forms. For proper varieties, this identifies Du Bois Hochschild homology dimensions with column sums of Du Bois numbers, and we also obtain degeneration of the associated Hodge--de Rham spectral sequence.

In the second subsection, we apply this degeneration to show that, for proper varieties over $\bbQ$, $\DK$ (non-canonically) splits off of $\Sigma\HC$ as a direct summand. We then describe the rational weight components of the two Deligne theories: Blanc's theory uses the derived Hodge filtration, while its Du Bois counterpart recovers periodized classical Deligne--Beilinson cohomology.
We also discuss the comparison with Levine's modified Deligne--Beilinson theory.

\medskip

From now on, $k$ is a field of characteristic zero. 

\subsection{Du Bois invariants and cdh descent}

For a (derived) scheme or a derived stack $X$ and a localizing invariant $F$, we write
\[F(X):=F(\Perf(X)).\] We note that in the geometric situation the zoo of Hochschild and cyclic theories is considerably smaller.

\begin{obs}
Let $X$ be a (derived) quasi-geometric $k$-stack satisfying connective perfect generation and $F^?$ be one of the pairs of invariants considered in the previous section. Then the comparison map $F^{\mathrm{E}}(X) \ra F^{\cl}(X)$ is an equivalence.
\end{obs}

\begin{proof}

By Example~\ref{Nisnevich}, for any qcqs derived scheme $X$, $\Uloc(X)\in \Motlocwg$. In fact, the same holds for quasi-geometric $k$-stacks satisfying connective perfect generation in the sense of
\cite[Definition 6.21]{LS25}. Indeed,
\cite[Theorem 6.28]{LS25} makes $\Perf(X)$ a $c$-category, so Remark~\ref{LevySosnilo} places its motive in $\Motlocwg_k$.
Observation~\ref{obscomparisonclassical} and Proposition~\ref{HHDBfinitarytruncating} then identify the versions of $\TC$, $\TC_{\DB}$, and the associated Hochschild and cyclic theories: for each one of the listed theories the comparison map $F^{\mathrm{E}} \ra F^{\cl}$ is an equivalence.    
\end{proof}
From now on, we will stay within assumptions above, therefore we omit these superscripts from notation.

\medskip

One of the properties that introducing Du Bois invariants aimed to achieve is cdh-descent. We start with the following statement.

\begin{theorem}\label{DBstacksdescent}
Let $F\ra F_{\DB}$ be the canonical comparison for any of the Du Bois localizing invariants constructed in the paper, over a field $k$ of characteristic zero.
\begin{enumerate}
\item
For a derived ANS $k$-stack $X$, restriction to the classical locus induces an equivalence
\[
F_{\DB}(X)\xrightarrow{\sim}F_{\DB}(X_{\cl}).
\]

\item
On classical Noetherian ANS $k$-stacks, $F_{\DB}$
satisfies cdh descent. In particular, the comparison admits a canonical factorization
\[
F\longrightarrow F_{\cdh}\longrightarrow F_{\DB}.
\]

\item
On tame classical ANS stacks of finite type over $k$, the map
\[
F_{\cdh}\longrightarrow F_{\DB}
\]
is an equivalence, exhibiting the Du Bois invariants as cdh-sheafifications of the corresponding localizing invariants.
\end{enumerate}
\end{theorem}

Here ANS stands for ``affine diagonal, nice stabilizers''
and is understood in the sense of
\cite[Section~5.2]{ES21}.
The cdh topology on classical Noetherian ANS stacks is generated by Nisnevich covers and abstract blowup squares (cf. Subsection~2.3), with the proper morphism $\widetilde X\to X$ required to be representable; see \cite[Section~5.2]{ES21}. The last statement in the theorem is inspired by \cite{Hoy20}; one can find more on the cdh-topology on tame stacks there.

\begin{remark}
One can also include $\K$ and $\K^{\inf}$ in this theorem by setting $\K_{\DB}:=\KH$ and $\K^{\inf}_{\DB}:=\K^{\inf}$. The invariants $\KH$ and $\K^{\inf}$ are truncating; see, for example, \cite[Example~5.2.7]{ES21}.
\end{remark}

\begin{proof}
The Du Bois invariants are truncating, by Propositions~\ref{TCDBtruncating},
\ref{HHDBfinitarytruncating}, and
\ref{Deligneproperties}.
Thus (1) and the descent statement in (2) follow from
\cite[Corollaries~5.2.3 and~5.2.6]{ES21}, respectively. The results are stated for invariants defined on $\Catperf_\bbS$, but, just as \cite[Remark 5.1.4]{ES21} suggests, the proofs of \cite{ES21}
extend to the $k$-linear setting.
The factorization in (2) follows from the universal property of cdh-sheafification.

For (3), Nisnevich
descent reduces to tame quotient stacks; see \cite[Section~5.2, especially (5.2.1)]{ES21}.
In characteristic zero, stacks of finite type in this class admit cdh covers by regular stacks, as in
\cite[Corollary~1.4]{Hoy20}.
On a regular tame quotient stack $U$, the comparison
\[
\K(U)\xrightarrow{\sim}\KH(U)
\]
is an equivalence by \cite[Theorem~1.3]{Hoy20}.
Hence $\DK(U)=0$, and the constructions of the Du Bois invariants give
\[
F(U)\xrightarrow{\sim}F_{\DB}(U).
\]
Thus $F\to F_{\DB}$ is a cdh-local equivalence. Since $F_{\DB}$ satisfies cdh descent by (2), it identifies $F_{\DB}$ with $F_{\cdh}$.
\end{proof}

\medskip

We now specialize to schemes. In view of Theorem~\ref{DBstacksdescent}(1), from now on we assume that all schemes are classical.

\begin{theorem}\label{thmcdhTCDB}
For a separated finite-type $k$-scheme $X$ there is a canonical equivalence
\[
\TC_{\DB}(X)\simeq (\HC^-(-/\bbQ))_{\cdh}(X)=:\HC^-_{\cdh}(X/\bbQ).
\]
The base change using \eqref{TCEtoHH} gives an equivalence of mixed complexes
\begin{equation}\label{HHDBcdhHH}
\HH_{\DB}(X/k)\simeq \HH_{\cdh}(X/k)
\end{equation}
and equivalences of $\HC^-(k/k)$-modules
\begin{equation}\label{eqTCDBcdhHCminus}
\begin{split}
\HC^-_{\DB}(X/k)&\simeq \HC^-_{\cdh}(X/k),\\
\DHC^-(X/k)&\simeq\fib\bigl(\HC^-(X/k)\to \HC^-_{\cdh}(X/k)\bigr).
\end{split}
\end{equation}
\end{theorem}

The idea behind the first statement of the theorem is, of course, well-known to experts, see, for example, \cite[Section 8.1]{EKS25} and \cite[Appendix A]{Loe26}.

We give an argument independent of Theorem~\ref{DBstacksdescent}.

\begin{proof}
By \cite[Theorem 4.6]{CHSW08}, $\K^{\mathrm{inf}}$ satisfies cdh descent. Moreover
$\K_{\cdh}\simeq \KH$, compatibly with $\K\to \KH$; see
\cite[Theorem 6.4]{Hae04} and \cite[Example 3.14]{CHSW08}. Applying the exact cdh localization to
$\K^{\mathrm{inf}}\to \K\to \HC^-(-/\bbQ)$ therefore gives
\[
(\HC^-(-/\bbQ))_{\cdh}
 \simeq\cof(\K^{\mathrm{inf}}\to \KH)\simeq \TC_{\DB}.
\]
The last identification uses the lower row of the motivic square of Corollary~\ref{multiplicativeinvariants}. Equivalently the right vertical map in that square gives an exact triangle in $\Mod_{\HC^-(k/\bbQ)}$
\[
\DK(X)\longrightarrow \HC^-(X/\bbQ)
 \longrightarrow(\HC^-(-/\bbQ))_{\cdh}(X)
\] exhibiting the equivalence above as \emph{an equivalence of $\HC^-(k/\bbQ)$-modules}.

Recall (see Section~\ref{secDBrecollections}) that cdh-sheafifications of localizing invariants do not change the values of the invariants on smooth schemes. $X$ admits a bounded cubical hyperresolution by smooth
$k$-schemes, so the cdh-sheafifications on $X$ 
are computed by finite homotopy limits. These commute with killing $u$ and with the base change in \eqref{TCEtoHH}. By Lemma~\ref{equivariantscalarreduction}, in $\Mod_k^{BS^1}$, one has
\[
\HH_{\DB}(X/k)
\simeq \TC_{\DB}(X)\otimes_{\HC^-(k/\bbQ)}k\\
\simeq \HC^-_{\cdh}(X/\bbQ)\otimes_{\HC^-(k/\bbQ)}k
 \simeq \HH_{\cdh}(X/k).
\]
and
\[
\HC^-_{\DB}(X/k)
\simeq \HH_{\DB}(X/k)^{hS^1}
\simeq \HH_{\cdh}(X/k)^{hS^1}
 \simeq \HC^-_{\cdh}(X/k).
\]
The last assertion follows.
\end{proof}

The smooth HKR decomposition gives an identification of underlying $k$-modules
\begin{equation}\label{eqcdhHKR}
\HH_{\DB}(X/k)\simeq
 \bigoplus_{p\geq0}\Gamma(X,\underline\Omega^p_{X/k})[p];
\end{equation}
the circle action corresponds to the de Rham differential on the right hand term. The latter has only finitely many nonzero summands; the contribution of the $p$-th summand to $\HH_{\DB, i}$ is
$\mathbb H^{p-i}(X,\underline\Omega^p)$.

Passing to dimensions, one sees what information about Hodge numbers is stored in $\HH_{\DB}$. We refer to Section~\ref{secDBrecollections} for the definitions of the Du Bois numbers and mixed Hodge numbers for a proper variety over an arbitrary characteristic zero field $k$. 

\begin{corollary}\label{DBHodgeformula}
Let $X/k$ be a proper variety. Then, for the Du Bois-Hochschild numbers of $\Perf(X)$, the following holds:
\begin{equation}\label{eqcolumnHodgesum}
hh_i(X):=\dim_k \HH_{\DB,i}(X/k)
 =\sum_{p-q=i}\underline h^{p,q}(X)
 =\sum_{\substack{n,a,b\\2a-n=i}}h_n^{a,b}(X).
\end{equation}
The Hodge--de Rham spectral sequence from $\HH_{\DB}(X/k)$
to $\HP_{\DB}(X/k)$ degenerates at $E_2$.
\end{corollary}

When $X$ is smooth, this recovers the usual noncommutative Hodge--de Rham spectral sequence for $\Perf(X)$, considered, for example, in \cite{Kal17}.

\begin{proof}
The dimension formula follows from \eqref{eqcdhHKR} and Proposition~\ref{numbersoverk}, with $a=p$ and $n=p+q$, so that $p-q=2a-n$.

Under the HKR comparison, periodizing the bounded Du Bois filtered complex identifies the $u$-adic spectral sequence, weight by weight and up to reindexing, with the Du Bois Hodge-to-de Rham spectral sequence.
The latter degenerates at $E_1$ by Proposition~\ref{DBHodgecharzero}; boundedness supplies convergence. Indexing the cyclic filtration by the cohomological degree of $u$, rather than its exponent, doubles the filtration indices. Thus the differential $d_r$ corresponds to
$d_{2r}$ in the usual indexing of the $S^1$-Tate spectral sequence, while the odd differentials vanish. The resulting spectral sequence therefore degenerates at $E_2$.

\end{proof}

\medskip

\subsection{DK-splitting and Deligne cohomology of proper singular varieties}

The same argument as above also identifies the canonical map
\[
\HP(X/k)\longrightarrow \HP_{\DB}(X/k)
\]
as an equivalence: ordinary periodic cyclic homology satisfies cdh descent in characteristic zero \cite[Corollary 3.13]{CHSW08}, and existence of bounded smooth hyperresolutions implies that cdh-sheafification commutes with the Tate construction. Theorem~\ref{char0GW} consequently identifies
$\HC^-_{\DB}(X/k)\to \HP(X/k)$ with its $\bbA^1$-homotopification map.

For proper $X/k$, Corollary~\ref{DBHodgeformula} also implies that the circle action on $\HH_{\DB}(X/k)$ is (noncanonically)
trivial.
Therefore the assumptions of Proposition~\ref{DHCsplitting} are satisfied by $T=\Perf(X)$. Moreover, for $k=\bbQ$, Theorem~\ref{thmcdhTCDB} identifies $\DHC^-(X/\bbQ)$ with $\DK(X)$.

We therefore have proved the following statement.

\begin{corollary}\label{DKsplittingvarieties}
Let $X/k$ be a proper variety. Then the triangle
\[
\DHC^-(X/k)\longrightarrow\Sigma \HC(X/k)
\longrightarrow\Sigma \HC_{\DB}(X/k)
\]
splits as a triangle of $k$-modules. Equivalently, its connecting map
\[
\HC_{\DB}(X/k)\longrightarrow \DHC^-(X/k)
\]
vanishes. In particular, for a proper variety $X$ over $\bbQ$, there is a split triangle of $\bbQ$-modules
\[
\DK(X)\longrightarrow\Sigma \HC(X/\bbQ)
\longrightarrow\Sigma \HC_{\DB}(X/\bbQ).
\]
\end{corollary}

This proves beneficial for computations of $\DK$ and, ultimately, algebraic K-theory of singular varieties; we will return to this elsewhere.

\medskip

We now justify our definition of Du Bois-Deligne cohomology (Definition~\ref{DBDelignedefinition}).

Let $X$ be a proper complex scheme. Using the notation of Section~\ref{secDBrecollections}, define
\[
\Gamma_{\mathcal D}^{\der}(X,\bbQ(p))
:=
\fib(
\Gamma(X^{\an},\bbQ(p))
\longrightarrow
\Gamma(X,\widehat{\dR}_{X/\bbC})/
F_{\der}^p
),
\]
and
\[
\Gamma_{\mathcal D, \DB}(X,\bbQ(p))
:=
\fib(
\Gamma(X^{\an},\bbQ(p))
\longrightarrow
\Gamma(X,\underline\Omega^\bullet_{X/\bbC})/F^p
).
\]
Here $\bbQ(p)=(2\pi i)^p\bbQ$, the quotients are cofibers, and the maps are supplied by the de Rham comparisons
\[
\Gamma(X,\widehat{\dR}_{X/\bbC})
\xrightarrow{\sim}
\Gamma(X,\underline\Omega^\bullet_{X/\bbC})
\simeq \Gamma(X^{\an},\bbC).
\]
For $p\leq0$, $F^p_{\der}$ and $F^p$ recover the whole complex.

The second fiber computes the classical Deligne--Beilinson cohomology of a proper singular variety, defined using smooth proper hypercovers
\cite[\S5, especially Definition 5.3]{EV88}.
The first is its derived Hodge counterpart. The filtered comparison induces a natural map
\[
\Gamma_{\mathcal D}^{\der}(X,\bbQ(p))
\longrightarrow
\Gamma_{\mathcal D, \DB}(X,\bbQ(p)),
\]
which is an equivalence when $X$ is smooth.

\begin{corollary}\label{DeligneDBrational}
For a proper complex scheme $X$, there are natural equivalences
\[
\begin{aligned}
H_{\mathcal D}(\Perf(X))_{\bbQ}
&\simeq
\prod_{p\in\bbZ}
\Gamma_{\mathcal D}^{\der}(X,\bbQ(p))[2p],\\
H_{\mathcal D,\DB}(\Perf(X))_{\bbQ}
&\simeq
\prod_{p\in\bbZ}
\Gamma_{\mathcal D, \DB}(X,\bbQ(p))[2p].
\end{aligned}
\]
Under these equivalences, the categorical comparison is the product of the maps induced by the filtered comparison from derived de Rham to Du Bois cohomology.
\end{corollary}

See Section 4.2 for the definitions of the left hand side invariants.

\begin{proof}
The completed characteristic-zero HKR comparison (see, for example,
\cite[Theorem 1.1]{Ant19}) identifies the negative-to-periodic cyclic homology map with
\[
(
\HC^-(X/\bbC)\longrightarrow \HP(X/\bbC)
)
\simeq
(\prod_{p\in\bbZ}
F_{\der}^p
\Gamma(X,\widehat{\dR}_{X/\bbC})[2p]
\longrightarrow
\prod_{p\in\bbZ}
\Gamma(X,\widehat{\dR}_{X/\bbC})[2p]).
\]
Theorem~\ref{thmcdhTCDB} similarly identifies the Du Bois cyclic homology map with
\[
(\HC^-_{\DB}(X/\bbC)\longrightarrow \HP_{\DB}(X/\bbC))
\simeq
(\prod_{p\in\bbZ}
F^p\Gamma(X,\underline\Omega^\bullet_{X/\bbC})[2p]
\longrightarrow
\prod_{p\in\bbZ}
\Gamma(X,\underline\Omega^\bullet_{X/\bbC})[2p]).
\]
These identifications are compatible with the natural comparison between the two theories.

Blanc's comparison with topological $\K$-theory
\cite[Proposition 4.32]{Bla16} and the rational topological Chern character give
\[
\K^{\mathrm{top}}(\Perf(X))_{\bbQ}
\simeq \KU(X^{\an})_{\bbQ}
\simeq
\prod_{p\in\bbZ}\Gamma(X^{\an},\bbQ(p))[2p].
\]
Boundedness of Betti cohomology implies that each cohomological degree gets contributions only from a finite number of factors of the right term product.
Compatibility with the maps to periodic cyclic homology can be checked on smooth proper schemes: the periodic comparison is an equivalence, and the existence of bounded smooth proper hyperresolutions allows one to use finite descent and proper GAGA.

Rationalization preserves the pullbacks defining
$H_{\mathcal D}$ and $H_{\mathcal D,\DB}$, whose cyclic homology terms are already rational. Products commute with pullbacks, so the displayed comparisons identify these pullbacks with the products of the defining weightwise fibers.
\end{proof}

Thus Blanc's version of Deligne cohomology computed on proper singular schemes is related to the derived Hodge filtration, whereas its Du Bois counterpart recovers periodized classical Deligne--Beilinson cohomology after rationalization.

\begin{remark}\label{LevineDeligne}
An intermediate theory is M.~Levine's modified Deligne--Beilinson cohomology, which uses ordinary K\"ahler differentials \cite{Lev92}.
For a projective complex scheme $X$, its rational defining complex is
\[
\Gamma_{\mathcal D^*}(X,\bbQ(p))
\simeq
\fib(\Gamma(X^{\an},\bbQ(p))
\longrightarrow
\Gamma(X^{\an},\Omega_{X^{\an}}^{<p})),
\]
where $\Omega_{X^{\an}}^{<p}$ is the quotient of the ordinary de Rham complex by its stupid filtration $F^p$.
The natural comparisons of differential-form complexes give
\begin{equation}\label{Deligneformscompare}
\Gamma_{\mathcal D}^{\der}(X,\bbQ(p))
\xrightarrow{\alpha^{\der}_{\mathrm{naive}}}
\Gamma_{\mathcal D^*}(X,\bbQ(p))
\xrightarrow{\alpha^{\mathrm{naive}}_{\DB}}
\Gamma_{\mathcal D,\DB}(X,\bbQ(p)).
\end{equation}
These theories use derived differential forms, ordinary K\"ahler forms, and the Du Bois complex, respectively;
the comparisons are equivalences for smooth proper $X$.
\end{remark}

\begin{remark}
In \cite[\S2, p.~125]{Lev92}, M.~Levine constructs a Chern character whose weight components, for a projective complex scheme $X$, have the form
\[
\ch^{\mathrm{Lev}}_{p,n}:
\K_n(X)_{\bbQ}
\longrightarrow
H_{\mathcal D^*}^{2p-n}(X,\bbQ(p)),
\qquad p,n\geq0.
\]
On the other hand, Proposition~\ref{Deligneproperties} combined with Corollary~\ref{DeligneDBrational} provides two Chern characters with weight components
\[
\begin{aligned}
\ch^{\der}_{p,n}:&
\K_n(X)_{\bbQ}
\longrightarrow
H_{\mathcal D,\der}^{2p-n}(X,\bbQ(p)),\\
\ch_{\DB,p,n}:&
\KH_n(X)_{\bbQ}
\longrightarrow
H_{\mathcal D,\DB}^{2p-n}(X,\bbQ(p)).
\end{aligned}
\]
Using the comparison maps \eqref{Deligneformscompare}
from Remark~\ref{LevineDeligne}, we expect that, with compatible Chern-character normalizations,
\[
H^{2p-n}(\alpha^{\der}_{\mathrm{naive}})
\circ\ch^{\der}_{p,n}
=
\ch^{\mathrm{Lev}}_{p,n}
\qquad\text{and}\qquad
H^{2p-n}(\alpha^{\mathrm{naive}}_{\DB})
\circ\ch^{\mathrm{Lev}}_{p,n}
=
\ch_{\DB,p,n}
\circ\bigl(\K_n(X)_{\bbQ}\longrightarrow \KH_n(X)_{\bbQ}\bigr).
\]
Thus Levine's Chern character should lift to the derived Hodge version of Deligne cohomology, while its image in Du Bois--Deligne cohomology should factor through homotopy $\K$-theory.
\end{remark}

\medskip

For convenience of future reference, we also note the following observation.

\begin{remark}
Theorem~\ref{DBstacksdescent} applies to Du Bois--Deligne cohomology. In particular, on separated finite-type $\bbC$-schemes, and more generally on tame classical ANS stacks of finite type over $\bbC$, the canonical comparison identifies
\[
(H_{\mathcal D})_{\cdh}\simeq H_{\mathcal D,\DB}.
\]
Thus Du Bois--Deligne cohomology is the cdh-sheafification of Blanc's Deligne theory.
For proper complex schemes, it can also be computed using smooth proper hyperresolutions.
\end{remark}

\medskip

\section{Fourier--Mukai invariance and mixed Hodge numbers}

The geometric comparison of the previous section produces numerical invariants of the perfect derived category of a proper variety. We begin to investigate how much information about mixed Hodge numbers these invariants retain.

In the first subsection, we introduce the Du Bois--Hochschild numbers, together with their support and center of mass. 
For a nonempty proper variety, we prove that this center of mass vanishes if and only if its cohomology in degree $n$ is pure of weight $n$ for every $n$.
Thus purity can be detected from the perfect derived category over any field of characteristic zero.
We also show that if a compact motive is the motive of a smooth proper category, then its ordinary Hochschild support, when nonempty, is symmetric about zero.

In the second subsection, we formulate the mixed Hodge numbers invariance question for Fourier--Mukai partners and record the equalities of column sums of Du Bois numbers provided by Du Bois Hochschild homology. In the third subsection, we apply these equalities to normal projective complex surfaces with quasi-rational singularities. We compute their mixed Hodge numbers in terms of the resolution and exceptional divisor, and obtain reconstruction results under additional assumptions. Using the categorical invariants introduced in Section~7, we will extend the reconstruction result in Section~8 to all connected normal projective complex surfaces.

\medskip

As before, $k$ is a field of characteristic zero.

\subsection{Du Bois-Hochschild numbers, purity, Hochschild support, and compact motives}

For a motive $M\in\Motloc_k$ such that $\HH^{?}_{\DB}(M/k)$ has finite-dimensional total homology, set $hh_i^{?}(M)=\dim_k \HH^{?}_{\DB,i}(M/k)$. For simplicity, we assume that the classical and Efimov versions produce the same numbers for $M$, so we can set $hh_i(M):=hh_i^?(M)$. As usual, given a distribution of numerical invariants over $\bbZ$, we can cook up a few new invariants reflecting the properties of this distribution and, ultimately, of the motive $M$. For example, one can define the \emph{Du Bois-Hochschild support interval} of a motive as  the smallest interval of integers containing
$\{i\in\bbZ:hh_i(M)\ne0\}$; its width, \emph{Du Bois amplitude} of $M$ provides an interesting dimension-type invariant. Similarly one can define the \emph{Hochschild support interval} and the \emph{amplitude} of a motive $M$ whenever it has finite-dimensional usual (say, classical) total Hochschild homology. We show that this is already an interesting invariant by proving the following statement inspired by work of Shklyarov \cite[Proposition 3.6 and Theorem 6.2]{Shk08}.

\begin{prop}
Suppose that $M$ is a compact motive: $M\in (\Motloc_k)^\omega$. Then $\HH^{cl}(M/k)$ is perfect, so the classical Hochschild support interval and the classical amplitude of $M$ are correctly defined. If, moreover, $M=\Uloc(T)$ for $T$ a smooth proper $k$-category, then the classical Hochschild support interval is empty or has zero as the center of symmetry, i.e. is of the form $[-d, d]$ for an integer $d\ge 0$. Finally, if $M=\Uloc(\Perf(X))$ for a smooth proper nonempty variety $X$, then $d\le \dim X$.
\end{prop}

\begin{proof}
By rigidity of $\Motloc_k$ \cite{Efi25b}, every compact motive is dualizable. The classical Hochschild realization
\[
\HH^{cl}(-/k):\Motloc_k\longrightarrow\operatorname{Mod}_k
\]
is symmetric monoidal, hence preserves dualizable objects.
Since the dualizable $k$-modules are precisely the perfect complexes, $\HH^{cl}(M/k)$ is perfect. 

If $T$ is smooth and proper, it is dualizable in
$\Catperf_k$, with dual $T^{op}$; see
\cite[Theorem 3.7]{BGT13} (their proof applies in the $k$-linear setting as well). Symmetric monoidal functors respect duals of dualizable objects, hence
\[
\HH^{cl}(T^{op}/k)
\simeq
\Hom_k(\HH^{cl}(T/k),k).
\]
This is the perfect Hochschild pairing studied by Shklyarov; see \cite[Theorem 6.2]{Shk08}.
Reversal of Hochschild chains gives an equivalence of \emph{underlying $k$-modules}
\[
\HH^{cl}(T^{op}/k)\simeq\HH^{cl}(T/k),
\]
hence
\[
\dim_k\HH^{cl}_i(T/k)
=
\dim_k\HH^{cl}_{-i}(T/k)
\]
for every $i$. Thus its support interval, when nonempty, is symmetric about zero.

Finally, for a smooth proper variety $X$, HKR gives
\[
\HH^{cl}_i(X/k)
\simeq
\bigoplus_{p-q=i}H^q(X,\Omega^p_{X/k}).
\]
This vanishes for $|i|>\dim X$ and $0\ne 1\in H^0(X,\mathcal{O}_X)$ unless $X$ is empty, contributing to $hh_0$.
\end{proof}

In particular, this gives a necessary condition for a compact motive to be a motive of a smooth proper category. This condition is obviously not vacuous: it fails for compact motives of the form $\Sigma^m\mathbbm 1$ for $m\ne 0$, which have the support interval $[m, m]$. We will further explore this story elsewhere.

\begin{remark}
To prove the first claim one could instead combine \cite[Proposition 1.9(i)]{Efi25a} with rigidity of $\Motloc_k$, like we did proving \eqref{basemotloc}. Similarly, combining it with \cite[Proposition 1.12(iv)]{Efi25a} gives the projection formula:
\[
M\otimes_k R^{\HH}_k(V)
\xrightarrow{\sim}
R^{\HH}_k\bigl(\HH^{cl}(M/k)\otimes_k V\bigr),
\qquad M\in\Motloc_k,\quad V\in\Mod_k.
\]
where $R^{\HH}_k:\Mod_k\lra \Motloc_k$ is the right adjoint to $\HH^\cl(-/k)$. In particular, the adjunction implies
\[
\Gamma(-\otimes_k R^{\HH}_k(k))
\simeq
\Gamma R^{\HH}_k(\HH^{cl}(-/k))
\simeq
\Hom_k(k,\HH^{cl}(-/k))
\simeq \HH^{cl}(-/k)
\]
providing an incarnation of Hochschild homology:
\[
\HH^{\cl}(-/k) \simeq \Gamma(-\otimes_kM_{\HH, k}) : \Motloc_k \lra \Mod_k
\]
for $M_{\HH, k}:= R^{\HH}_k(k)$ (of course, existence of such $M_{\HH, k}$, at least for $\Sp$-valued Hochschild homology, is a formal corollary of rigidity of $\Motloc_k$) and a canonical $k$-action on $M_{\HH, k}$. The same construction works for an arbitrary base $E_\infty$-ring spectrum.
\end{remark}

\medskip

If the total dimension of $\HH_{\DB}(M/k)$ is nonzero, one can also define its \emph{center of mass} (or "mean value") by the usual formula:
\[
\mu(M)=\frac{\sum_i ihh_i(M)}{\sum_i hh_i(M)}\in \bbQ.
\]
This assigns to a motive a rational number which is an interesting invariant on its own. In particular, it measures "non-purity of the motive" in the sense that will be made clear below. 

When considering purity, it is convenient to eliminate the denominator of $\mu$. Put
\[
\delta(M):=\sum_i ihh_i(M) \in \bbZ.
\]

We give the following definition in the spirit of Section~\ref{secDBrecollections}.

\begin{defi}
Let $X$ be a proper $k$-variety. Equip
$\mathbb H^n(X,\widehat{\dR}_{X/k})$
with the weight filtration transported from
$H^n_{\DB}(X/k)$ through the canonical de Rham comparison of Section~\ref{secDBrecollections}.
We say that $X$ \emph{has pure cohomology} if
\[
\Gr_m^W\mathbb H^n(X,\widehat{\dR}_{X/k})=0
\qquad\text{for every }n\text{ and every }m\ne n.
\]
Equivalently, for every $n$,
\[
W_{n-1}\mathbb H^n(X,\widehat{\dR}_{X/k})=0,
\qquad
W_n\mathbb H^n(X,\widehat{\dR}_{X/k})
=\mathbb H^n(X,\widehat{\dR}_{X/k}).
\]
\end{defi}

By Proposition~\ref{DBHodgecharzero}, this property is preserved and reflected by extensions of the ground field.
When $k=\bbC$, it is equivalent to the Betti cohomology
$H^n(X^{\an},\bbQ)$ being pure of weight $n$
for every $n$.

\begin{theorem}\label{thmpurity}
Let $X$ be a nonempty proper $k$-variety. Put $\delta(X) = \delta(\Uloc(\Perf(X)))$ and $\mu(X) = \mu(\Uloc(\Perf(X)))$. Then
\begin{equation}\label{eqweightdefect}
\delta(X)\le 0
\end{equation}
and the following conditions are equivalent:
\begin{enumerate}
\item $\mathbb{H}^n(X,\widehat{\dR}_{X/k})$ is pure of weight $n$ for every $n$;
\item $\delta(X)=0$;
\item $\mu(X) = 0$;
\item $hh_i(X)=hh_{-i}(X)$ for every $i$.
\end{enumerate}
Consequently the property of having pure cohomology of $X/k$ is Morita-invariant, i.e. it depends only on its perfect derived category $\Perf(X)$.
\end{theorem}

\begin{proof}

We substitute \eqref{eqcolumnHodgesum}:
\[
\begin{split}
\delta = \sum_i ihh_i =\sum_{n,a,b}(2a-n)h_n^{a,b}.\end{split}
\]
Write $(a,b)$ for an unordered pair. Recall that Hodge symmetry, as stated in Proposition~\ref{numbersoverk}, provides $h^{a,b}_n=h^{b,a}_n=:h^{(a,b)}_n$. When $a\ne b$, the couple ($h^{a,b}_n,h^{b,a}_n$) contributes $h^{(a,b)}_n$ to $hh_{2a-n}$ and to $hh_{2b-n}$ and makes the total contribution of
\[(2a-n+2b-n)h^{(a,b)}_n \]
to $\delta$. The contribution of $h^{a,a}_n$ to $\delta$ is $(2a-n)h^{(a,a)}_n$ via the summand $hh_{2a-n}$. In both cases, the contribution to $\delta$ is nonpositive, since $a+b\le n$ (Proposition~\ref{numbersoverk}); it is zero iff $a+b=n$ or $h^{(a,b)}_n=0$. Thus $\delta\leq0$, with equality precisely when every nonzero $h_n^{a,b}$ satisfies $a+b=n$. This proves equivalence between (1) and (2) and  equivalence between (1) and (4). The equivalence between (2) and (3) is obvious.
\end{proof}

\begin{remark}
The Hodge diamond of a smooth proper complex variety has symmetries coming from complex conjugation and Poincar\'e duality. For proper singular varieties, non-purity forces the corresponding symmetries of the Du Bois numbers to fail, even though Hodge symmetry still holds on each pure weight quotient.
Above we have seen that Du Bois--Hochschild numbers detect the failure of the "horizontal" symmetry of Du Bois numbers.
In the section 8, we will make use of topological K-theory to measure the failure of Poincare duality.
\end{remark}

\subsection{Fourier-Mukai partners and the mixed Hodge numbers invariance question}

The rest of the paper is built around the following question.

\begin{ques}\label{Hodgerecoveryquestion}
Let $X/k$ be a connected proper variety. To what extent can one recover the mixed Hodge numbers of $X$ (in the sense of Definition~\ref{numbersdefi}) from its perfect derived category $\Perf(X)$?
\end{ques}

Following tradition, we will also consider the following formulation.

\begin{ques}\label{FMquestion}
Let $X$ and $Y$ be connected proper Fourier-Mukai partners over $k$. Is it true that $h^{a,b}_n(X)=h^{a,b}_n(Y)$ for every $(a,b,n)\in \bbN^3$?
\end{ques}

Here Fourier-Mukai partners are understood in the sense of Definition~\ref{defiFM}; see also Proposition~\ref{propderivedequivalence} for equivalent definitions.

\medskip

Corollary~\ref{DBHodgeformula} makes a significant step in the direction of Question~\ref{Hodgerecoveryquestion}. Here is its Fourier--Mukai formulation.

\begin{theorem}\label{FMcolumns}
Let $X$ and $Y$ be proper Fourier-Mukai partners over $k$. Then for every $i\in\bbZ$
\[\sum_{\substack{n,a,b\\2a-n=i}}h_n^{a,b}(X)
 =\sum_{p-q=i}\underline h^{p,q}(X)
 =\sum_{p-q=i}\underline h^{p,q}(Y)
 =\sum_{\substack{n,a,b\\2a-n=i}}h_n^{a,b}(Y) .\]
\end{theorem}

\begin{proof}
Follows from Corollary~\ref{DBHodgeformula}.
\end{proof}

\begin{remark}
Note that the theorem above does not require the full strength of the Du Bois localizing invariants theory. It suffices to show that there exists a $\Mod_{k^{hS^1}}$-valued invariant $\HC^{-,\mathrm{lazy}}_{\DB}(-/k)$ defined on ttg motives (see \ref{defittg}) which has properties of a localizing truncating invariant on this subcategory of $\Motloc_k$ and naturally agrees with $\HC^-(-/k)$ on smooth proper $k$-schemes as a $\Mod_{k^{hS^1}}$-valued functor.
By the argument of \cite[Appendix A]{LT19}, motives of arbitrary proper $k$-schemes are ttg, hence one would have a well-defined $k[[u]]$-module $\HC^{-,\mathrm{lazy}}_{\DB}(X/k)$ for such schemes depending only on $\Perf(X)$ and isomorphic to $\HC^{-}_{\cdh}(X/k)$ as a $k[[u]]$-module. By defining $\HH_{\DB}^{\mathrm{lazy}}(-/k):=\HC^{-,\mathrm{lazy}}_{\DB}(-/k)\otimes_{k^{hS^1}}k$, one then recovers cdh-sheafified Hochschild homology, arriving at Corollary~\ref{DBHodgeformula} and Theorem~\ref{FMcolumns}.    
\end{remark}

\medskip

\subsection{Du Bois-Hochschild numbers of surfaces}

To simplify the exposition, from now on, we work with projective varieties over complex numbers. All the cohomology groups and mixed Hodge structures are, by default, with $\bbQ$-coefficients.

Below we show that Theorem~\ref{FMcolumns} already allows to address Question~\ref{FMquestion} for a large class of complex projective surfaces. In the next section, we consider other invariants that together with Corollary~\ref{DBHodgeformula}  allow us to attack the Question~\ref{FMquestion} in higher dimensions.

\medskip

Let $X$ be a connected normal projective complex surface, and let
$\pi:\widetilde X\to X$ be a projective resolution which is an isomorphism over the smooth locus. Suppose that the reduced exceptional divisor
\[
E=\coprod_{s\in\Sing(X)}E_s
\]
has simple normal crossings. Write
\begin{equation}\label{curveHodgenumbers}
L=\sum_s b_1(\Gamma_{E_s})=h^{0,0}_1(E),\qquad
G=\sum_{E_j\subset E}g(E_j)=h^{1,0}_1(E),
\end{equation}
where the second sum runs over the irreducible components of $E$.
Here $\Gamma_{E_s}$ is the dual graph (as in e.g.  \cite{Wei01}): its vertices correspond to the irreducible components of $E_s$, and its edges to their intersection points. Thus $L$ counts the independent cycles in these graphs, while $G$ is the sum of the genera of the exceptional components. For smooth $X$, both sums are zero.

We use Arapura--Hiatt's notion of quasi-rational singularities.
By \cite[Corollary 5.4]{AH26}, a normal projective surface $X$
has quasi-rational singularities if and only if every exceptional component is contracted by the Albanese map of $\widetilde X$.
We take this characterization as our definition here; its equivalence with the intrinsic definition
\cite[Definition 5.1]{AH26} shows independence of the resolution.

\begin{lemma}\label{qratvanishing}
Let $X$ be a connected normal projective surface with quasi-rational singularities and resolution as above. Then
\begin{equation}\label{eqsurfacerestriction}
H^1(\widetilde X,\bbQ)\longrightarrow H^1(E,\bbQ)
\quad\text{is zero}.
\end{equation}
\end{lemma}

\begin{proof}
The Albanese map induces an isomorphism
\[
H^1(\operatorname{Alb}(\widetilde X),\bbQ)
\xrightarrow{\sim}H^1(\widetilde X,\bbQ).
\]
Since it contracts every $E_j$, the restriction
$H^1(\widetilde X,\bbQ)\to H^1(E_j,\bbQ)$ is zero.

The normalization $\coprod_j E_j\to E$ gives an exact sequence of mixed Hodge structures
\[
0\longrightarrow
\bigoplus_s H^1(\Gamma_{E_s},\bbQ)
\longrightarrow H^1(E,\bbQ)
\longrightarrow\bigoplus_j H^1(E_j,\bbQ)
\longrightarrow0.
\]
The left term is pure of weight zero and the right term is pure of weight one. Hence restriction to the components identifies
\[
\Gr_1^W H^1(E,\bbQ)
\simeq\bigoplus_j H^1(E_j,\bbQ).
\]
By the vanishing above, the image of
$H^1(\widetilde X,\bbQ)\to H^1(E,\bbQ)$ lies in
$W_0H^1(E,\bbQ)$, hence it vanishes, due to strictness of the weight filtration, since $H^1(\widetilde X,\bbQ)$ is pure of weight one.
\end{proof}

\begin{prop}\label{propsurfaceformulas}
Let $X$ be a connected normal projective surface with quasi-rational singularities and resolution as above. Put
\[
q=h^{1,0}(\widetilde X),\qquad
p_g=h^{2,0}(\widetilde X).
\]
As in the proof of Theorem~\ref{thmpurity}, write $(a,b)$
for an unordered pair and
$h_n^{(a,b)}=h_n^{a,b}=h_n^{b,a}$.
The mixed Hodge numbers of $X$ are
\[
\begin{array}{c|c|c}
n&(a,b)&h_n^{(a,b)}(X)\\\hline
0&(0,0)&1\\
1&(0,1)&q\\
2&(0,0)&L\\
2&(0,1)&G\\
2&(0,2)&p_g\\
2&(1,1)&h_2^{1,1}(X)\\
3&(1,2)&q\\
4&(2,2)&1
\end{array}
\]
In particular, $q$
and $p_g$ do not depend on the choice of resolution.

All other mixed Hodge numbers are zero. 

For $X$ we have
\begin{equation}\label{eqsurfacediagonals}
\begin{aligned}
hh_2&=p_g,\\
hh_1=hh_{-1}&=2q,\\
hh_0&=2+h_2^{1,1}(X)+G,\\
hh_{-2}&=p_g+G+L,
\end{aligned}
\end{equation}
and $hh_i=0$ for $|i|>2$.

Equivalently, the mixed Hodge numbers are expressed in terms of the $hh_i$ and $G$ by
\begin{equation}\label{eqsurfacemixedhodgenumers}
\begin{aligned}
h_0^{0,0}(X)=h_4^{2,2}(X)&=1,\\
h_1^{(0,1)}(X)=h_3^{(1,2)}(X)&=\frac{hh_1}{2},\\
h_2^{(0,2)}(X)&=hh_2,\\
h_2^{(0,1)}(X)&=G,\\
h_2^{0,0}(X)&=hh_{-2}-hh_2-G=L,\\
h_2^{1,1}(X)&=hh_0-2-G.
\end{aligned}
\end{equation}
\end{prop}

\begin{proof}
The abstract blowup square
\[
\begin{tikzcd}
E \arrow[r,hook] \arrow[d]
  & \widetilde X \arrow[d,"\pi"]\\
\Sing(X) \arrow[r,hook]
  & X
\end{tikzcd}
\]
gives an exact sequence of mixed Hodge structures

\begin{equation}\label{lesMHSsurface}
\begin{aligned}
0\longrightarrow H^1(X)
&\longrightarrow H^1(\widetilde X)
\xrightarrow{0} H^1(E)
\longrightarrow H^2(X)
\longrightarrow H^2(\widetilde X)\\
&\xrightarrow{\phi} H^2(E)
\xrightarrow{0} H^3(X)
\longrightarrow H^3(\widetilde X)
\longrightarrow0\longrightarrow H^4(X)\simeq H^4(\widetilde X)
\end{aligned}
\end{equation}

with the last term isomorphic to $\bbQ(-2)$.

Connectedness of exceptional fibers gives surjectivity of $H^0(\Sing(X))\ra H^0(E)$, hence the initial zero.

The intersection matrix of the exceptional curves is negative definite. Restricting their divisor classes to the components of $E$ therefore shows that
\[
H^2(\widetilde X)\longrightarrow
H^2(E)\simeq\bigoplus_j\bbQ(-1)
\]
is surjective, hence the next map is zero. The other zero in \eqref{lesMHSsurface} is provided by Lemma~\ref{qratvanishing}.

It remains to understand the four-term exact sequence in the middle. 
$H^2(\widetilde X)$ is pure of weight two and $H^2(E)$
has only type $(1,1)$, hence $\ker(\phi)$ has $(2,0)$ and $(0,2)$ parts, both of dimension $p_g$.

Hodge symmetry and Poincar\'e duality on
$\widetilde X$ give
\[
h^{1,0}(\widetilde X)=h^{0,1}(\widetilde X)
=h^{2,1}(\widetilde X)=h^{1,2}(\widetilde X)=q.
\]
Combining these observations with \eqref{curveHodgenumbers} proves the table. Applying
\eqref{eqcolumnHodgesum} gives
\eqref{eqsurfacediagonals}; rearranging gives
\eqref{eqsurfacemixedhodgenumers}.
\end{proof}

\begin{corollary}\label{corsurfacereconstruction}
The mixed Hodge numbers are completely determined by the Du Bois-Hochschild numbers via \eqref{eqsurfacemixedhodgenumers} within either of the following classes of connected normal projective complex surfaces:
\begin{enumerate}
\item surfaces for which every irreducible exceptional curve is rational;
\item surfaces with quasi-rational singularities for which every exceptional dual graph is a tree.
\end{enumerate}
Thus two Fourier--Mukai partners belonging to the same one of these classes have identical mixed Hodge numbers.
\end{corollary}

\begin{proof}

In (2), $L=0$, hence Proposition~\ref{propsurfaceformulas} gives $G=hh_{-2}-hh_2$. Rational exceptional components imply quasi-rationality by \cite[Corollary 5.4]{AH26}, hence Proposition~\ref{propsurfaceformulas} applies also in (1), where $G=0$. In both cases, substituting $G$ in \eqref{eqsurfacemixedhodgenumers} provides formulas for the mixed Hodge numbers.
\end{proof}

Note that (1) allows cycles of rational curves in the exceptional divisor; in particular, one can allow cusp singularities \cite[Example 5.5]{AH26}.

\section{Topological K-homology and the Baum--Fulton--MacPherson Todd class}
\label{sectopologicalcomparison}

The Du Bois--Hochschild numbers record column sums of Du Bois numbers, but these sums need not determine the individual Du Bois or mixed Hodge numbers. In this section, we construct additional invariants using topological K-theory and Poincar\'e duality.

The first subsection studies the canonical comparison map from the topological K-theory of $\Perf(X)$ to that of $\Coh(X)$. Using T.~Moulinos's relative topological K-theory, we identify this map with cap product by the Borel--Moore $\KU$-homology class of the structure sheaf.
Under rational Atiyah-Hirzebruch isomorphisms, it becomes the cap product with the Baum--Fulton--MacPherson Todd class. The same comparison expresses the categorical index pairing in terms of this class.

In the second subsection, we introduce the topological defect numbers as the dimensions of the kernels of the comparison map on rational homotopy groups. Its block triangular form relates these numbers to the failure of Poincar\'e duality. For projective rational homology manifolds, Poincar\'e duality forces the defect numbers to vanish. In even complex dimension, however, the symmetrized index pairing recovers the topological signature. These Fourier--Mukai invariants, which are of independent interest, supply the additional equations needed in Section~8.

\subsection{Topological K-theory and Poincar\'e duality}

We consider a proper $\bbC$-scheme $X$.

By the characterization of $\Coh(X)$ in terms of $\Perf(X)$ used in Proposition~\ref{propderivedequivalence}, the map
\begin{equation}\label{eqTheta}
\Theta_X:\K^{\mathrm{top}}(\Perf(X))\longrightarrow
\K^{\mathrm{top}}(\Coh(X))
\end{equation}
is determined by $\Perf(X)\in \Catperf_\bbC$.
Its cofiber, $\K^{\mathrm{top}}(\Sg(X))$ (see Example~\ref{Sgcategory}), is already extremely useful (see Remark~\ref{Ktopsgvstopdefect}), but the map itself provides an even more powerful invariant.

We use relative topological $\K$-theory, introduced by T.~Moulinos \cite{Mou19}. For $B\in\Sch^\sft_\bbC$, it assigns to a $\Perf(B)^\otimes$-linear category $T$ a sheaf
$\K^{\mathrm{top}}_B(T)$ of $\KU$-modules on $B^{\an}$.
For $B=\Spec\bbC$, this recovers Blanc's topological
$\K$-theory \cite[Proposition~7.6]{Mou19}, while
\[
\K^{\mathrm{top}}_B(\Perf(B))\simeq\KU_B,
\]
where $\KU_B$ denotes the constant sheaf associated to $\KU$.
The construction is lax symmetric monoidal; in particular, it carries algebra and module structures to the corresponding structures on sheaves. The following  statement describes the global sections of the relative topological K-theory over a proper base.

\begin{lemma}\label{lemrelativeKprojective}
Let $f:X\to B$ be a proper morphism in $\Sch^\sft_\bbC$. Consider $T\in \Catperf_{\Perf(X)^\otimes}$ and write $f_*T\in\Catperf_{\Perf(B)^\otimes}$ for the same category with its
$\Perf(B)$-linear structure induced by $f^*$. The natural comparison
\[
\K^{\mathrm{top}}_B(f_*T)
\longrightarrow f^\an_*\K^{\mathrm{top}}_X(T)
\]
is an equivalence compatible with products in $\K^\mathrm{top}$ and the induced algebra and module structures.

In particular, if $X\in \Sch^\sft_\bbC$ is proper, then
\[
\K^{\mathrm{top}}(T)
\xrightarrow{\sim}
\Gamma(X^{\an},\K^{\mathrm{top}}_X(T))
\]
for every $\Perf(X)$-linear category $T$, compatibly with these structures.
\end{lemma}

\begin{proof}

The underlying equivalence is recorded in \cite[Theorem~2.12]{GS22} and attributed to \cite{Mou19}.
The proof of \cite[Proposition~7.8]{Mou19}, stated there for
$T=\Perf(X)$, indeed applies for a general $T$ with minor modifications.

We use the motivic description of relative topological $\K$-theory from \cite[Section~7.3, especially Proposition 7.7]{Mou19} and that
\begin{equation}\label{KGLequi}
KGL_B(f_*T)\simeq f_*KGL_X(T).
\end{equation}
Indeed, by construction of $KGL$, the equivalence
\[
(f_*T)\otimes_{\Perf(B)}\Perf(Y)
\simeq
f_*\bigl(T\otimes_{\Perf(X)}\Perf(X\times_B Y)\bigr)
\]
for every smooth $B$-scheme $Y$
implies identification \eqref{KGLequi}. The proper-pushforward compatibility of Betti realization used in \cite[Proof of Proposition~7.8]{Mou19}
therefore gives
\[
\begin{aligned}
\K^{\mathrm{top}}_B(f_*T)
&\simeq \operatorname{Betti}_B(KGL_B(f_*T))\\
&\simeq \operatorname{Betti}_B(f_*KGL_X(T))\\
&\simeq f^\an_*\operatorname{Betti}_X(KGL_X(T))
 \simeq f^\an_*\K^{\mathrm{top}}_X(T),
\end{aligned}
\]
compatibly with products in $\K^{\mathrm{top}}$,
algebra and module structures, due to symmetric monoidality of Betti realization. The last claim follows by taking $B=\Spec\bbC$.
\end{proof}

Recall that the tensor product equips $\Coh(X)$ with a $\Perf(X)^\otimes$-module structure, providing an example of the situation above. The Hom bifunctor also gives an
index pairing on
$\pi_0\K^{\mathrm{top}}(\Perf(X))$. We denote it by
$\chi_X^{\mathrm{top}}$. The following theorem identifies the geometric meaning of the present structures.

\begin{theorem}\label{proptopcomparison}(cf. \cite[Lemma 4.2]{PT26})
\label{reqtopcomparison}
Let $X$ be a projective complex variety. Denote by $p:X\to\Spec\bbC$ the structure morphism. The comparison equivalences
\begin{equation}\label{Ktopcomparison}
\K^{\mathrm{top}}(\Perf(X))\simeq\KU(X^{\an}),
\qquad
\K^{\mathrm{top}}(\Coh(X))\simeq\KU^{\mathrm{BM}}(X^{\an})
\end{equation}
identify the action induced by $\Perf(X)^\otimes$-module structure on $\Coh(X)$ with the topological cap-product action. They carry the
class $[\mathcal O_X]\in\pi_0\K^\mathrm{top}(\Coh(X))$ to its Baum--Fulton--MacPherson class and respect proper pushforward to a point.

Consequently,
for $a\in H^{\mathrm{ev}/\mathrm{odd}}(X)$,
\begin{equation}\label{eqThetaTodd}
\Theta_{X,\bbQ}(a)=a\cap\td_*(X),
\qquad \td_*(X):=\td_*([\mathcal O_X]).
\end{equation}
For $a,b\in H^{\mathrm{ev}}(X,\bbQ)$, the index pairing is
\begin{equation}\label{eqtopologicalindex}
\chi^{\mathrm{top}}_{X,\bbQ}(a,b)
=p^\an_*\bigl((a^\vee\smile b)\cap\td_*(X)\bigr),
\qquad a^\vee_{2j}=(-1)^j a_{2j}.
\end{equation}
\end{theorem}

\eqref{eqThetaTodd} and \eqref{eqtopologicalindex} above use the Atiyah-Hirzebruch identifications
\[
\pi_0\KU(X^\an)_{\bbQ} \simeq H^{\mathrm{ev}}(X^{\an},\bbQ), \qquad
\pi_1\KU(X^\an)_{\bbQ}
\simeq H^{\mathrm{odd}}(X^{\an},\bbQ),
\]
\[
\pi_0\KU^\mathrm{BM}(X^\an)_{\bbQ}
\simeq H_{\mathrm{ev}}(X^{\an},\bbQ), \qquad
\pi_1\KU^{\mathrm{BM}}(X^\an)_{\bbQ}
\simeq H_{\mathrm{odd}}(X^{\an},\bbQ).
\]
The maps are the rational cohomological and homological Chern characters.

\begin{proof}
The equivalences \eqref{Ktopcomparison} are well known (\cite[Proposition~4.32]{Bla16}, \cite[Theorem~2.10]{HLP20}). In what follows, we essentially reinterpret \cite[Sections 4-5]{BFM79} in categorical terms.
We first verify compatibility with the module action.

Choose a closed immersion $i:X\hookrightarrow\mathbb P^N$
with complement $j:U\hookrightarrow\mathbb P^N$. We have
\[
\begin{gathered}
i^\an_*\K^{\mathrm{top}}_X(\Coh(X))
\simeq
\K^{\mathrm{top}}_{\bbP^N}(i_*\Coh(X))
\\
\simeq
\fib\bigl(
\K^{\mathrm{top}}_{\bbP^N}(\Perf(\bbP^N))
\ra
\K^{\mathrm{top}}_{\bbP^N}(j_*\Perf(U))
\bigr)
\simeq
\fib\bigl(\KU_{\bbP^N}\ra j^\an_*\KU_U\bigr),
\end{gathered}
\]
where the first equivalence is provided by Lemma~\ref{lemrelativeKprojective}, the second by localization and d\'evissage, and the last by \cite[Lemma 4.4]{PT26}.

Write $\KU^{\mathrm{BM}}_X:=(p^\an)^!\KU$
for the $\KU$-valued dualizing sheaf, whose global sections compute $\KU^{\mathrm{BM}}(X^\an)$.
The complex tangent bundles of $\mathbb P^N$ and $U$
give compatible $\KU$-orientations. Poincar\'e duality and Bott periodicity identify
\[
\KU_{\mathbb P^N}\simeq\KU^{\mathrm{BM}}_{\mathbb P^N},
\qquad
\KU_U\simeq\KU^{\mathrm{BM}}_U,
\]
as $\KU$-module sheaves, compatibly with restriction.
The localization triangle for the dualizing sheaf then gives
\[
\fib(\KU_{\mathbb P^N}\longrightarrow j^\an_*\KU_U)
\simeq i^\an_*\KU^{\mathrm{BM}}_X.
\]

All these identifications are
$\KU_{\bbP^N}$-linear. Since $i^\an_*$ is fully faithful and $(i^\an)^*\KU_{\bbP^N}\simeq\KU_X$, they give a
$\KU_X$-linear equivalence
\[
\K^{\mathrm{top}}_X(\Coh(X))
\simeq\KU^{\mathrm{BM}}_X.
\]
Applying Lemma~\ref{lemrelativeKprojective} to the map of $\Perf(X)^\otimes$-modules $\Perf(X) \ra \Coh(X)$ identifies the two $\KU(X^\an)$-module structures on $\KU^{\mathrm{BM}}(X^\an)$.

Choose a finite locally free resolution of $i_*\mathcal O_X$
on $\mathbb P^N$. After analytification, it is a complex of topological vector bundles exact on $U^\an$, and hence defines a class in $\KU^0((\mathbb P^N)^\an,U^\an)$.
Under the duality identification above, this class maps exactly to
\[
[\mathcal O_X]_{\mathrm{BFM}}
\in\KU^{\mathrm{BM}}_0(X^\an).
\]
Here we denote by $[-]_{\mathrm{BFM}}$ the Baum--Fulton--MacPherson construction of the topological $\K$-homology class of a coherent sheaf; see \cite[Sections~4--5]{BFM79}. The proof of \cite[Lemma~2.9]{HLP20} shows that the equivalence in \eqref{Ktopcomparison} is given by the same construction, is independent of the embedding, and commutes with proper pushforward. In particular, the square
\[
\begin{tikzcd}
\K^{\mathrm{top}}(\Coh(X))
  \arrow[r,"\sim"]
  \arrow[d,"{R\Gamma(X,-)}"']
&
\KU^{\mathrm{BM}}(X^\an)
  \arrow[d,"{p^\an_*}"]
\\
\K^{\mathrm{top}}(\bbC)
  \arrow[r,"\sim"]
&
\KU
\end{tikzcd}
\]
commutes.

The composite of the coherent-sheaf-to-$\K$-homology map with the homological Chern character is the Baum--Fulton--MacPherson Riemann--Roch transformation.
Consequently,
\[
\operatorname{ch}_*
\bigl([\mathcal O_X]_{\mathrm{BFM}}\bigr)
=\td_*([\mathcal O_X])=\td_*(X).
\]

Since the inclusion $\Perf(X)\to\Coh(X)$ is $\Perf(X)$-linear and sends the unit to $\mathcal O_X$, under \eqref{Ktopcomparison}, $\Theta_X$
is cap product with the Baum--Fulton--MacPherson class of $\mathcal O_X$. Compatibility of the cohomological and homological Chern characters with cap products proves \eqref{eqThetaTodd}.

Finally, applying the compatibilities above to the index pairing
\[
\Hom_{\Perf(X)}(P,Q)\simeq\Gamma(X,P^\vee\otimes Q)
\]
shows that, for
$x,y\in\pi_0\KU(X^\an)$,
\[
\chi_X^{\mathrm{top}}(x,y)
=
p^\an_*\bigl(x^\vee\cap\Theta_X(y)\bigr).
\] 
Under the Chern character, duality corresponds to the involution
$a\mapsto a^\vee$ on even cohomology defined by
\[
a^\vee=\sum_j(-1)^j a_{2j},
\qquad a=\sum_j a_{2j},\quad a_{2j}\in H^{2j}(X^\an,\bbQ).
\]
Applying \eqref{eqThetaTodd}
gives \eqref{eqtopologicalindex}.
\end{proof}

\subsection{Topological defect numbers and the signature pairing}

\begin{defi}\label{defitopdefect}
For a projective complex variety $X$, put
\begin{equation}\label{eqtopologicalkerneldimensions}
d_i(X):=\dim_{\bbQ}\ker\bigl(\pi_i\Theta_X\otimes\bbQ\bigr).
\end{equation}
We will refer to them as the \emph{topological defect numbers} (cf. Theorem~\ref{proptopdefect}).
\end{defi}

Theorem~\ref{proptopcomparison} proves the following.

\begin{corollary}
The topological defect numbers are finite, invariant under Fourier--Mukai equivalence, and depend only on the parity of $i$.
More precisely, $d_0(X)=d_{\mathrm{ev}}(X)$ and $d_1(X)=d_\mathrm{odd}(X)$ are the dimensions of the kernels of
\[
-\cap\td_*(X):H^{\mathrm{ev}}(X,\bbQ)\to H_{\mathrm{ev}}(X,\bbQ),
\qquad
-\cap\td_*(X):H^{\mathrm{odd}}(X,\bbQ)\to H_{\mathrm{odd}}(X,\bbQ),
\]
respectively.
\end{corollary}

For an irreducible $n$-dimensional $X$, the top component of
$\td_*(X)$ is $[X]$.
These maps are block triangular, with diagonal blocks
\begin{equation}\label{capblocks}
-\cap[X]:H^j(X,\bbQ)\longrightarrow H_{2n-j}(X,\bbQ).
\end{equation}
The two outer blocks, corresponding to $j=0$ and $j=2n$, are invertible. We will use this description to compute the kernel dimensions for surfaces and threefolds.

We note that considering cokernels instead of kernels in Definition~\ref{defitopdefect} gives nothing new.

\begin{obs}\label{kercoker}
Let $X$ be an irreducible projective complex variety. Then, for every $i$,
\[
\dim(\ker(\pi_i\Theta_X\otimes\bbQ))=\dim(\coker(\pi_i\Theta_X\otimes\bbQ))
\]
\end{obs}

\begin{proof}
By Theorem~\ref{proptopcomparison}, the source and target of $\pi_{i}\Theta_X\otimes\bbQ$ identify with $H^{\mathrm{ev/odd}}(X,\bbQ)$ and $H_{\mathrm{ev/odd}}(X,\bbQ)$ for the same parity, respectively. The source and the target for a given $i$ have the same finite dimension, hence the kernel and cokernel have the same dimension
$d_i(X)$.    
\end{proof}

\medskip

The symmetrized index pairing gives another handy invariant. 

Recall (\cite{BM83}, \cite[Section~3]{PP25a}) that an equidimensional complex variety $X$ of dimension $n$ is a rational homology manifold precisely when the canonical comparison
\[
\bbQ_X[n]\longrightarrow IC_X
\]
is an equivalence.
If $X$ is proper, this identifies ordinary cohomology with intersection cohomology, compatibly with duality.
In particular, cap product with the fundamental class gives isomorphisms
\[
-\cap[X]:H^j(X,\bbQ)\xrightarrow{\sim}H_{2n-j}(X,\bbQ).
\]
The corresponding perfect pairing is the ordinary cup-product pairing
\[
H^j(X,\bbQ)\otimes H^{2n-j}(X,\bbQ)
\longrightarrow\bbQ,
\qquad
a\otimes b\longmapsto\int_X a\smile b.
\]

The above provides all that is necessary to define the topological signature as a derived invariant as in \cite{Ant26}. 
The following is an adaptation of the argument given there for smooth projective varieties to singular projective rational homology manifolds.

\begin{prop}\label{propcategoricalsignature}
Let $X$ be a connected projective complex rational homology manifold of complex dimension $2d$. If $S_X(a,b):=\chi_X^{\mathrm{top}}(a,b)
+\chi_X^{\mathrm{top}}(b,a)$ is the symmetrized index pairing on
$\pi_0\K^{\mathrm{top}}(\Perf(X))\otimes_{\bbZ}\mathbb R$, then
\[
\sig(S_X)=(-1)^d\sig_{\mathrm{top}}(X),
\]
where $\sig_{\mathrm{top}}(X)$ is the signature of the intersection pairing on $H^{2d}(X,\mathbb R)$. In particular, both signatures are invariant under Fourier--Mukai equivalence within this class of varieties.
\end{prop}

We will refer to $\sig_{\mathrm{top}}(X)$ as the \emph{topological signature} of $X$.

\begin{proof}
By the Poincar\'e duality recalled above, there is a unique class
$t_X\in H^{\mathrm{ev}}(X,\bbQ)$ such that
\[
\td_*(X)=t_X\cap[X].
\]
Its degree-zero component is $1$, since the top-dimensional component of $\td_*(X)$ is $[X]$.

Via the Chern character, identify
\[
\pi_0\K^{\mathrm{top}}(\Perf(X))\otimes_{\bbZ}\mathbb R
\simeq H^{\mathrm{ev}}(X,\mathbb R),
\]
and write $\chi_{X,\mathbb R}^{\mathrm{top}}$ for the scalar extension of the index pairing in these coordinates.
By Theorem~\ref{proptopcomparison},
\[
\chi_{X,\mathbb R}^{\mathrm{top}}(a,b)
=p_*\bigl((a^\vee\smile b)\cap\td_*(X)\bigr)=\int_X a^\vee\smile b\smile t_X=:\int_X a^\vee b\,t_X.
\]

Since $\dim_{\bbC}X=2d$, the involution $(-)^\vee$ acts trivially on $H^{4d}$, hence
\[
\chi_{X,\mathbb R}^{\mathrm{top}}(b,a)=\int_X b^\vee a\,t_X=\int_X(b^\vee a\,t_X)^\vee=\int_X a^\vee b\,t_X^\vee,
\]
and
\[
S_X(a,b)=\int_Xa^\vee b\,(t_X+t_X^\vee).
\]

The class
\[
v:=\frac{t_X+t_X^\vee}{2}-1
\]
has positive degree and is therefore nilpotent, hence one can define
\[
B=((t_X+t_X^\vee)/2)^{1/2} = (1+v)^{1/2} =\sum_{r\geq0}\binom{1/2}{r}v^r,
\]
since the binomial expansion here has only a finite number of nonzero terms. This also shows that $B$ is invertible; moreover, $B^\vee=B$.

Hence we have
\[
S_X(a,b)=2\int_X(Ba)^\vee(Bb).
\]
Since $B$ is invertible, \[\sig(S_X)=\sig\bigl((a,b)\mapsto\int_Xa^\vee b\bigr).\]

The summands $H^{2j}(X,\mathbb R)$ and
$H^{4d-2j}(X,\mathbb R)$ for $j<d$ are isotropic and pair perfectly with one another. Hence these summands do not contribute to the signature and only the middle-dimensional pairing remains.
Since $a^\vee=(-1)^da$ for $a\in H^{2d}(X)$, we have
\[
\sig(S_X)=(-1)^d\sig_{\mathrm{top}}(X).
\]
Fourier-Mukai equivalences preserve the index pairing and hence its symmetrization, proving the final claim.
\end{proof}

\medskip

\section{Hodge numbers of Fourier--Mukai partners: surfaces, threefolds, and fourfolds}
\label{sechodgereconstruction}

We combine the Du Bois--Hochschild numbers with the invariants of Section~7 coming from topological K-theory to recover mixed Hodge numbers for several classes of projective varieties. In each case, geometry additionally restricts possible Hodge numbers, so that we arrive at a simple system of linear equations.

For connected normal surfaces, the topological defect numbers provide exactly the additional information needed to recover all mixed Hodge numbers. For connected threefolds with rational singularities, the even topological defect equals the
$\bbQ$-factoriality defect. Together with the Du Bois--Hochschild numbers and the irregularity, it again determines all mixed Hodge numbers.
For connected rational homology fourfolds, cohomology is pure, and the topological signature provides the additional equation needed for reconstruction once $h^{0,1}$ and $h^{0,2}$ are fixed. In each case, the reconstruction formulas we obtain imply the corresponding Fourier--Mukai invariance statements. We also discuss some particular classes of varieties to which they apply.

The final subsection considers the limitations of what can be recovered from the underlying localizing motive alone. All the Du Bois invariants factor through this motive, but we give examples showing that individual Hodge numbers do not. This illustrates the importance of retaining categorical information beyond the localizing motive in the Hodge numbers invariance question.

\medskip

Unless otherwise specified, all varieties in this section are projective over $\bbC$, and cohomology has rational coefficients.

We use $(-)^*$ for the $\bbQ$-linear dual vector space.

We write $h_n^{(a,b)}=h_n^{a,b}=h_n^{b,a}$ for unordered pairs as in Proposition~\ref{propsurfaceformulas}.

For background on intersection cohomology, its Hodge structures, Poincar\'e duality, and the Hodge--Lefschetz theorems used below, we refer to \cite[Sections~1.9 and~2.1]{dCM09}.

\subsection{Normal surfaces}
\label{secnormalsurfacereconstruction}

Having introduced the additional invariants $d_0,d_1$, we can significantly extend Corollary~\ref{corsurfacereconstruction}.

We start by relating these invariants to mixed Hodge numbers of a surface.

\begin{prop}\label{propsurfacetopologicaldefects}
Let $X$ be a connected normal projective complex surface. Then
\begin{equation}\label{eqsurfacetopologicaldefects}
\begin{aligned}
d_0(X)&=\dim_{\bbQ}W_1H^2(X)
       =h_2^{0,0}(X)+2h_2^{(0,1)}(X),\\
d_1(X)&=b_3(X)-b_1(X)
       =2h_3^{(1,2)}(X)-2h_1^{(0,1)}(X).
\end{aligned}
\end{equation}
\end{prop}

Note that, since $X$ is connected and normal, it is irreducible.

\begin{proof}
Take a projective resolution $\pi:\widetilde X\to X$, an isomorphism away from singularities, with a reduced simple normal crossings exceptional divisor $E$. As in the proof of Proposition~\ref{propsurfaceformulas}, we use the long exact sequence of mixed Hodge structures \eqref{lesMHSsurface} coming from the corresponding abstract blowup square; the difference is that the restriction
$H^1(\widetilde X)\to H^1(E)$ now does not necessarily vanish. Once again, negative definiteness of the exceptional intersection matrix makes
$H^2(\widetilde X)\to H^2(E)$ surjective. Consequently
\begin{equation}\label{eqsurfacecohomology}
\begin{gathered}
H^1(X)\hookrightarrow H^1(\widetilde X),\qquad
H^3(X)\xrightarrow{\sim}H^3(\widetilde X),\qquad
H^4(X)\simeq\bbQ(-2),\\
\ker\bigl(H^2(X)\to H^2(\widetilde X)\bigr)
\simeq\operatorname{coker}\bigl(H^1(\widetilde X)\to H^1(E)\bigr).
\end{gathered}
\end{equation}
By strictness, the kernel here contains $W_1H^2(X)$, and it is exactly $W_1H^2(X)$, since the cokernel has nothing in weight two. We see also that $H^1(X)$ and $H^3(X)$ are pure of weights 1 and 3, respectively.

Now we analyze the diagonal blocks in \eqref{capblocks}.
The projection formula factors the cap product as
\[
H^j(X)\xrightarrow{\pi^*}H^j(\widetilde X)
\xrightarrow{-\cap[\widetilde X]}H_{4-j}(\widetilde X)
\xrightarrow{\pi_*}H_{4-j}(X).
\]
Here the middle map is the Poincar\'e duality isomorphism and the last map for $j_0$ is dual to the first map for $4-j_0$. Hence, by \eqref{eqsurfacecohomology}, the cap product is injective for $j=1$ and surjective for $j=3$. The triangular form of
\eqref{eqThetaTodd} therefore implies
\[
\begin{aligned}
d_1
&=\dim_\bbQ\ker\bigl(
H^{\mathrm{odd}}(X)\xrightarrow{-\cap[X]}H_{\mathrm{odd}}(X)
\bigr)
=\dim_\bbQ\ker\bigl(
H^3(X)\xrightarrow{-\cap[X]}H_1(X)
\bigr)\\
&=b_3(X)-b_1(X)
=2h_3^{(1,2)}(X)-2h_1^{(0,1)}(X).
\end{aligned}
\]
where the last equality follows from the purity of $H^{1/3}$ proven above. It remains to understand the kernel of the cap product on $H^\mathrm{ev}$.

The outer blocks in \eqref{capblocks}, corresponding to $j=0,4$, are invertible. We claim that
\[
\ker\bigl(H^2(X)\xrightarrow{-\cap[X]}H_2(X)\bigr)
=W_1H^2(X).
\]
Indeed, the long exact sequence gives 
\[
\operatorname{im}\bigl(\pi^*:H^2(X)\to H^2(\widetilde X)\bigr)
=
\ker\bigl(H^2(\widetilde X)\to H^2(E)\bigr)
=
\langle[E_j]\rangle^\perp,
\]
where $[E_j]\in H^2(\widetilde X)$ are the classes of the irreducible exceptional divisors, and orthogonality is taken with respect to the intersection pairing on $\widetilde X$. The negative definite exceptional intersection matrix makes the restriction of this pairing to
$\langle[E_j]\rangle$ nondegenerate.
Poincar\'e duality therefore gives an orthogonal decomposition
\[
H^2(\widetilde X)
=
\langle[E_j]\rangle
\oplus
\langle[E_j]\rangle^\perp,
\]
and the pairing on the second summand is nondegenerate as well.

For $\alpha,\beta\in H^2(X)$, the projection formula gives
\[
\int_X\alpha\smile\beta
=
\int_{\widetilde X}\pi^*\alpha\smile\pi^*\beta.
\]
If $\alpha\cap[X]=0$, the left-hand side vanishes for every $\beta$.
Thus $\pi^*\alpha$ is orthogonal to
$\operatorname{im}\pi^*=\langle[E_j]\rangle^\perp$.
Since $\pi^*\alpha$ itself belongs to this subspace, nondegeneracy of its intersection pairing forces
$\pi^*\alpha=0$. Hence
\[
\ker\bigl(H^2(X)\xrightarrow{-\cap[X]}H_2(X)\bigr) = \ker\pi^* = W_1H^2(X).
\]
Consequently,
\[
d_0=\dim_\bbQ\ker(H^{\mathrm{ev}}(X)\xrightarrow{(-)\cap[X]} H_{\mathrm{ev}}(X))=\dim_\bbQ(W_1H^2(X))=h_2^{0,0}(X)+2h_2^{(0,1)}(X).
\]
\end{proof}

Adding the obtained equations
to the ones describing Du Bois-Hochschild numbers gives a system of linear equations with exactly one solution, proving the following.

\begin{theorem}\label{thmnormalsurfacereconstruction}
Let $X$ be a connected normal projective complex surface. Its mixed Hodge numbers are determined by $d_0(X),d_1(X)$ and $hh_i(X)$:
\begin{equation}\label{eqnormalsurfaceinversion}
\begin{aligned}
h_0^{0,0}(X)=h_4^{2,2}(X)&=1,\\
h_1^{(0,1)}(X)&=\frac{hh_1}{2}-\frac{d_1(X)}{4},\\
h_3^{(1,2)}(X)&=\frac{hh_1}{2}+\frac{d_1(X)}{4},\\
h_2^{(0,2)}(X)&=hh_2,\\
h_2^{(0,1)}(X)&=d_0(X)-hh_{-2}+hh_2,\\
h_2^{0,0}(X)&=2(hh_{-2}-hh_2)-d_0(X),\\
h_2^{1,1}(X)&=hh_0-2-d_0(X)+hh_{-2}-hh_2.
\end{aligned}
\end{equation}
All other mixed Hodge numbers vanish. In particular, connected normal projective surfaces which are Fourier--Mukai partners have identical mixed Hodge numbers.
\end{theorem}

\begin{proof}
Combining Proposition~\ref{propsurfacetopologicaldefects} with Corollary~\ref{DBHodgeformula} gives the formulas. Their right-hand sides are defined in terms of $\Perf(X)$, hence are shared by Fourier-Mukai partners.
\end{proof}

\begin{remark}
For a surface with quasi-rational singularities, Proposition~\ref{propsurfaceformulas} gives $d_1(X)=0$ (in fact, this is provided by Lemma~\ref{qratvanishing}), and
\[
d_0(X)=L+2G.
\]
Together with $hh_{-2}-hh_2=L+G$ (Proposition~\ref{propsurfaceformulas}), this recovers
\[
G=d_0(X)-hh_{-2}+hh_2,\qquad
L=2(hh_{-2}-hh_2)-d_0(X).
\]
Thus for surfaces with quasi-rational singularities, $hh_i$ together with $d_0$ recover all the mixed Hodge numbers. Under the additional assumptions of Corollary~\ref{corsurfacereconstruction}, either $G=0$ or $L=0$, forcing $d_0$ to equal $hh_{-2}-hh_2$ or $2(hh_{-2}-hh_2)$, respectively.
\end{remark}

\medskip

\subsection{Threefolds with rational singularities}

For a connected projective complex threefold $X$ with rational singularities, write
\[
q=h^1(X,\mathcal O_X),\qquad
\sigma(X)=\dim_{\bbQ}(\Cl(X)/\Pic(X))_{\bbQ},
\]
for the irregularity and the $\bbQ$-factoriality defect (see \cite{Kaw88,PP25b}), respectively.
We first give a categorical description of the defect, motivated by \cite{KPS21}.

\begin{theorem}\label{proptopdefect}
Let $X$ be a connected projective complex threefold with rational singularities. Then
\[
d_0(X)=\sigma(X).
\]
In particular, $\sigma(X)$ is invariant under Fourier--Mukai equivalence within this class.
\end{theorem}

Note that, since $X$ here is connected and has rational singularities, it is normal and hence, in fact, irreducible.

\begin{proof}

The proof of \cite[Theorem~E, Section~7]{PP25b} shows that
\[
H^2(X)\hookrightarrow IH^2(X),
\qquad
H^4(X)\xrightarrow{\sim}IH^4(X),
\]
for the usual comparison maps and that
\[
b_4(X)-b_2(X)=\sigma(X)
\]
(for isolated rational singularities and $H^2(X,\mathcal O_X)=0$, this was proved in \cite[Theorem~3.2]{NS95}; Park--Popa's result applies in the greater generality considered here). We use these comparisons to analyze the diagonal cap-product blocks in \eqref{capblocks}.
 
The cap product with $[X]$ can be described as the composite
\begin{equation}\label{capviaIH}
H^j(X)\longrightarrow IH^j(X)
\xrightarrow{\sim}IH^{6-j}(X)^*
\longrightarrow H^{6-j}(X)^*
\simeq H_{6-j}(X),
\end{equation}
using intersection Poincar\'e duality and the compatibility of the comparison maps with Verdier duality.

For $j=2$, the first arrow is injective and the other arrows are isomorphisms. For $j=4$, the third arrow is surjective and the other arrows are isomorphisms. Thus $-\cap[X]$ is injective on $H^2$ and surjective on $H^4$, with a kernel of dimension $b_4(X)-b_2(X)=\sigma(X)$. Finally, $-\cap[X]$ is an isomorphism on $H^0$ and $H^6$, hence
\[
d_0(X)=\dim\ker\bigl(H^4(X)\to H_2(X)\bigr)
=b_4(X)-b_2(X)=\sigma(X).
\]\end{proof}
By \cite[Theorem~A]{PP25a}, applied to rational singularities,
\begin{equation}\label{ParkPopasymmetries}
\underline h^{0,j}(X)=\underline h^{j,0}(X)
=\underline h^{3,3-j}(X)=\underline h^{3-j,3}(X),
\qquad 0\leq j\leq3.
\end{equation}
In particular,
$\underline h^{0,0}=\underline h^{3,3}=1$ and
\[
\underline h^{3,2}=\underline h^{2,3}=\underline h^{1,0}=\underline h^{0,1}=h^1(X,\mathcal O_X)=q.
\]
Here the second to the last equality holds, since rational singularities are Du Bois.
Combining these equalities with the column-sum formula
\eqref{eqcolumnHodgesum} gives
\[
\begin{aligned}
hh_3=hh_{-3}&=\underline h^{3,0},&
hh_2=hh_{-2}&=2\underline h^{2,0},\\
hh_1&=2q+\underline h^{2,1},&
hh_{-1}&=2q+\underline h^{1,2},\\
hh_0&=2+\underline h^{1,1}+\underline h^{2,2}.
\end{aligned}
\]
Finally, \cite[Theorem~E]{PP25b} gives
\[
\sigma(X)=\underline h^{2,2}(X)-\underline h^{1,1}(X).
\]
Thus the Hodge--Du Bois numbers are recovered by
\begin{equation}\label{eqthreefoldinversion}
\begin{aligned}
\underline h^{3,0}&=hh_3,&
\underline h^{2,0}&=\frac{hh_2}{2},\\
\underline h^{2,1}&=hh_1-2q,&
\underline h^{1,2}&=hh_{-1}-2q,\\
\underline h^{1,1}&=\frac{hh_0-2-\sigma(X)}{2},&
\underline h^{2,2}&=\frac{hh_0-2+\sigma(X)}{2},
\end{aligned}
\end{equation}
together with $\underline h^{1,0}=q$, $\underline h^{0,0}=1$, and the symmetries \eqref{ParkPopasymmetries}. It remains to recover the weight grading, splitting the Du Bois numbers into sums of mixed Hodge numbers.

We use the exact sequence of mixed Hodge structures \cite[(7.2)]{PP25b}:
\[
0\longrightarrow H^2(X)\longrightarrow IH^2(X)
\longrightarrow V^H\longrightarrow H^3(X)
\longrightarrow IH^3(X)\longrightarrow0.
\]
Here $V^H$ is pure of weight two and $H^i(X)\simeq IH^i(X)$ for $i\ne2,3$.
It follows that $H^i(X)$ is pure of weight $i$ unless $i=3$, while
$H^3(X)$ has only weights 2 and 3, with
$\Gr_3^WH^3(X)\simeq IH^3(X)$.
The equality $\underline h^{0,3}=\underline h^{3,0}$, together with Hodge symmetry on $\Gr_3^WH^3(X)$, forces the $(0,2)$ part of
$\Gr_2^WH^3(X)$ to vanish, hence its $(2,0)$ part vanishes as well by Hodge symmetry. Thus
\begin{equation}\label{eqthreefoldweights}
\begin{aligned}
h_3^{1,1}(X)&=hh_{-1}-hh_1,\\
h_3^{(0,3)}(X)&=hh_3,\qquad
h_3^{(1,2)}(X)=hh_1-2q.
\end{aligned}
\end{equation}
All other mixed Hodge numbers of $H^3(X)$ vanish. In the remaining degrees, purity identifies $h_{p+q}^{p,q}(X)$ with
$\underline h^{p,q}(X)$.

By Theorem~\ref{FMcolumns} and Theorem~\ref{proptopdefect}, the defect and the Du Bois-Hochschild numbers are derived invariants. The equations \eqref{eqthreefoldinversion} and \eqref{eqthreefoldweights} therefore prove the following.

\begin{theorem}\label{thmthreefoldreconstruction}
Let $X,Y$ be connected projective complex threefolds with rational singularities which are Fourier--Mukai partners. If
\[
h^1(X,\mathcal O_X)=h^1(Y,\mathcal O_Y),
\]
then $h_n^{a,b}(X)=h_n^{a,b}(Y)$ for every triple $(a,b,n)$.
\end{theorem}

In particular, the theorem applies when both irregularities vanish.
This already includes many examples of interest. For instance, Fano threefolds with klt singularities have rational singularities and vanishing irregularity by
\cite[Theorems~5.22 and~2.70]{KM98}.
Complete intersection threefolds with rational singularities provide another class of examples: the vanishing of
$H^1(X,\mathcal O_X)$ follows from the Koszul resolution.
Calabi--Yau threefolds with rational singularities and
$H^1(X,\mathcal O_X)=0$ provide another interesting class;
these assumptions are satisfied, for example, by the three-dimensional Calabi--Yau factors in the Beauville--Bogomolov decomposition in the singular projective setting \cite[Definition~1.4 and Theorem~1.5]{HP19}.

The irregularity condition can be checked on an arbitrary resolution $\widetilde X\to X$: indeed, by rationality of the singularities,
\begin{equation}\label{rationalsingirr}
H^1(X,\mathcal O_X)
\simeq H^1(\widetilde X,\mathcal O_{\widetilde X})
\end{equation}
Therefore threefolds with rational singularities admitting rationally connected smooth projective resolutions provide another broad class of examples (see, for example, \cite[Proposition~2.5(b)]{Deb24}). This includes rational, stably rational, and unirational threefolds with rational singularities
\cite[(1)]{Deb24}.
It also includes all normal projective toric threefolds:
they are rational and have rational singularities
\cite[Theorem~11.4.2]{CLS11}.

\begin{remark}
Theorem~\ref{thmthreefoldreconstruction} applies to several of the singular Fourier--Mukai partners mentioned in Subsection~2.1. For the expected-dimensional Pfaffian--Grassmannian threefolds of \cite[Theorem~11.4]{Kuz06} (see Example~\ref{examplePfaffian}) and the intersections of two Grassmannians of \cite[Theorem~1.1]{BCP20}, the irregularities vanish. Thus, whenever both partners have rational singularities, all their mixed Hodge numbers agree. We do not know of a classical geometric argument explaining these equalities of mixed Hodge numbers in the stated generality.

The theorem also applies to projective flops of threefolds with terminal Gorenstein singularities \cite[Theorem~1.1]{Che02}, without assuming that the irregularity vanishes. Indeed, terminal singularities are rational, so the irregularities can be computed on smooth projective resolutions, where they agree by birational invariance. When both sides are $\bbQ$-factorial, Park--Popa \cite[Corollary~8.1]{PP25b} prove that their Du Bois numbers agree. As shown above, the only possible contribution below the expected weight is the $(1,1)$-part of $H^3$, of dimension
$\underline h^{1,2}-\underline h^{2,1}$. Thus the result of Park--Popa gives another proof that the mixed Hodge numbers agree in this case.
\end{remark}

\medskip

\begin{remark}\label{Ktopsgvstopdefect}
    Our first proof of Theorem~\ref{thmthreefoldreconstruction} used the topological K-theory of the singularity category in the spirit of \cite{KPS21, PS21}. We finish this subsection by a statement confirming that, in the setting of Theorem~\ref{thmthreefoldreconstruction}, $d_0$ and $\dim_\bbQ\K^{\mathrm{top}}_{-1}(\Sg(X))_\bbQ$ add the same amount of information to the information provided by the Du Bois-Hochschild numbers.
\end{remark} 

\begin{prop}
Let $X$ be a connected projective complex threefold with rational singularities. Then
\[\dim_{\bbQ}\K^{\mathrm{top}}_{-1}(\Sg(X))_{\bbQ} =\sigma(X)+d_1(X)\]
and
\[
d_1(X)=hh_{-1}(X)-hh_1(X).
\]
\end{prop}

\begin{proof}
Applying $\K^{\mathrm{top}}$ to Example~\ref{Sgcategory} gives a short exact sequence
\[
0\longrightarrow
\coker(\pi_{-1}\Theta_X\otimes\bbQ)
\longrightarrow
\K^{\mathrm{top}}_{-1}(\Sg(X))_\bbQ
\longrightarrow
\ker(\pi_{-2}\Theta_X\otimes\bbQ)
\longrightarrow0.
\]
By Observation~\ref{kercoker}, the cokernel on the left has dimension
$d_1(X)$. Since $d_{-2}(X)=d_0(X)=\sigma(X)$ by Theorem~\ref{proptopdefect}, this proves the first equality.

For the second equality, consider the odd-degree blocks
\eqref{capblocks}. The blocks
\[
H^1(X)\longrightarrow H_5(X),
\qquad
H^5(X)\longrightarrow H_1(X)
\]
are isomorphisms, since ordinary and intersection cohomology agree in these degrees. Thus,
\[
d_1(X)=
\dim_\bbQ\ker\bigl(H^3(X)\xrightarrow{-\cap[X]}H_3(X)\bigr).
\]
The middle cap-product map factors as
\[
H^3(X)\longrightarrow IH^3(X)
\xrightarrow{\sim}IH^3(X)^*
\longrightarrow H^3(X)^*=H_3(X);
\]
this is \eqref{capviaIH}.
The first map is surjective with kernel $W_2H^3(X)$, as established above, and the last map is its dual, hence injective. Consequently, by \eqref{eqthreefoldweights},
\[
d_1(X)=\dim_\bbQ W_2H^3(X)
=hh_{-1}(X)-hh_1(X).
\]
\end{proof}

\subsection{Rational homology fourfolds}

The results of this section were inspired by B.~Antieau's blog post \cite{Ant26} proving the Hodge numbers invariance conjecture for smooth projective complex fourfolds. Aside from the derived invariance of Hochschild homology allowing one to recover column sums of Hodge numbers, the proof essentially rests on three ingredients:
the derived invariance of the topological signature of even-dimensional varieties, the Hodge index theorem (see e.g. \cite[Theorem 6.33]{Voi02}) relating it to Hodge numbers, and derived invariance of $h^{0,1}$ and $h^{0,2}$ of fourfolds (due to \cite{PS11} and \cite{Abu17}, respectively). 
As in the lower-dimensional cases, Corollary~\ref{DBHodgeformula} provides a counterpart to the Hochschild column-sum formula in the singular setting, suggesting that one should look for a class of fourfolds to which the argument of \cite{Ant26} extends.
We observe that the first two ingredients mentioned above extend to projective rational homology manifolds.

Let $X$ be a connected projective complex rational homology manifold of dimension four. As recalled in Subsection~7.2, the rational homology manifold assumption identifies ordinary cohomology with intersection cohomology
\[
H^k(X,\bbQ)\xrightarrow{\sim}IH^k(X,\bbQ),
\]
compatibly with duality. This is already enough to get the first ingredient -- derived invariance of topological signature in the form of Proposition~\ref{propcategoricalsignature}.
The identification also respects Hodge structures and the action of an ample class. The intersection cohomology of a projective variety is pure and satisfies hard Lefschetz and the Hodge--Riemann bilinear relations (see
\cite[\S1.9, items~1--2, and Theorem~5.2.1]{dCM09} or \cite[Theorem~2.2.3]{dCM05}), hence the same holds for the ordinary cohomology. This is enough to get the second ingredient.

We write
\begin{equation}\label{puritynotation}
h^{p,q}(X):=h_{p+q}^{p,q}(X)=\underline h^{p,q}(X).
\end{equation}
All mixed Hodge numbers outside weight $k$ in $H^k(X)$ vanish. Hodge symmetry and Poincar\'e duality give
\[
h^{p,q}(X)=h^{q,p}(X)=h^{4-p,4-q}(X).
\]

\begin{obs}[The Hodge index theorem for rational homology evenfolds]
Let $X$ be a connected projective complex rational homology manifold of complex dimension $2d$ and $\sig_{\mathrm{top}}(X)$ the signature of the intersection pairing on $H^{2d}(X, \mathbb{R})$ (as in Proposition~\ref{propcategoricalsignature}). Then 
\[
\sig_{\mathrm{top}}(X)=\sum_{p,q}(-1)^q h^{p,q}(X).
\]
\end{obs}

\begin{proof}
We refer to the proof of the classical Hodge index theorem \cite[Theorem~6.33]{Voi02}. The proof only uses the Lefschetz decomposition and Hodge--Riemann relations \cite[Theorem~6.32]{Voi02}, which are both available for intersection cohomology (\cite[Theorem~2.2.3(c)--(d)]{dCM05}).
As recalled above, ordinary cohomology agrees here with intersection cohomology, compatibly with Hodge structures, the ample-class action, and the intersection pairing. Thus Voisin's proof goes through.
\end{proof}
For a fourfold, by Hodge symmetry and Poincar\'e duality the signature formula simplifies to
\begin{equation}\label{eqfourfoldsignature}
\sig_{\mathrm{top}}(X)=\sum_{p,q}(-1)^q h^{p,q}(X)=2+4h^{0,2}+2h^{0,4}-2h^{1,1}-2h^{1,3}+h^{2,2}.
\end{equation}

Now we have everything we need to prove the following.

\begin{theorem}\label{thmfourfoldreconstruction}
Let $X$ be a connected projective complex rational homology fourfold. Then $H^k(X)$ is pure of weight $k$ for every $k$, and, in notation \eqref{puritynotation},
\begin{equation}\label{eqfourfoldinversion}
\begin{aligned}
h^{0,0}(X)&=1,\\
h^{0,3}(X)&=\frac{hh_3(X)}{2},\\
h^{0,4}(X)&=hh_4(X),\\
h^{1,2}(X)&=\frac{hh_1(X)}{2}-h^{0,1}(X),\\
h^{1,3}(X)&=hh_2(X)-2h^{0,2}(X),\\
h^{1,1}(X)
&=\frac{
hh_0(X)-2hh_2(X)+2hh_4(X)
+8h^{0,2}(X)-\sig_{\mathrm{top}}(X)
}{4},\\
h^{2,2}(X)
&=\frac{
hh_0(X)+2hh_2(X)-2hh_4(X)
-8h^{0,2}(X)+\sig_{\mathrm{top}}(X)-4
}{2}.
\end{aligned}
\end{equation}
Together with $h^{0,1}(X)$, $h^{0,2}(X)$, and the symmetries
\[
h^{p,q}(X)=h^{q,p}(X)=h^{4-p,4-q}(X),
\]
these formulas determine all mixed Hodge numbers of $X$.

In particular, if $X$ and $Y$ are Fourier--Mukai partners within this class and
\[
h^{0,1}(X)=h^{0,1}(Y),
\qquad
h^{0,2}(X)=h^{0,2}(Y),
\]
then their mixed Hodge numbers agree.
\end{theorem}

Note that Poincar\'e duality makes all the diagonal blocks \eqref{capblocks} invertible, so the topological defect numbers $d_i(X)$ vanish. It also implies that $X$ is actually irreducible: the dimensions of the source and target of
\[
-\cap[X]:H^0(X,\bbQ)\xrightarrow{\sim}H_8(X,\bbQ)
\]
equal the numbers of connected components and irreducible components of $X$, respectively.

\begin{proof}
Purity was established above. Combining \eqref{eqcolumnHodgesum}, Hodge symmetry, and Poincar\'e duality gives
\[
\begin{aligned}
hh_4&=h^{0,4},&
hh_3&=2h^{0,3},\\
hh_2&=2h^{0,2}+h^{1,3},&
hh_1&=2h^{0,1}+2h^{1,2},\\
hh_0&=2+2h^{1,1}+h^{2,2}.
\end{aligned}
\]
Adding \eqref{eqfourfoldsignature} gives a system of linear equations which, for fixed $h^{0,1}$ and $h^{0,2}$, has exactly one solution.

The Du Bois--Hochschild numbers are invariant under Fourier--Mukai equivalence by Theorem~\ref{FMcolumns}, and the topological signature is invariant by Proposition~\ref{propcategoricalsignature}.
The final claim follows.
\end{proof}

Recall that an irreducible complex fourfold is a rational homology manifold iff the link of every point has the rational homology of $S^7$. Thus this is an analytically local condition, which holds, in particular, for finite quotient singularities and for isolated ADE hypersurface singularities(\cite{Web96}). More general criteria for varieties to be rational homology manifolds can be found in \cite{PP25a}.

There are several interesting classes to which Theorem~\ref{thmfourfoldreconstruction} applies with
\[
h^{0,1}(X)=h^{0,2}(X)=0.
\]
We describe three such classes below. All varieties considered are connected and projective over $\bbC$. Their singularities are rational, hence Du Bois, so
\[
h^{0,i}(X)
=\underline h^{0,i}(X)
=\dim_\bbC H^i(X,\mathcal O_X).
\]

\begin{enumerate}
\item
\emph{Quasismooth weighted complete-intersection fourfolds.}
Let
\[
X=X_{d_1,\ldots,d_c}
\subset\bbP(a_0,\ldots,a_{4+c})
\]
be a well-formed quasismooth weighted complete intersection fourfold (we refer to \cite{BK16} for the terminology used here).
Such an $X$ has only finite quotient singularities by \cite[Theorem~3.1.6]{Dol82}, and is therefore a rational homology manifold. Moreover,
\cite[Theorem~3.2.4(iii)]{Dol82} gives
\[
H^i(X,\mathcal O_X)=0,\qquad 0<i<4.
\]
In particular, $h^{0,1}$ and $h^{0,2}$ vanish.

This class includes Fano, Calabi--Yau, and general-type fourfolds. For more details and extensive collections of hypersurface ($c=1$) examples in these three subclasses, see \cite{BK16}.

\item
\emph{Complete-intersection fourfolds with isolated ADE singularities.}
Let $X\subset\bbP^{4+c}$ be a complete intersection with only isolated ADE hypersurface singularities.
The Koszul resolution gives
\[
H^1(X,\mathcal O_X)=H^2(X,\mathcal O_X)=0.
\]
Thus the theorem applies to this class as well.

\item
\emph{Fano fourfolds with quotient singularities.}
Let $X$ be a Fano fourfold with only finite quotient singularities. These singularities make $X$ a rational homology manifold and are klt. Kawamata--Viehweg vanishing, applied to $\mathcal O_X$, therefore gives
\[
H^i(X,\mathcal O_X)=0,\qquad i>0;
\]
see \cite[Theorem~2.70]{KM98}.
Consequently the theorem applies here without any complete-intersection assumption.
\end{enumerate}

In particular, Fourier--Mukai partners belonging to any of these classes (possibly different ones) have identical mixed Hodge numbers.

\subsection{Localizing motives and categorical invariants}

The Morita-invariance statements of Sections 5-6
have the following stronger version.

\begin{obs}
The Du~Bois Hochschild and cyclic homology theories are determined by the localizing motive
$\Uloc(X):=\Uloc(\Perf(X))$.
Over $\bbC$, the same holds for both Deligne cohomology theories of Section~4.
In particular, for a proper $k$-variety $X$, the numbers
\[
hh_i(X)=\dim_k\HH_{\DB,i}(X/k)
=\sum_{p-q=i}\underline h^{p,q}(X)
=\sum_{\substack{n,a,b\\2a-n=i}}h_n^{a,b}(X)
\]
depend only on $\Uloc(X)$.
\end{obs}

\begin{proof}
These invariants factor through $\Uloc$ by their constructions. The numerical equalities follow from Corollary~\ref{DBHodgeformula}.
\end{proof}

One distinction between the last two sections and the rest of the paper is that the invariants introduced in Section~7
and used in Section~8 are, a priori, ``categorical'' rather than ``motivic'': even though, for a complex projective variety $X$,
\[
\K^{\mathrm{top}}(\Coh(X))\simeq
\underline{\operatorname{Hom}}_{KU}
\bigl(\K^{\mathrm{top}}(\Perf(X)),KU\bigr)
\]
is determined by the localizing motive of $X$, there is no clear reason why the $\Perf$--$\Coh$ comparison map or the index pairing would be motivic.
The discussion below, especially Example~\ref{examplemotivesforgetHodge}, illustrates why retaining information beyond the underlying localizing motive is important for the Hodge numbers invariance question.

\medskip

The functor $\Uloc$ is known to lose some important information.

\begin{example}[Beilinson, Orlov]
Let $k$ be a field, and let $X/k$ be obtained from a product of projective spaces by successive blowups along
products of projective spaces
.
Beilinson's exceptional collections \cite{Bei78} and Orlov's blowup formula
\cite[Theorem~4.3 and Corollary~4.4]{Orl93}
give a full exceptional collection on $X$.
 Consequently, writing
$\chi=\operatorname{rk}_{\bbZ}K_0(X)$, we have
\[
\Uloc(X)\simeq\mathbbm1_k^{\oplus \chi}
\simeq
\Uloc(\coprod_{i=1}^{\chi}\Spec k).
\]
In particular, the localizing motive determines neither dimension nor the number of connected components.
\end{example}

The example above shows that even under the irreducibility assumption on $X$, the localizing motive loses substantial information.
We finish the paper by using the construction above to demonstrate two more manifestations of information loss by $\Uloc$ in the geometric setting.

\begin{example}
By Anel--To\"en's theorem \cite[Corollary~7.1]{AT09}, a smooth projective complex variety has at most countably many smooth projective Fourier--Mukai partners. In fact, they prove that in a certain precise sense, the functor $\Perf$ restricted to smooth projective varieties has zero-dimensional fibers.

In contrast, a fixed localizing motive can be shared by positive-dimensional families of smooth projective varieties. For $n\ge5$, unordered sets of $n$ distinct points in $\bbP^2_\bbC$, considered up to projective automorphism, have moduli of dimension $2n-8$, describing the associated blowups together with their blowdown to $\bbP^2$; forgetting the blowdown identifies at most countably many members for any fixed isomorphism class of surfaces.
Nevertheless, all these surfaces have localizing motive
\[
\mathbbm1_\bbC^{\oplus(n+3)}.
\]
\end{example}

For a fixed $n$, however, all the blowups above have the same Hodge numbers (in fact, the same Chow motive).
The following example shows that equality of Hodge numbers does not follow from equivalence of localizing motives. 

\begin{example}\label{examplemotivesforgetHodge}
The smooth projective fivefolds
$$ \bbP^1_k\times\bbP^2_k\times\bbP^2_k, \qquad \operatorname{Bl}_{\Spec k\sqcup\Spec k\sqcup \Spec k}\bbP_k^5 $$
both admit full exceptional collections of length $18$. Thus they have the same localizing motive $\mathbbm 1_k^{\oplus 18}$ and all the additive invariants agree on these varieties. Nevertheless, when the base field is of characteristic zero, their $h^{1,1}$ are 3 and 4, respectively.
\end{example}

This demonstrates that individual Hodge numbers cannot be recovered from the underlying localizing motive alone, even for smooth projective varieties.

\end{document}